\documentclass{conm-p-l}

\usepackage{amssymb}
\usepackage{amsfonts}
\usepackage{amscd}
\usepackage[mathscr]{eucal}
\usepackage{verbatim}
\usepackage{epsfig}

\input xy
\xyoption{all}

\newtheorem{thm}{Theorem}[subsection]
\newtheorem{lem}[thm]{Lemma}
\newtheorem{prop}[thm]{Proposition}
\newtheorem{cor}[thm]{Corollary}
\newtheorem{rem}[thm]{Remark}

\newtheorem{sth}{Theorem}[thm]

\newtheorem{scor}[sth]{Corollary}

\theoremstyle{definition}

 \newtheorem{exams}[thm]{Examples}

 \newtheorem{ack}{Acknowledgments}  

\numberwithin{equation}{thm}

\newcommand{\Cref}[1]{Corollary~\textup{\ref{#1}}}
\newcommand{\Dref}[1]{Definition~\textup{\ref{#1}}}

\newcommand{\Lref}[1]{Lemma~\textup{\ref{#1}}}
\newcommand{\Pref}[1]{Proposition~\textup{\ref{#1}}}
\newcommand{\Rref}[1]{Remark~\textup{\ref{#1}}}
\newcommand{\Sref}[1]{Section~\textup{\ref{#1}}}
\newcommand{\Ssref}[1]{Subsection~\textup{\ref{#1}}}
\newcommand{\Tref}[1]{Theorem~\textup{\ref{#1}}}

\def\bilap#1{\hbox to 0pt{\hss#1\hss}}
 \def\Rarrow#1{\bilap{\hbox to#1{\rightarrowfill}}}
 \def\Larrow#1{\bilap{\hbox to#1{\leftarrowfill}}}
 \def\Equals#1{\bilap
                  {\raise 4pt\hbox
                    {\vrule width#1 height.5pt}%
                   \kern-#1\raise 1pt\hbox
                    {\vrule width#1 height.5pt}%
                  }}

\newcommand{\EQAL}[1]%
{\,\begin{picture}(#1,0)%
\put(0,3){\line(1,0){#1}}%
\put(0,1){\line(1,0){#1}}%
\end{picture}\,}%

\newcommand{\vlto}[1]%
{\,\begin{picture}(#1,3)%
\put(0,2){\vector(1,0){#1}}%
\end{picture}\,}%

\newcommand{\vllarrow}[1]%
{\,\begin{picture}(#1,3)%
\put(#1,2){\vector(-1,0){#1}}%
\end{picture}\,}%

\newcommand{\dirlm}[1]%
  {
     {\lim\hskip-1.58em\lower.65ex
       \hbox{$
                {}_{\stackrel{\lower1ex\hbox
                                        {$\scriptstyle -\!\!\!\longrightarrow$}
                                      }{\vbox to0pt{\vss\vskip.6ex
                                            \hbox{$\scriptstyle{}^{#1}$}\vss}}
                   }
            $}
     }
\:}

\newcommand{\subdirlm}[1]%
  {
     {\lim\hskip-1.5em\lower.6ex
       \hbox{$
                   {}_{\stackrel{\lower1ex\hbox
                                           {$\scriptstyle\longrightarrow$}
                                }{ ^{#1} }
                      }
             $}
     }
\:}

\newcommand{\inlm}[1]%
   {
      {\lim\hskip-1.58em\lower.65ex
        \hbox{$
                 {}_{\stackrel{\lower1ex\hbox
                                        {$\scriptstyle \longleftarrow\!\!\<-$}
                              }{\vbox to0pt{\vss\vskip.6ex
                                            \hbox{$\scriptstyle{}^{#1}$}\vss}}
                    }
             $}
      }
\:}

\def\hz#1{{\hbox to 0pt{#1}}}
\def\Iso{\vbox to 0pt{\vss\hbox{$\widetilde{\phantom{nn}}$}\vskip-7pt}}

\def\>{\mspace {1mu}}
\def\<{\mspace{-1mu}}
\def\({{\textup(}}
\def\){{\textup)}}

\newcommand{\fm}{{\mathfrak{m}}}
\newcommand{\fn}{{\mathfrak{n}}}
\newcommand{\fp}{{\mathfrak{p}}}

\newcommand{\X}{{\mathscr X}}

\newcommand{\Y}{{\mathscr Y}}
\newcommand{\Z}{{\mathscr Z}}
\newcommand{\V}{{\mathscr V}}
\newcommand{\U}{{\mathscr U}}
\newcommand{\W}{{\mathscr W}}
\newcommand{\I}{{\mathscr I}}
\newcommand{\J}{{\mathscr J}}

\newcommand{\Spec}{{\mathrm {Spec}}}
\newcommand{\Spf}{{\mathrm {Spf}}}
\newcommand{\cm}{{\mathrm {CM}}}
\newcommand{\res}[1]{{\mathrm {res}}_{#1}}
\newcommand{\wres}[1]{\widetilde{\mathrm {res}}_{#1}}

\newcommand{\A}{{\mathcal A}}
\newcommand{\B}{{\mathcal B}}
\newcommand{\C}{{\mathcal C}}
\newcommand{\cD}{{\mathcal D}}
\newcommand{\De}{{\Delta}}

\newcommand{\E}{{\mathcal E}}
\newcommand{\F}{{\mathcal F}}
\newcommand{\G}{{\mathcal G}}

\newcommand{\Hr}{{\mathrm H}}
\newcommand{\Hp}[2]{\Hr_{#1}^{\<\prime\>#2}}
\newcommand{\Rr}{{\mathrm R}}

\newcommand{\cH}{{\mathcal H}}

\newcommand{\cO}{{\mathcal O}}

\newcommand{\cR}{{\mathcal R}}

\newcommand{\D}{{\mathbf D}}

\newcommand{\bbP}{{\mathbb P}}
\newcommand{\bbG}{{\mathbb G}}
\newcommand{\lbbG}{\overline{\mathbb G}}
\newcommand{\bbF}{{\mathbb F}}
\newcommand{\bbFc}{{\bbF_c}}
\newcommand{\cbbF}{{{\bbF}^{\,*}}}
\newcommand{\rbbF}{{{\bbF}^{\,r}}}
\newcommand{\cbbFc}{{\bbF^{\,*}_c}}
\newcommand{\rbbFc}{{\bbF^{\,r}_c}}
\newcommand{\Coz}{{\mathrm {Coz}}}
\newcommand{\Cozs}[1]{{\mathrm {Coz}}_{#1}}

\newcommand{\Ed}[1]{E_{#1}}
\newcommand{\fE}{f^E}

\newcommand{\Dqc}{\D_{\mkern-1.5mu\mathrm {qc}}}
\newcommand{\wDqc}{ \widetilde
         {\vbox to6.5pt{\vss\hbox{$\mathbf D$}}}
   _{\mkern-1.5mu\mathrm {qc}} }
\newcommand{\wDc}{ \widetilde
         {\vbox to6.5pt{\vss\hbox{$\mathbf D$}}}
   _{\mkern-1.5mu\mathrm {c}} }
\newcommand{\wDqcp}{\wDqc^{\lower.5ex\hbox{$\scriptstyle+$}}}
\newcommand{\wDcp}{\wDc^{*}}

\newcommand{\Dqct}{\D_{\mkern-1.5mu\mathrm{qct}}}

\newcommand{\Dc}{\D_{\mkern-1.5mu\mathrm c}}

\newcommand{\R}{{\mathbf R}}

\newcommand{\bL}{{\mathbf L}}
\newcommand{\Hom}{{\mathrm {Hom}}}

\newcommand{\Homb}{{\mathrm {Hom}}^{\bullet}}

\newcommand{\Aqct}{\A_{\mathrm {qct}}\<}

\newcommand{\fs}{f^!}
\newcommand{\fc}{{\boldsymbol{F}}}
\newcommand{\gc}{{\boldsymbol{G}}}
\newcommand{\hc}{{\boldsymbol{H}}}
\newcommand{\fca}{{\boldsymbol{F_1}}}
\newcommand{\fcb}{{\boldsymbol{F_2}}}
\newcommand{\fcc}{{\boldsymbol{F_3}}}
\newcommand{\gca}{{\boldsymbol{G_1}}}
\newcommand{\gcb}{{\boldsymbol{G_2}}}
\newcommand{\eca}{{\boldsymbol{E_1}}}
\newcommand{\ecb}{{\boldsymbol{E_2}}}
\newcommand{\bmu}{{\boldsymbol{\mu}}}
\newcommand{\bnu}[1]{{\boldsymbol{\nu_{#1}}}}
\newcommand{\fb}{{\boldsymbol{f}}}
\newcommand{\gb}{{\boldsymbol{g}}}
\newcommand{\hb}{{\boldsymbol{h}}}

\newcommand{\gs}{g^!}
\newcommand{\ga}[1]{\gamma^!_{\hbox{$\scriptstyle{#1}$}}}

\newcommand{\gae}[1]{\gamma^E_{\hbox{$\scriptstyle{#1}$}}}
\newcommand{\gas}[1]{\lambda_{\hbox{$\scriptstyle{#1}$}}}
\newcommand{\gab}{{\boldsymbol{\gamma}}^{\boldsymbol{!}}}
\newcommand{\sh}[1]{{#1}^{\sharp}}
\newcommand{\shr}[1]{{#1}^{\boldsymbol{!}}}
\newcommand{\sha}[1]{{#1}^{\boldsymbol{\sharp}}}
\newcommand{\shrbr}[1]{{#1}^{(!)}}
\newcommand{\shrbrb}[1]{{#1}^{\boldsymbol{(!)}}}
\newcommand{\nat}[1]{{#1}^{\natural}}
\newcommand{\natb}[1]{{#1}^{\boldsymbol{\natural}}}

\newcommand{\Tr}[1]{{\mathrm {Tr}}_{{#1}}}
\newcommand{\vTr}[1]{{\mathrm {Tr}}^*_{{#1}}}
\newcommand{\ttr}[1]{{\mathrm {\tau}}_{#1}}

\newcommand{\sHom}{\cH{om}}
\newcommand{\sHomb}{\cH{om}^{\bullet}}
\newcommand{\iGp}[1]{{\varGamma_{\<\!#1}'}}
\newcommand{\iG}[1]{{\varGamma_{\<\!#1}^{\phantom\prime}}}

\newcommand{\set}{\!:=}
\newcommand{\lra}{\longrightarrow}
\newcommand{\iso}%
{{\mkern8mu\longrightarrow \mkern-25.5mu{}^\sim\mkern17mu}}
\newcommand{\osi}%
{{\mkern8mu\longleftarrow \mkern-24.5mu{}^\sim\mkern16mu}}
\newcommand{\Otimes}{\underset
  {\vbox to 0pt {\vskip-1ex\hbox{$\scriptscriptstyle=$}\vss}}
    {\otimes}\vadjust{\kern.4pt}}

\newcommand{\smcirc}%
  {{\raise.15ex\hbox to.7em{$\hss \scriptstyle\circ\hss$}}} 

\newcommand{\hd}{\pmb{-}}

\title[Duality for Cousin complexes]%
{Duality for Cousin Complexes}

\author{Pramathanath Sastry}

\thanks{The author would like to thank the Mathematisches Forschungsinstitut 
at Oberwolfach, Germany and the Banff International Research Station at
Banff, Alberta, Canada for their help at a critical period of this research.
This work was also funded by the Ganita Lab at the University of Toronto.}

\date{\today}

\begin{document}
{\bf\hfill \today}

\begin{abstract}{We relate the variance theory for Cousin complexes
$\sha{\boldsymbol{-}}$ developed by Lipman, Nayak and the author
to Grothendieck duality for Cousin complexes. Specifically
for a Cousin complex $\F$ on $(\Y,\,\De)$---with $\De$ a codimension
function on a formal scheme $\Y$ (noetherian, universally catenary)---and
a pseudo-finite type map $f\colon (\X,\,\De') \to (\Y,\,\De)$ of such
pairs of schemes with codimension functions, we show there is a derived
category map $\ga{f}(\F)\colon \sh{f}\F \to \fs\F$, functorial
in $\F \in \Cozs{\De}(\Y)$, inducing a functorial isomorphism
$\sh{f}\F \simeq E(\sh{f}\F) \iso E(\fs\F)$
(where $E$ is the Cousin functor on $(\X,\,\De')$). The map $\ga{f}(\F)$
is itself an isomorphism if (and clearly only if) $\fs\F$ is
Cohen-Macaulay on $(\X,\,\De')$---which will be so, for example,
whenever the complex $\F$ is injective or whenever 
the map $f$ is flat. For a fixed Cousin
complex  $\F$ on $(\Y,\,\De)$, $\ga{f}(\F)$ is an isomorphism for every 
map $f$ with target $(\Y,\,\De)$ if and only if $\F$ is a complex
of (appropriate) injectives. For a fixed map $f$, the functorial map $\ga{f}$ 
is an isomorphism of functors if and only if $f$ is flat.

We also generalize the Residue Theorem of Grothendieck for
residual complexes to Cousin complexes by defining a functorial 
{\em Trace Map} of graded
$\cO_\Y$--modules $\Tr{f}(\F)\colon f_*\sh{f}\F \to \F$ (a sum of local
residues) such that when $f$ is pseudo-proper, $\Tr{f}(\F)$ is a map
of complexes and the pair $(\sh{f}\F,\,\Tr{f}(\F))$ represents the functor
$\Hom(f_*\G,\,\F)$ of Cousin complexes $\G$ on $(\X,\,\De')$.
}
\end{abstract}

%

\maketitle

\setcounter{tocdepth}{2}
\tableofcontents
\newpage

\section{\bf Introduction} This paper integrates the variance theory
of Cousin complexes---worked out by Lipman, Nayak and the author in 
\cite{lns}---with the variance theory of the twisted inverse
image (``upper shriek") occurring in Grothendieck duality for
noetherian formal schemes in the form obtained by Alonso, Jerem{\'\i}as
and Lipman in \cite{dfs} (with a very important input from
Nayak \cite{suresh}).

It is useful for this introduction to use the catchall symbols
$\sha{\boldsymbol{-}}$ and $\shr{\boldsymbol{-}}$ to denote in
one stroke the entire variance theory for Cousin complexes in
\cite{lns} and the variance theory of ``upper shriek" respectively.
First, let us point to two results in this paper which can be stated
entirely in terms of Grothendieck duality, i.e., entirely in the
framework of $\shr{\boldsymbol{-}}$. To fix ideas we will (for
this paragraph) restrict ourselves to ordinary noetherian
excellent schemes admitting codimension functions\footnote{i.e. an
integer valued function $\De$ on the given scheme $Y$ such that
$\De(y')= \De(y)+1$ if $y'$ is an immediate specialization of $y$
(cf.~\cite[\S\S\,2.1]{lns}).}.\,The first one, viz.~\Tref{thm:flat},
says that if $f\colon X\to Y$ is a separated finite type map and $\De$ is
a codimension function on $Y$, then the twisted inverse image
functor $\fs$ takes Cohen-Macaulay complexes (with respect to
$\De$) to Cohen-Macaulay complexes (with respect to $\sh{f}\De$, where
$\sh{f}\De$ is as in \cite[Example\,2.1.2]{lns}) if and only if $f$ is flat. 
The second result (viz.~ \Tref{thm:gorenstein}) concerns Gorenstein 
complexes on $(Y,\De)$, i.e., complexes $\F$ which are Cohen-Macaulay with 
respect to $\De$ and such that the associated Cousin complex $E_{\De}\F$ is
a complex of injective quasi-coherent $\cO_Y$--modules (cf.~\cite[p.\,248]{RD}
where Hartshorne unnecessarily restricts himself to bounded complexes).
We show that $\F$ is Gorenstein with respect to $\De$ if and only
if $\fs\F$ is Cohen-Macaulay with respect to $\sh{f}\De$ for
every separated finite type map $f$ whose target is $Y$.

The above results are not our main results, but they are the ones that
could be stated without referring to the constructions in \cite{lns}
and therefore provide a kinder gentler introduction for the lay reader
to our work.  However, the proofs of the just stated results (which are valid
for formal schemes) require an understanding of the interrelationships
between $\sha{\boldsymbol{-}}$ and $\shr{\boldsymbol{-}}$ (the main theme
of this paper).

The constructions of $\sha{\boldsymbol{-}}$ and $\shr{\boldsymbol{-}}$
are in spirit and outlook antithetical. The approach to Grothendieck
duality in \cite{dfs} and \cite{suresh} is global, top down and 
holistic---in the spirit of Deligne and Verdier (cf.~\cite{deligne},
\cite{SGA4} and \cite{f!}). One begins by providing $\shr{\boldsymbol{-}}$
for pseudo-proper maps by global methods; and then works one's way downward
(via the flat base change theorem \cite[Theorem\,7.4]{dfs})
to pseudo-finite type maps which are composites of
compactifiable maps by showing that $\shr{\boldsymbol{-}}$ is local. 
This is not straightforward---the stumbling block being the fact that in the
category of formal schemes a closed subscheme of an open subscheme of a 
scheme $\X$ need not be an open subscheme of a closed subscheme of $\X$. 
The ``localness" of $\shr{\boldsymbol{-}}$ is proved by
Nayak \cite{suresh} via a surprising twist of Deligne's
argument for ordinary schemes.

At the other extreme is the construction of the variance theory 
$\sha{\boldsymbol{-}}$. The outlook from the outset is punctual (complete
local rings!). And reductionist in the following sense; the initial search 
is for the basic irreducible units (the atoms) of such a theory at the
level of formal neighborhoods of points and sheaves supported 
there---i.e. at the level of complete local rings and zero-dimensional modules 
(cf.~\cite{I-C}).  In such an approach it has to be an article of faith 
that these local (punctual) constructions somehow link up and 
give a global canonical variance theory. This faith is not misplaced and 
in \cite{lns} the core result is that given data of the form 
$$
D=(\X\xrightarrow{f}\Y,\,\De,\F)
$$ 
where $f$ is a morphism in the category $\bbF$ of \cite{lns} (which is
included in the category of formal schemes with codimension functions
on them, and essentially pseudo-finite type maps), $\De$ is a codimension
function on $\Y$ and $\F$ is a Cousin complex on $(\Y,\,\De)$, there is a
Cousin complex $\sh{f}\F$ on $(\X,\,\sh{f}\De)$ which is functorial in
complexes $\F$ with a fixed codimension function $\De$. Here $\sh{f}\De$
is the codimension function on $\X$ induced by $f$ as in \cite[2.1.2]{lns}.
The basic units for this construction are modules on complete
local rings $(R,\,\fm_R)$ such that $\Gamma_{\fm_R}M=M$ (i.e. $M$ is
a zero-dimensional $R$--module). In \cite{I-C}, I-Chiau Huang worked out
a variance theory ${\boldsymbol{-}}_{\boldsymbol{\#}}$ for such
objects with respect to local homomorphisms $(R,\fm_R)\xrightarrow{\varphi}
(S,\fm_S)$ of complete local rings such that the residue field extension
$k_R\to k_S$ is finitely generated. For $D$ as above and a point $x\in\X$
with $f(x)=y$, one has a map $\varphi^x\colon \widehat{\cO}_{\Y,y} \to
\widehat{\cO}_{\X,x}$ and the functor $\sh{f}$ is so constructed that---among
other properties that it enjoys---it satisfies the relation $(\sh{f}\F)(x)
=\varphi^x_{\#}\F(f(x))$. The construction is done in such a way
that $\sha{\boldsymbol{-}}$ is a variance theory, i.e. it is a 
pseudofunctor. 

What then do these theories have to do with each other?\footnote{I am
aware that this might be come under the heading of ``setting up
a straw man" but this was indeed my original confusion.} For data $D$ with
$f$ a composite of compactifiable maps (so that $f^!$ exists)
there seems no obvious way to compare $\sh{f}\F$ and $\fs\F$. And
the hope (admittedly faint) that the two might be abstractly
isomorphic in $\D(\X)$ is dashed by the counter-example
obtained by setting (with $k$ a field) $\X=\Spec{(k)}$, $\Y=\Spec{(k[T])}$,
$f$ the map given by the natural map $(T\mapsto 0) \colon k[T] \to k$, 
and $\F=f_*{k^{\sim}}$. The natural codimension function $\De$ on $\Y$ for
the Cousin complex $\F$ is the one which gives the closed
point of $\Y$ value zero. In this instance $\sh{f}\F = k^{\sim}$
and $\fs\F = k^{\sim}\oplus k^{\sim}[1]$. In fact $\fs\F$ is not
even Cohen-Macaulay (with respect to the codimension function $\De'=0$)
and so is far from being isomorphic to $\sh{f}\F$. However, note that
$E(\fs\F)$ is indeed isomorphic to $\sh{f}\F$, where $E$ is the Cousin
functor associated with $\De'(=\sh{f}\De)$.

We show that this is true in general for data $D$. 
We also investigate when
$\sh{f}\F$ and $\fs\F$ are isomorphic. The bridge between the
two theories is a derived category comparison map $\sh{f}\F \to \fs\F$ 
which on applying the (appropriate) Cousin functor $E$ transforms 
to an isomorphism. The key
to obtaining this map is the observation that if the point $x$
is closed in its fiber, then with $M=\F(y)$, $R={\widehat{\cO}}_{\Y,y}$
and $S={\widehat{\cO}_{\X,x}}$ there is an $R$--linear
punctual trace map $\Tr{f,x}(\F)\colon \varphi^x_{\#}M \to M$ which
induces an isomorphism $\Phi\colon \varphi^x_{\#}M \iso \Hom^c_R(S,\,M)$
\cite[Chapter\,7]{I-C}.
This allows one to have a map of graded $\cO_\Y$--modules
\begin{equation*}\label{eq:(TR)}\tag{TR}
\Tr{f}(\F)\colon f_*\sh{f}\F \to \F.
\end{equation*}

Our principal results are the following.

1) (The Trace Theorem) If $f$ in $D$ is pseudo-proper then 
the {\em trace map} $\Tr{f}(\F)$ is a map of complexes
(cf.~\Tref{thm:tr-thm}(a)). The crucial step here is the 
proof of the trace theorem for (relative) projective space over
ordinary schemes. This takes up a considerable amount of time, and
involves digressions into residues and the residue theorem for 
projective space. Once the trace theorem is established for
projective space, then the proof given in \cite{RD} for
residual complexes applies without change to our situation.
The punctual trace $\Tr{f,x}$ has a transitivity property
which now gives the following result (cf.~\Tref{thm:tr-thm}(b)): 
let $g:\W\to \X$ be a second pseudo-proper map. Then identifying 
$\sh{(fg)}\F$ with $\sh{g}\sh{f}\F$ (a part of the variance theory for 
$\sha{\boldsymbol{-}}$) we have
$$
\Tr{fg}\F = \Tr{f}(\F)\circ f_*\Tr{g}(\sh{f}\F).
$$

2) First suppose $f$ as above is pseudo-proper. 
Let $(f^!,\,\ttr{f})$ be the
dualizing pair given by Grothendieck duality. By the universal 
property of $(f^!,\,\ttr{f})$, the map $\Tr{f}(\F)$ gives rise
to a map 
$$
\ga{f}(\F)\colon \sh{f}\F \to \fs\F
$$ 
such that $\ttr{f}(\F)\circ\R f_*(\ga{f})=\Tr{f}(\F)$ (cf.~\eqref{eq:compare}). 
We are implicitly identifying $\R{f_*}\F$ with $f_*\F$ since $\F$, being a 
Cousin complex, is flasque.
Next suppose $f$ compactifiable. Then $\ga{f}$ can be defined as
$v^*\ga{\bar{f}}$ where $v$ is an open immersion, ${\bar f}$ a
pseudo-proper map and $f={\bar f}v$. We show that this is independent
of the compactification $(v,{\bar f})$ of $f$. We actually show more.
Suppose $f$ is a composite of compactifiable maps. Then applying the
just described process repeatedly one can define $\ga{f}$. We show 
in \Tref{thm:compare} that
this is independent of the factorization of $f$ as a composite of
compactifiables, and that this comparison map respects pseudofunctoriality
(in a sense made precise by \Tref{thm:compare}(b)). We point out
that in the ``classical" situation (i.e., ordinary schemes \dots), a separated
finite type map is always compactifiable by a theorem of Nagata
\cite{nagata}. (See also \cite{nagata2} and \cite{nagata3}.)

3) Suppose $f$ in $D$ is a composite of compactifiable maps. Set
$f^E\F = E(\fs\F)$ where $E=E_{\sh{f}\De}$ is the Cousin
associated with the codimension function $\sh{f}\De$ on $\X$.
Applying $E$ to $\ga{f}$ we obtain a map of Cousin complexes on
$(X,\,\sh{f}\De)$,
$$
\gae{f}(\F) \colon \sh{f}\F \to f^E\F.
$$
\Tref{thm:heart'} states that $\gae{f}(\F)$ is an isomorphism. 
The proof is not straightforward. The crucial step is \Tref{thm:heart}
which states that if $f$ is {\em smooth} then the map $\ga{f}(\F)$ is an 
isomorphism (note that the last is a statement about $\ga{f}$ which is
stronger than the corresponding statement about $\gae{f}$). 
The proof here is involved and in the end
amounts understanding Grothendieck duality for smooth maps
as well as the compatibility of local and global duality in terms
of endomorphisms of residual complexes. This result on smooth maps
takes up all of \Sref{s:smooth}. Now a general $f$ of the type we are
considering can be locally factored as a 
closed immersion followed by a smooth map, whence
we are reduced to showing $\gae{f}(\F)$ is an isomorphism when
$f$ is a closed immersion. This is easily proven (cf.~\Cref{cor:closed-imm}).

4) A consequence of the circle of ideas above is \Tref{thm:CM}
which states that $\ga{f}(\F)$
is an isomorphism if (and clearly only if) $\fs\F$ is Cohen-Macaulay
with respect to $\sh{f}\De$. We use ideas from Suominen \cite{suominen}
who shows that (with a fixed codimension function) the category of 
Cohen-Macaulay complexes is equivalent to the category of Cousin complexes.
Now $\fs\F$ being Cohen-Macaulay can be viewed as a condition on
$\F$ (if $f$ is allowed to vary) or as a condition on $f$ (if $\F$ is
allowed to vary within its codimension class). 
For fixed $\F$ we show that $\ga{f}(\F)$ is an isomorphism for every composite 
of compactifiable maps $f$ with target $(\Y,\,\De)$ if and only if 
$\F$ is a complex of $\Aqct(\Y)$--injectives.
Here $\Aqct(\Y)$ is the category of quasi-coherent $\cO_\Y$ modules
$\F$ which are {\it torsion}, i.e., satisfying
$\iGp{\Y}\F=\F$ where $\iGp{\Y}$ is as in the last paragraph of
\cite[1.2.1]{dfs}.  In particular $\ga{f}(\F)$
is an isomorphism when $\F$ is residual. 
At the other end, for a fixed map $f\colon (\X,\,\De')\to (\Y,\,\De)$, 
we show that $\ga{f}(\F)$ is an isomorphism
of complexes---for all Cousin complexes $\F$ with respect to $\De$---if and
only if $f$ is {\em flat}. (Cf. \Tref{thm:injective} and \Tref{thm:flat}.)
We have already stated these results (at the beginning of this introduction)
entirely in the framework of $\shr{\boldsymbol{-}}$.

5) We show (cf.~\Tref{thm:u-prop}) that given data of the form 
$((\X,\,\De')\xrightarrow{f}(\Y,\,\De),\,\F)$ where $f$ is a 
pseudo-proper map and $\F$ is a Cousin complex on
$(\Y,\,\De)$, the pair $(\sh{f}\F,\,\Tr{f}(\F))$ represents the functor
$\Hom_\Y(f_*\G,\,\F)$ of Cousin complexes on $(\X,\,\De')$, i.e., the pair
$(\sh{f}\F,\,\Tr{f}(\F))$ induces an isomorphism
$$
\Hom_{\De'}(\G,\,\sh{f}\F) \iso \Hom_\Y(f_*\G,\,\F).
$$
(Here the left side is the group of morphisms in $\Cozs{\De'}$
between $\G$ and $\sh{f}\F$.) In other words the variance theory 
$\sha{\boldsymbol{-}}$, together with the trace maps \eqref{eq:(TR)} 
for pseudo-proper maps $f$, is a true duality theory for Cousin complexes.

One consequence is this. Suppose $f\colon (\X,\,\De')\to (\Y,\,\De)$
is a {\em finite} map of schemes and for $\F\in \Cozs{\De}(\Y)$,
$f^\flat\F \in \Cozs{\De'}(\X)$ is the unique complex of quasi-coherent
$\cO_\X$--modules satisfying $f_*f^\flat\F = \sHomb_{\Y}(f_*\cO_\X,\,\F)$,
then a certain obvious isomorphism of the graded $\cO_\X$ modules
$\sh{f}\F$ and $f^\flat\F$ is an isomorphism of complexes 
(cf.~\Cref{cor:finite}).

6) In certain circumstances we can define the twisted inverse image
of $f$ even if $f$ is not a composite of open immersions and pseudo-proper
maps. In \Sref{s:(!)} we imitate the theory given in \cite{RD} 
to obtain this when
$\Y$ in the data $D$ has a bounded residual complex on it. For this
we apply the variance theory $\sha{\boldsymbol{-}}$ of \cite{lns}
to residual complexes. In this case we define $\fs\E$ via residual
complexes for objects 
$\E$ in $\D(\Y)$ such that $\R\iGp{\Y}\E$ is in $\D^+(\Y)$ and such
that there exists a complex $\G$ with {\em coherent} homology sheaves 
satisfying $\R\iGp{\Y}\E\simeq \R\iGp{\Y}\G$. We also show 
that this construction is indeed $\fs\E$ when $f$ is a composite of
compactifiable maps (cf.~\Tref{thm:(!)!}).

\subsection{Conventions}\label{ss:conventions}
In addition to the notations and conventions in
\cite[1.4]{lns} we use the following in this paper.
Minor differences in notation between this paper and
\cite{lns} are also noted (see item \ref{item:diff} below).

\begin{enumerate}
\item\label{item:diff} We use the upright symbol $\Hr^i$ for
cohomology rather than the slanted symbol $H^i$ used in \cite{lns}.
Further if $\X$ is a scheme and $\I$ an $\cO_\X$ ideal then
we use $\Gamma_\I$ for the functor $\dirlm{n}\Hom_{\cO_\X}(\cO_\X/\I^n,-)$
and the symbol $\iG{\I}$ for the sheafified version
$\dirlm{n}\sHom_{\cO_\X}(\cO_\X/\I^n,-)$. In \cite{lns} the
slanted version $\iG{}$ is the one used throughout (the context determining
the interpretation to be placed).

\item\label{item:wdqc} For a locally noetherian formal
scheme $\X$,
$$
\wDqc(\X)\set \R\iGp{\X}^{-1}(\Dqc(\X))
$$
is the triangulated subcategory of $\D(\X)$ whose objects are those complexes
$\F$ such that $\R\iGp{\X}\F \in \Dqc(\X)$---or equivalently, 
$\R\iGp{\X}\F \in \Dqct(\X)$. 
\item For $\X$ as above,
$\Dc^*(\X)$ is the essential image of $\D_{\mathrm c}(\X)$ under
$\R\iGp{\X}$, i.e., $\Dc^*(\X)$ is the full subcategory
of $\D(\X)$ such that $\E\in \Dc^*(\X) \Leftrightarrow \E\simeq
\R\iGp{\X}\F$ with $\F\in \Dc(\X)$.
\item For $\X$ as above,
$$
\wDcp(\X) \set \R\iGp{\X}^{-1}(\Dc^*(\X)\cap \D^+(\X))
$$
is the triangulated subcategory of $\D(\X)$ whose objects are complexes
$\F$ such that $\R\iGp{\X}\F \in \Dc^*(\X)\cap \D^+(\X)$.
\item If $f:\X\to \Y$ is a morphism in $\bbF$ which is pseudo-proper,
then we denote the resulting trace map by 
$\ttr{f}:\R f_*f^! \to {\boldsymbol 1}_{\D(\Y)}$.
(See \cite[Theorem\,6.1 (a)]{dfs}).
\item The category $\bbG$ has as objects noetherian formal
schemes, and its morphisms are maps of noetherian formal
schemes which are composites of compactifiable maps,
or equivalently, composites of open immersions and
pseudo-proper maps.
\item The category $\cbbF$ has the same objects as $\bbF$
but its morphisms are composites of compactifiable maps, or
equivalently, composites of open immersions and pseudo-proper
maps, i.e., $\cbbF=\bbF\cap \bbG$. Note that if $f\colon \X\to \Y$
is a map in $\bbG$ and $\Y\in \bbF$ then $\X\in \bbF$
(cf.~\cite[2.1.2]{lns}) whence $f$ is a map in
$\cbbF$. The category $\cbbFc$ has 
the same objects as $\bbFc$, and its morphisms are maps in $\bbFc$ 
such that the underlying map of formal schemes is in $\cbbF$.
\item The category $\rbbF$ is the full subcategory of $\bbF$
consisting of objects admitting a {\em bounded} residual complex
\cite[\S\S\,9.1]{lns}. Note that if $\X\xrightarrow{f}\Y$ is a
map in $\bbF$ and $\Y$ is an object in $\rbbF$ then $\X$ is an
object in $\rbbF$ (cf.~\cite[Prop.\,9.1.4]{lns}). The category
$\rbbFc$ is the full subcategory of $\bbFc$ consisting of
objects $(\X,\,\De)$ such that $\X\in \rbbF$. Note that
if $(\X,\,\De')\xrightarrow{f} (\Y,\,\De)$ is a map in $\bbFc$
and $(\Y,\,\De) \in \rbbFc$, then $(\X,\,\De')\in \rbbFc$.
\item For $(\X,\,\De)\in \bbFc$, $\cm(\X,\,\De)$ is the full
subcategory of $\Dqct^+(\X)$ consisting of complexes which
are Cohen-Macaulay with respect to $\De$, i.e. if $\D_{{\mathrm{CM}}}^+(\X;\De)$
is as in \cite[\S\S\,3.3]{lns} then $\cm(\X,\,\De)\set 
\D^+_{{\mathrm{CM}}}(\X;\De) \cap \Dqct(\X)$.  
The category $\cm^*(\X,\,\De)$ is the full subcategory
of $\cm(\X,\,\De)$ given by $\cm^*(\X,\,\De) = \cm(\X,\,\De)\cap \Dc^*(\X)$.
\item If $\natb{\boldsymbol{-}}$ is a contravariant
pseudofunctor \cite[\S\,4]{lns}
on a category ${\mathfrak{C}}$ then, for a pair of maps 
$X\xrightarrow{f}Y\xrightarrow{g}Z$ in ${\mathfrak{C}}$, we write
$\nat{C}_{f,g}$ for the resulting isomorphism of functors
$\nat{f}\nat{g}\iso \nat{(gf)}$. A similar convention applies for
covariant pseudofunctors.
\item The important ``pseudofunctor" from the point of view
of Grothendieck duality is the twisted inverse image ``pseudofunctor"
(see \Sref{s:suresh}).  Unfortunately
when we deal with non-ordinary schemes, this is not a pseudofunctor.
Indeed, if $\X$ is a noetherian scheme, then ${\boldsymbol{1}}_\X^!
= \R\iGp{\X}$ which is not in general isomorphic to 
${\boldsymbol{1}}_{\wDqcp(\X)}$.
The twisted inverse image is a pre-pseudofunctor in the sense
of Lipman.
A (contravariant) pre-pseudofunctor $\natb{\boldsymbol{-}}$ is data
of the form 
$$\left(\nat{(-)},\,\nat{(-)},\,\nat{C}_{(-),(-)},\,\nat{\delta}_{(-)}\right)$$
satisfying the requirements of \cite[\S\,4]{lns} for a pseudofunctor, with
the exception that $\delta_X$ is no longer required to
be isomorphism of functors for 
$X\in{\mathfrak{C}}$. We will often abuse terminology and
refer to pre-pseudofunctors as pseudofunctors. 
\item Let $f:\X\to \Y$ be a closed immersion of formal schemes (see
\cite[p.\,442]{gd}). Let $\I$ be the kernel of the surjective
map $\cO_\Y\twoheadrightarrow f_*\cO_\X$, and let $\overline{\Y}$ be
the ringed space $(\Y,\,\cO_\Y/\I)$. The natural map $\X\to \overline{\Y}$
is denoted $\bar{f}$. Note that $\bar{f}$ is {\it flat} and $f$
factors naturally as $\X\xrightarrow{\bar{f}}\overline{\Y}\xrightarrow{i}\Y$.
\item\label{item:etale} A {\it formally {\'e}tale} map of formal schemes
is a formally smooth map which is of relative dimension $0$. 
Such a map is called 
{\it {\'e}tale} if it is essentially of pseudo-finite type. Equivalently an
{\'e}tale map is a smooth map of relative dimension $0$. Relative
dimension here is as in \cite[Definition\,2.6.2]{lns}.
\item A map $(\X,\,\De')\xrightarrow{f}(\Y,\,\De)$ in $\bbFc$ (or in
$\cbbFc$, $\rbbFc$) is {\em smooth, {\'e}tale, pseudo-proper etc.}, if
the underlying map $\X\to \Y$ of formal schemes is smooth, {\'e}tale,
pseudo-proper etc.
\item For $(\X,\,\De) \in \bbFc$ and $\F, \G \in \Cozs{\De}(\X)$
$$
\Hom_\De(\F,\,\G)\set \Hom_{\Cozs\De}(\F,\,\G).
$$
\item For $\X\in \bbF$, $Q_\X$ will denote all the localization functors
from subcategories of ${\mathbf{K}}(\X)$ to $\D(\X)$.  The source of $Q_\X$
will be clear from the context. If the subcategory
is $\Cozs{\De}(\X)$ for a codimension function $\De$ on $\X$, then
we may sometimes restrict the target of $Q_\X$ to $\cm(\X,\,\De)$ (so
that $Q_\X$ becomes an equivalence between $\Cozs{\De}(\X)$ and
$\cm(\X,\De)$).
This will also be clear from the context.
\item If $R$ is a local ring, $\fm_R$ will denote its maximal ideal
and $k_R$ the residue field $R/\fm_R$. A $0$--dimensional $R$--module $M$
is (as in \cite{I-C}) a module satisfying $\Gamma_{\fm_R}M = M$.
\item Let $x$ be a point on a noetherian (formal) scheme $\X$,
and let $M$ be a $0$--dimensional $\cO_{\X,x}$ module. Then 
$i_xM$ will denote the sky-scraper sheaf on $\X$ whose sections
are $M$ over open sets containing $x$ and zero otherwise. Note that
$i_xM$ is a quasi-coherent $\cO_\X$--module (cf.~\cite[Lemma\,2.3.5]{lns}).
\item\label{item:sharp} As in \cite{lns}, for a complete local ring $A$, 
$A_{\sharp}$ will denote the category of $0$-dimensional $A$-modules.
The pseudofunctor on complete local rings given in 
\cite[Theorem\,4.3.1]{lns} will be denoted ${\boldsymbol{{-}_{\sharp}}}$.
\end{enumerate}


\section{\bf Traces}\label{s:trace}

In this section we define a trace map
$$
\Tr{f}(\F): f_*\sh{f}\F \longrightarrow \F
$$
associated with the data $(f,\F)$ where $f:(\X,\,\De')\to (\Y,\,\De)$
is a morphism in $\bbFc$ and $\F$ is a Cousin complex on $(\Y,\,\De)$.
More precisely $\Tr{f}(\F)$ is a map from ${\tt{Fgt}}f_*\sh{f}\F\to
{\tt{Fgt}}\F$ where ${\tt{Fgt}}$ is the forgetful functor from complexes
to graded objects. The most important result in this section is the
{\em Trace Theorem}, i.e., \Tref{thm:tr-thm}, which asserts that $\Tr{f}(\F)$
is compatible with coboundary maps on $f_*\sh{f}\F$ and $\F$ (i.e.,
$\Tr{f}(\F)$ is a map of complexes) {\em when $f$ is pseudo-proper}.
Related results may be found in \cite{IC2}.

Before we treat traces, we begin with more preliminary material 
involving local rings associated with a point on a formal scheme.

\subsection{Local rings}\label{ss:localrings} If $X$ is an ordinary scheme and
$x$ a point on $X$, then there are two naturally occurring
local rings associated with the point $x$---the stalk
$\cO_x$ of the structure sheaf $\cO_X$ at $x$ and the
completion ${\widehat{\cO}}_x$ of $\cO_x$ at its maximal
ideal. For a formal (non-ordinary) scheme there are 
four local rings that one can legitimately attach to
a point as we shall see.

To simplify the discussion let $\X=\Spf{(A,I)}$, where
$(A,I)$ is a noetherian adic ring. For $f\in A$, as is standard
in such situations, we set $A_{\{f\}}$ equal to the $I$--adic
completion of $A_f$. The points of $\X$ are in one-to-one
correspondence with open prime ideals of $A$, and the closed
points of $\X$ correspond to maximal ideals of $A$ (which are
necessarily open in the $I$--adic topology, 
since $A$ is complete in the $I$-adic topology,
whence $I$ is contained in the Jacobson radical of $A$). Let
$x$ be a point on $\X$, and let $\fp\subset A$ be the open
prime ideal corresponding to $x$. We then have an obvious
local ring associated with $x$, viz.~$A_{\fp}$. However, this
ring is not the stalk of $\cO_\X$ at $x$. One verifies that
the following formula gives the stalk at $x$:
$$
A_{\{\fp\}} \set \dirlm{f\notin\fp}A_{\{f\}} = \cO_{\X,x}.
$$
We have a canonical map 
$$
A_\fp \longrightarrow A_{\{\fp\}}
$$
whose $I$--adic completion is an isomorphism (see proof of
\cite[Lemma\,7.1.1]{dfs}). Denote the common $I$-adic completion
of $A_\fp$ and $A_{\{\fp\}}$ by $A_{[\fp]}$ and the common
$\fp$--adic completion of $A_\fp$, $A_{\{\fp\}}$ and $A_{[\fp]}$
by $\widehat{A}_\fp$. We then have faithfully flat inclusions
of local rings
$$
A_{\fp} \subset A_{\{\fp\}} \subset A_{[\fp]} \subset \widehat{A}_\fp.
$$
We also have formal schemes 
$$
\X_{[\fp]} \set \Spf{(A_{[\fp]},IA_{[\fp]})}
$$
and
$$
\X^*_\fp \set \Spf{(\widehat{A}_\fp,\fp\widehat{A}_\fp)}
$$
together with natural {\'e}tale maps (see 
\ref{ss:conventions}\eqref{item:etale})
$$
\X^*_\fp \longrightarrow \X_{[\fp]} \longrightarrow \X.
$$
Note that the map $\X_{[\fp]}\to \X$ is {\it{adic}}, but the
map $\X^*_\fp \to \X_{[\fp]}$ (as well as the above composite)
need not be so.

\subsection{Trace at the graded level} Let 
${\boldsymbol{-}}_{\boldsymbol{\sharp}}$ be as in 
\ref{ss:conventions}\eqref{item:sharp}. Suppose $\varphi: R\to S$ is
a morphism in ${\mathfrak C}_{\mathbf{rf}}$---i.e., $R$, $S$ are complete
noetherian local rings and $\varphi$ is a local homomorphism with
$k_S$ a finite $k_R$--algebra---and suppose
$M$ is an object in $R_{\sharp}$. According to \cite[Chapter\,7]{I-C},
we have an $R$--linear trace map, functorial in $M\in R_{\sharp}$
\stepcounter{thm}
\begin{equation*}\label{map:ic-trace}\tag{\thethm}
\Tr{S/R,M}: \varphi_{\sharp}M \longrightarrow M
\end{equation*}
such that the induced map
\stepcounter{thm}
\begin{equation*}\label{map:ic-phi}\tag{\thethm}
\Phi_{S/R,M} :\varphi_{\sharp}M\longrightarrow \Hom_R^c(S,\,M)
\end{equation*}
is an isomorphism of $S$--modules. Moreover if $\varphi$ is {\em surjective},
\eqref{map:ic-phi} recovers the isomorphism in 
\cite[Theorem\,4.3.1\,I(ii)]{lns}. Note that we are using the fact that
$R_{\sharp}=R_{\#}$, $S_{\sharp}=S_{\#}$ and $\varphi_{\sharp}=\varphi_{\#}$
where ${\boldsymbol{{-}_{\#}}}$ is Huang's pseudofunctor on ${\mathfrak C}$,
the category of of complete noetherian local rings and local homomorphisms
(cf. proof of \cite[Theorem\,4.3.1]{lns}, especially Lemma\,4.3.2). Further,
by the commutative diagram on the top of \cite[p.\,51]{I-C}, if
$\xi:S\to T$ is a second morphism in ${\mathfrak C}_{\mathbf{rf}}$ then
the diagram
$$
\xymatrix{
\xi_{\sharp}\varphi_{\sharp}M \ar[rr]^{\Iso}_{C_{\sharp}^{\varphi,\xi}}
\ar[d]_{\Phi_{T/S}, \Phi_{S/R}} & & (\xi\varphi)_{\sharp}M 
\ar[d]^{\Phi_{T/R}} \\
\Hom_S^c(T,\,\Hom_R^c(S,\,M)) \ar[rr]^{\Iso}_{\mathrm{natural}}
& & \Hom_R^c(T,\,M)
}
$$
commutes. We are using the fact that the transcendence degrees of the
residue field extension $k_R\to k_S$ and $k_S\to k_T$ are zero, whence
$C^{\varphi,\xi}_{\sharp} = C^{\varphi,\xi}_{\#}$ by the construction
of ${\boldsymbol{{-}_{\sharp}}}$ from ${\boldsymbol{{-}_{\#}}}$ via
\cite[Lemma\,4.3.2]{lns} in the proof of [{\em Ibid},\,Theorem\,4.3.1].
As a consequence the following diagram commutes
\stepcounter{thm}
\begin{equation*}\label{diag:tr-tr}\tag{\thethm}
\xymatrix{
\xi_{\sharp}\varphi_{\sharp}M \ar[rr]^{\Iso}_{C^{\varphi,\xi}_{\sharp}}
\ar[d]_{\Tr{T/S}} & & (\xi\varphi)_{\sharp}M \ar[d]^{\Tr{T/R}} \\
\varphi_{\sharp}M \ar[rr]_{\Tr{S/R}} & & M
}
\end{equation*}

Next consider datum of the form
$$
D = \left((\X,\,\De')\xrightarrow{f} (\Y,\,\De),\,\F,\,x\right)
$$
where $f$ is a map in $\bbFc$, $\F$ is an object in $\Cozs{\De}(\Y)$
and $x$ is a point on $\X$. Associated to $D$ we define a {\em punctual
trace} at $x$, obviously functorial in $\F$,
\stepcounter{thm}
\begin{equation*}\label{map:tr-x}\tag{\thethm}
\Tr{f,x}(\F) = \begin{cases}
{\Tr{S/R,M}:\varphi_{\sharp}M\to M \qquad{\text{if}}\quad\De'(x)=
\De(f(x))} \\
0 \qquad{\text{otherwise}}\end{cases}
\end{equation*}
where, in the first case, $R$ and $S$ are the completions of the local
rings at $x$ and $y=f(x)$ respectively, $\varphi: R\to S$ is the natural
map induced by $f$, and $M$ is $\F(y)$. We remind the reader that the
condition $\De'(x)=\De(f(x))$ implies that $x$ is closed in the
fiber of $f$ over $f(x))$ (and when $f$ is of pseudo-finite type,
the conditions are equivalent). Varying $x$, we get a map of graded 
$\cO_\Y$--modules 
\stepcounter{thm}
\begin{equation*}\label{map:tr-f}\tag{\thethm}
\Tr{f}(\F):f_*\sh{f}\F \to \F
\end{equation*}
given by 
$$
\Tr{f}(\F) = \sum_{x\in\X} i_{f(x)}\Tr{f,x}(\F).
$$
More precisely, as we noted before, the source of $\Tr{f}(\F)$ is
$f_*{\tt{Fgt}}_\X\sh{f}\F$ and the target is ${\tt{Fgt}}_\Y\F$ where
${\tt{Fgt}}$ is the forgetful functor from complexes to graded modules.

If $g:(\W,\,\De'') \to (\X,\,\De')$ is a second map in $\bbFc$ then from
\eqref{diag:tr-tr} we see that the diagram
\stepcounter{thm}
\begin{equation*}\label{diag:tr-tr-f}\tag{\thethm}
\xymatrix{
(fg)_*\sh{g}\sh{f}\F \ar@{=}[d] \ar[rr]^{\Iso}_{\sh{C}_{g,f}}
& & (fg)_*\sh{(fg)}\F \ar[dd]^{\Tr{fg}} \\
f_*g_*\sh{g}\sh{f}\F \ar[d]_{f_*\Tr{g}} & & \\
f_*\sh{f}\F \ar[rr]_{\Tr{f}} & & \F
}
\end{equation*}
commutes.

We end this subsection by discussing $\Tr{f}$ when $f\colon \X\to \Y$
is a closed immersion. By \cite[\S\S\,8.2, Def.\,8.3.1 and \S\S\,8.4]{lns}, 
for $\F\in \Cozs{\De}(\Y)$ we have an isomorphism
\stepcounter{thm}
\begin{equation*}\label{iso:closed-imm}\tag{\thethm}
f_*\sh{f}\F \iso \sHomb_\Y(f_*\cO_\X,\,\F).
\end{equation*}
We then have
\begin{lem}\label{lem:trThm-climm}
Let $f\colon (\X,\,\De')\to (\Y,\,\De)$ be a closed immersion in $\bbFc$.
Let ${\boldsymbol{e}}\colon \sHomb_\Y(f_*\cO_\X,\,\F) \to \F$ be 
``evaluation at 1". Then $Tr{f}(\F)$ is the composite
$$
f_*\sh{f}\F \xrightarrow{\eqref{iso:closed-imm}} \sHomb_\Y(f_*\cO_\X,\,\F)
\xrightarrow{\boldsymbol{e}} \F.
$$
In particular $\Tr{f}(\F)$ is a map of complexes.
\end{lem}

\proof
This follows from \cite[\S\S\,8.2 (75)]{lns} and from the fact that the map
in \cite[Thm.\,4.3.1\,I(ii)]{lns} is the map
\eqref{map:ic-phi} when $\varphi$ is a surjective map. 
\qed

\subsection{Relative projective space}\label{ss:proj-space} 
The key step for proving the general Trace Theorem
is the Trace Theorem for relative projective
spaces $\pi: \bbP^d_Y \to Y$, after which the proof in 
\cite[VII\,\S\,2,pp.\,369--373]{RD} applies mutatis mutandis. 

Throughout this subsection $Y$ is a fixed {\it ordinary}
noetherian scheme. Let 
$$
\pi=\pi_Y: \bbP\set\bbP^d_Y \longrightarrow Y
$$
be the relative projective space of fiber dimension $d$ over $Y$, i.e.,
$\pi$ is the first projection in the decomposition $\bbP^d=Y\times_{\mathbb Z}
\bbP^d_{\mathbb Z}$. There is a well known isomorphism
\stepcounter{thm}
\begin{equation*}\label{iso:int-pi}\tag{\thethm}
\int_\pi:\Rr^d\pi_*\omega_\pi \iso \cO_Y
\end{equation*}
(cf. \cite[2.1.12]{ega3} or \cite[p.152,\,Theorem 3.4]{RD}); the 
generating section $\sigma$ of
$\Rr^d\pi_*\omega_\pi$ corresponding to the standard section $1$
of $\cO_Y$ described as follows.  We have 
$\bbP={\mathbf{Proj}}\left(\cO_Y[T_0,\ldots,T_d]\right)$. Let 
${\mathscr U}=(U_i)_{i=0}^d$ be the open cover of $\bbP$
given by $U_i=\{T_i\ne 0\}$. On $U_0\cap\ldots\cap U_d$
we have inhomogeneous coordinates $t_i=T_i/T_0$, $i=1,\dots,d$
whence a section
$$
\check{\sigma}_T\set \frac{{\mathrm{d}}t_1\wedge\ldots\wedge{\mathrm{d}}
t_d}{t_1\dots t_d}
\in \omega_\pi(U_0,\dots,U_d).
$$
We have an isomorphism
$$
\Rr^d\pi_*\omega_\pi \iso \Hr^d(\pi_*\check{\mathcal C}^\bullet
({\mathscr U},\,\,\omega_\pi))
$$
and $\check{\sigma}_T$ has a natural image in the right side as a
\v{C}ech cohomology class. Let $\sigma$ be the corresponding
element on the left side.  The section $\sigma$ does not depend 
on the choice of homogeneous coordinates $T_0,\dots,T_d$ of $\bbP$ 
(cf. \cite[p.34,\,Lemma\,2.3.1]{conrad}) and is the sought after
section.

It is well known (and easily verified from the description above) that
$\int_\pi$ is compatible with arbitrary base changes $Y'\to Y$.
If $Y$ is affine, say $Y=\Spec(A)$, then define
$$
\int_{\bbP/Y} \colon \Hr^d(\bbP,\,\omega_\pi) 
\longrightarrow A.
$$
as the global section of $\int_\pi$.

Since the global trace map $\Tr{\pi}$ is built out of punctual
traces, we look into these traces in somewhat greater detail now.
Suppose $\varphi:R\to S$ is a smooth residually finite map of relative
dimension $d$ between complete local rings and $M$ is an object in
$R_{\sharp}$. From \cite[Theorem\,4.3.1,\,I(i)]{lns} we have
\stepcounter{thm}
\begin{equation*}\label{iso:sharp-h}\tag{\thethm}
\varphi_{\sharp}M \iso \Hr^d_{\fm_S}(M\otimes\omega_{S/R})
\end{equation*}
and this isomorphism is functorial in $M$. This gives us a residue map
\stepcounter{thm}
\begin{equation*}\label{map:res-ic}\tag{\thethm}
\res{S/R,M}:\Hr^d_{\fm_S}(M\otimes\omega_{S/R}) \to M
\end{equation*}
defined by the commutativity of
$$
\xymatrix{
\Hr^d_{\fm_S}(M\otimes\omega_{S/R}) \ar[rd]^{\res{S/R,M}} & \\
& M \\
\varphi_{\sharp}M \ar[uu]^{\eqref{iso:sharp-h}}_\wr 
\ar[ur]_{\,\Tr{S/R,M}} &
}
$$
where $\Tr{S/R,M}$ is as in \eqref{map:ic-trace}. If $S$ is of
the form
$$
S=R[t_1,\dots,t_d]_{\fp}\widehat{{}}
$$
where $t_1,\dots,t_d$ are independent variables over $R$,
then Huang has explicit formulae for $\res{S/R}$ (cf. 
\cite[p.42,\,(7.1)]{I-C} and \cite[bottom of p.21]{I-C}) based on
an iterated form of the Tate residue.

Recall again that in this subsection $Y$ is a fixed noetherian
ordinary scheme.

\begin{prop}\label{prop:residue-thm} Suppose $Y$ is an affine
noetherian scheme, $x$ a point on $\bbP$ closed in the fiber
of $\pi$ over $y=\pi(x)$ and $M$ a zero-dimensional $\cO_{Y,y}$--module.
Denote $\widehat{\cO}_{Y,y}$, $\widehat{\cO}_{X,x}$ and $i_yM$ by
$R$, $S$ and $\F$ respectively. Let $W$ be the closure of $x$ in $\bbP$.
Then the diagram below---with unlabelled arrows the natural maps---commutes.
$$
\xymatrix{
\Hr^d_x(\pi^*\F\otimes\omega_\pi) \ar@{=}[d] & 
\Hr^d_{\fm_S}(M\otimes\omega_{S/R}) \ar[l]_{\Iso} \ar[dr]^{\res{S/R,M}} & \\
\Hr^d_W(\pi^*\F\otimes\omega_\pi) \ar[d]  & & M \\
\Hr^d(\bbP,\,\pi^*\F\otimes\omega_\pi) & 
M\otimes\Hr^d(\bbP,\,\omega_\pi) \ar[l]_{\,\,\Iso} \ar[ur]_{1\otimes\int_{\bbP/Y}}
&  
}
$$
\end{prop}

\proof First we elaborate on the equality on the top of the left column
in the above diagram. Let $Y'=\Spec{R}$, $y'\in Y'$ the closed point of
$Y'$, $Y'\xrightarrow{g}Y$, ${\bbP}'\set{\bbP}^d_{Y'}\xrightarrow{g'} {\bbP}$
the natural affine maps, $\pi'\colon {\bbP}'\to Y'$ the natural projection
and $F$ the common fiber ${\pi'}^{-1}(y') = \pi^{-1}(y)$. Since the
fiber $F$ is shared by ${\bbP}$ and ${\bbP}'$ therefore we have a 
point $x'\in \bbP'$ corresponding to $x\in \bbP$. Let $\G={\pi}^*\F\otimes
\omega_\pi$ and $\G'={g'}^*\G$ ($={\pi'}^*(g^*\F)\otimes\omega_{\pi'}$).
It is not hard to see that
$$
\Hr^d_{x'}(\G') = \Hr^d_x(\G).
$$
More generally we remind the reader that if $(B,\eta)$ is a local ring,
$N$ a $B$--module, and $({\widehat{B}},\,{\widehat{\eta}})$ the completion
of $(B,\,\eta)$ at $\eta$, then for every integer $i$ we have
$\Hr^i_\eta(N) = \Hr^i_\eta(N)\otimes_B{\widehat{B}}=\Hr^i_{\widehat{\eta}}
(N\otimes_B{\widehat{B}})$.\footnote{In \cite[\S\S\,4.2, (18)]{lns} we 
made an elaborate distinction between these modules, necessary,
to make sure that our book-keeping was correct when constructing the
pseudofunctor $\sha{\boldsymbol{-}}$.}

Now, $\iG{F}\G'=\G'$ and ${\underline{\Hr}}^i_F(\G')\set\Hr^i(\iG{F}\G')\simeq
{\pi'}^*i_{y'}(\Hr^i_{\fm_{y'}}(M))\otimes \omega_{\pi'} = 0$ for 
$i>0$, since $M$ is a zero dimensional $R$--module.  One can therefore
find an $\cO_{Y'}$--injective resolution $\G'\to \I$ such that
$\iG{F}\I=\I$, i.e.~an injective resolution (topologically)
supported on $F$ (by replacing any injective resolution $\J$ by $\I=\iG{F}\J$). 
Moreover, $g'$ being affine, $g'_*\I$ is a resolution of 
$g'_*\G'=g'_*{g'}^*(\pi^*\F\otimes\omega_{\pi})=g'_*{g'}^*(\pi^*\F)\otimes
\omega_\pi = \G$ (the last equality follows from the fact that $M$ is
a zero dimensional $R$-module, whence $\F=g_*g^*\F$ giving $g'_*{g'}^*(\pi^*\F)
= \pi^*\F$). Let $z\in {\bbP}$ be such that $\pi(z)=y$, and let
$z'\in {\bbP}'$ be the point lying over $y'$ corresponding to $z$.
Since the local rings $\cO_{\bbP,z}$ and $\cO_{\bbP',z}$ have the same
completion for such pairs of points---i.e.~the natural map $\cO_{\bbP,z}
\to \cO_{\bbP',z'}$ transforms into an equality on completing---and since
$\iG{F}\I=\I$, therefore it is not hard to see that $g'_*\I$ consists of 
$\cO_Y$--injectives. Thus $g'_*\I$ is an injective resolution of $\G$. 
We therefore have
\begin{equation*}
\begin{split}
\Hr^d_W(\G) = \Hr^d(\Gamma_Wg'_*\I) & = \Hr^d(\Gamma_{{g'}^{-1}W}\I) \\
& = \Hr^d(\Gamma_{x'}\I) \\
& = \Hr_{x'}^d(\G') \\
& = \Hr_x^d(\G)
\end{split}
\end{equation*}
which explains the top of the left column.

Let $\wres{S/R,M}\colon 
\Hr^d_{\fm_S}(M\otimes\omega_{S/R}) \to M$ be the composite
\begin{align*}
\Hr^d_{\fm_S}(M\otimes\omega_{S/R}) & \iso  
\Hr^d_x(\pi^*\F\otimes \omega_\pi) \\
& \,\,= {\phantom{X}} \Hr^d_W(\pi^*\F\otimes\omega_\pi) \\
& \longrightarrow\,\,\, \Hr^d(\bbP,\,\pi^*\F\otimes\omega_\pi) \\
& \iso \, M\otimes\Hr^d(\bbP,\,\omega_\pi) \\
& \longrightarrow \,\,\,M
\end{align*}
where the last arrow is $1\otimes\int_{\bbP/Y}$.
We have to show that $\res{S/R}=\wres{S/R}$.

Since $\int_\pi$ is compatible with the flat base change $\Spec(R)\to Y$,
and so are all other
functorially defined arrows in the diagram, we assume 
without loss of generality that $Y=\Spec(R)$. Assume---again without
loss of generality---that $x$ lies in the open subscheme $U_0=
\Spec(R[T_1/T_0,\dots,T_d/T_0])$ of $\bbP={\mathrm{Proj}}(R[T_0,\dots,T_d])$.
We therefore have 
$$
S=R[T_1/T_0,\dots,T_d/T_0]_{\fp}\widehat{{}}
$$ 
for some maximal ideal $\fp$ of the polynomial ring $R[T_1/T_0,\dots,T_d/T_0]$
and \cite[p.42,\,(7.1)]{I-C} applies. If the natural extension $k_R\to k_S$
is trivial then $S=R[[t_1,\dots,t_d]]$ for some analytically independent
variables $t_1,\dots,t_d$ over $R$ and the formula in {\it loc.cit.}
reduces to 
\begin{equation}\label{eq:ic-stn}
\res{S/R}\begin{bmatrix}
m\otimes \mathrm{d}t_1\wedge\dots\wedge{\mathrm{d}}t_d \\
t_1^{\alpha_1},\dots,t_d^{\alpha_d}
\end{bmatrix} 
= \begin{cases}m \qquad {\text{if}}\,\, \alpha_1=\alpha_2=\dots=\alpha_d=1 \\
0 \qquad \,\text{otherwise}\end{cases} 
\end{equation}

Suppose we are in this situation, i.e., $k_R=k_S$ or equivalently
$S=R[[t_1,\dots,t_d]]$. Since the two maps, $\int_\pi$ and $\res{S/R}$, are
compatible with change of homogeneous coordinates $T_0,\dots,T_d$
(cf. \cite[p.49,(7.7)]{I-C} for the latter assertion), we may assume
without loss of generality that $t_i=T_i/T_0$ for $i=1,\dots,d$.
Let ${\mathbf{t}}=(t_1,\dots,t_d)$, and let $Z\hookrightarrow
U_0$ be the closed subscheme given by the
vanishing of the $t_i$'s. We define 
$$
\wres{\mathbf{t}}:\Hr^d_{{\mathbf{t}}S}(\omega_{S/R})
\longrightarrow R
$$
by the commutativity of
$$
\xymatrix{
\Hr^d_Z(\bbP,\,\omega_\pi)\ar[d] \ar[rr]^{\Iso} & & 
\Hr^d_{{\mathbf t}S}(\omega_{S/R})
\ar[d]^{\wres{\mathbf{t}}} \\
\Hr^d(\bbP,\,\omega_\pi) \ar[rr]_{\int_{\bbP/Y}} & & R
}
$$
Lipman's proof \cite[p.\,75,\,Prop.\,(8.5)]{ast-117} of his 
``Residue Theorem for Projective Space" applies without change to
our situation\footnote{In fact since $(1:0:\dots:0)$ are
the ``homogeneous coordinates" of the $R$--valued point $\Spec(R)
\iso Z \hookrightarrow \bbP$, only that part of the proof of {\it loc.cit.}
which concerns itself with rational points applies, i.e. the
part that begins at the bottom of p.77.} and we
conclude that
\begin{equation}\label{eq:joe-stn}
\wres{\mathbf{t}}\begin{bmatrix}
p({\mathbf{t}}){\mathrm{d}}t_1\wedge\dots\wedge{\mathrm{d}}t_d \\
t_1^{\alpha_1},\dots,t_d^{\alpha_d}
\end{bmatrix}
= p_{\alpha_1-1,\dots,\alpha_d-1}
\end{equation}
where the right-side is the coefficient of 
$t_1^{\alpha_1-1}\dots t_d^{\alpha_d-1}$ in the power series $p({\mathbf{t}})$.
It is straightforward to check that the diagram
$$
\xymatrix{
\Hr^d_{{\mathbf{t}}S}(M\otimes\omega_{S/R}) \ar[rr]^{\Iso} 
& & \Hr^d_{\fm_S}(M\otimes\omega_{S/R})
\ar[d]^{\wres{S/R,M}} \\
M\otimes\Hr^d_{{\mathbf{t}}S}(\omega_{S/R}) \ar[u]_{\wr}
\ar[rr]^{1\otimes\wres{\mathbf{t}}}
& & M
}
$$
commutes. Whence, on comparing \eqref{eq:ic-stn} with \eqref{eq:joe-stn},
we see that $\wres{S/R}=\res{S/R}$ in this case, i.e. when the natural
extension $k_R\to k_S$ is trivial.

One can reduce to the above case by making a flat base change $R\to R'$
with $R'$ a complete noetherian local ring satisfying $\fm_RR'=\fm_R$.
In greater detail, suppose $k_S=k_R[\theta_1,\dots,\theta_m]$. We can
find $R'$ as above such that for each $i=1,\dots,m$, the minimal
polynomial of $\theta_i$ over $k_R$ splits into a product of linear
factors over $k_{R'}$. Let $\bbP'=R'\otimes_R\bbP$ and let the
resulting projections be $p:\bbP'\to \bbP$ and $\pi':\bbP'\to \Spec(R')$.
Let $\{x_1,\dots,x_n\}=p^{-1}(x)$ and for $i=1,\dots,m$ let $S_i'$
be the completion of $\cO_{\bbP',x_i}$. If $M'=R\otimes_RM$, then
according to \cite[p.47,\,Lemma\,(7.6)]{I-C} we have
\begin{equation}\label{eq:ic-bc}
\sum_{i=1}^n\res{S'_i/R',M'}\circ\tau = \tau_o\circ\res{S/R,M}
\end{equation}
where $\tau$ is the natural map
$$
\Hr^d_{\fm_S}(M\otimes\omega_{S/R}) \longrightarrow
R'\otimes\Hr^d_{\fm_S}(M\otimes\omega_{S/R}) =
\bigoplus_{i=1}^n\Hr^d_{\fm_{S_i}}(M'\otimes\omega_{S'_i/R'})
$$
and $\tau_o$ is the natural map
$$
M \to R'\otimes_RM=M'.
$$

On the other hand, since $\int_{\bbP/Y}$ is compatible with base
change, clearly we have
\begin{equation}\label{eq:pi-bc}
\sum_{i=1}^n\wres{S'_i/R',M'}\circ\tau = \tau_o\circ\wres{S/R,M}.
\end{equation}
Since $k_{S'_i}=k_{R'}$, we have $\wres{S'_i/R'}=\res{S'_i/R'}$.
Now $\tau_o$ is an injective map and hence \eqref{eq:ic-bc} and
\eqref{eq:pi-bc} give the Proposition.
\qed

\begin{prop}{\em(Trace Theorem for Projective Space.)}\label{prop:proj-sp}
Let $Y$ be in $\bbF$ and $\F$ in $\Cozs{\De}(Y)$ for some codimension
function $\De$ on $Y$. Then the map of graded $\cO_Y$--modules
$$
\Tr{\pi}(\F)\colon \pi_*\sh{\pi}\F \longrightarrow \F
$$
is a map of complexes
\end{prop}

\proof Our strategy is as follows. We will show that there is a complex
$\E=\E_\F$ of quasi-coherent $\cO_Y$--modules associated with $\sh{\pi}\F$
together with maps of complexes (the second one an isomorphism)
$$
\alpha\colon\pi_*\sh{\pi}\F \to \E
$$
and
$$
\beta\colon\E \iso \F\otimes{\mathrm R}^d\pi_*\omega_\pi
$$
such that the diagram of graded $\cO_Y$--modules
\begin{equation}\label{diag:alpha-beta}
\xymatrix{
\pi_*\sh{\pi}\F \ar[r]^\alpha \ar[d]_{\Tr{\pi}(\F)} & \E \ar[d]^{\beta}\\
\F & \ar[l]^-{1\otimes\int_\pi} \F\otimes{\mathrm R}^d\pi_*\omega_\pi
}
\end{equation}
commutes. Since all arrows other than $\Tr{\pi}(\F)$ are maps of
complexes, this would prove that $\Tr{\pi}(\F)$ is a map of complexes.

The complex $\E$ is best described as the $E_1$ term of a 
spectral sequence associated to a natural filtration of the
complex $\pi_*\sh{\pi}\F$. Recall from \cite[\S\S\,10.2]{lns} (especially
the discussion following the proof of 10.2.4 in {\it Ibid.})
that we have a decreasing filtration $\{F^p\}_{p\in{\mathbb Z}}$ of 
subcomplexes of $\sh{\pi}\F$ given by
$$
F^p=\sh{\pi}(\sigma_{\ge p}\F),
$$
and recall that $F^p/F^{p+1}=\sh{\pi}(\F^p[-p])$ \cite[(89)]{lns}.
$F^\bullet$ induces a decreasing filtration $\{\Phi^p\}_{p\in{\mathbb Z}}$ on
$\pi_*\sh{\pi}\F$ obtained by applying $\pi_*$ to $F^\bullet$. Let
$\{E^{p,q}_r\}$ be the associated spectral sequence. Then $\E$ is the
complex
$$
 \dots \to E_1^{p-1,0} \xrightarrow{\partial_\E^{p-1}} E_1^{p,0}
\xrightarrow{\partial_\E^p} E_1^{p+1,0} \to \dots
$$
where $\partial_\E$ is the natural coboundary on the $E_1$ term of
the spectral sequence $\{E^{p,q}_r\}$.

To say more (i.e.~to define $\alpha$ and $\beta$) we need to deconstruct
the above definition of $\E$. We view $\sh{\pi}\F$ (resp. $\pi_*\sh{\pi}\F$)
as a bigraded complex (not necessarily a double complex!) whose $p$--th
column is the complex $F^p/F^{p+1}$ (resp. $\Phi^p/\Phi^{p+1}$). In
greater detail let
$$
\A^{p,q} \set (F^p/F^{p+1})^{p+q}.
$$
Then, with $\De'=\sh{\pi}\De$, we have
$$
\A^{p,q}=\bigoplus_{\substack{\De(\pi(x))=p\\ \De'(x)=p+q}} i_x(\sh{\pi}\F)(x).
$$ 
In \cite[\S\S\,10.1]{lns} the above sheaf is denoted $\E^{p+q,p}$. One has
the following decomposition of $\cO_{\bbP}$--modules
$$
(\sh{\pi}\F)^n = \bigoplus_{p+q=n}\A^{p,q} \qquad \qquad 
F^p=\bigoplus_{r\ge p}\A^{r,q}
$$
as well as the associated decomposition of $\cO_Y$--modules
$$
(\pi_*\sh{\pi}\F)^n = \bigoplus_{p+q=n}\pi_*\A^{p,q} \qquad \qquad
\Phi^p = \bigoplus_{r\ge p}\pi_*\A^{r,q}.
$$
Note that the map
$$
\partial^{p,q,k}\colon \A^{p,q} \to \A^{p+k,q-k+1}
$$
induced by the coboundary map on $\sh{\F}$ is such that $\partial^{p,q,k}=0$ 
if $k$ is negative, i.e.~there are no arrows with a westward component
in the bigraded complex (cf. \cite[Lemma 10.2.3]{lns}).
$$
\xymatrix{
\A^{p,q+1} & & & & \\
\A^{p,q} \ar[u]^{\partial^{p,q,0}} \ar[r]^{\partial^{p,q,1}} 
\ar[rrd]^{\partial^{p,q,2}} \ar@{.>}[rrrdd] \ar[rrrrddd]_{\partial^{p,q,k}} 
& \A^{p+1,q} & & & \\
& & \A^{p+2,q-1} \ar@{.}[rrdd] & &\\
& & & & \\
& & & & \A^{p+k,q-k+1} \\
}
$$
Recall that the $E_1$ term of $\{E^{p,q}_r\}$ has a simple description
which in our case translates to having
$$
\E^p = \Hr^p(\Phi^p/\Phi^{p+1})
$$
such that the coboundary map $\partial_\E^p$ is the connecting homomorphism 
associated to the short exact sequence of complexes
\begin{equation}\label{cplx:connE1}
0 \to \Phi^{p+1}/\Phi^{p+2} \to \Phi^p/\Phi^{p+2} \to \Phi^p/\Phi^{p+1}
\to 0.
\end{equation}

Note that $\Phi^p/\Phi^{p+1} = \pi_*(F^p/F^{p+1})$ is the complex whose
$(p+q)$-th homogeneous piece is $\pi_*\A^{p,q}$ and whose coboundary
on $(p+q)$-cochains is $\pi_*\partial^{p,q,0}$. The complex
$\Phi^p/\Phi^{p+2}$ is a two column double complex\footnote{Which means,
for this proof, the grids anti-commute.}; its left column---the $p$-th
column---being 
$$
\Phi^p/\Phi^{p+1}=\pi_*\sh{\pi}(\F^p[-p])
$$
and its right column---the $(p+1)$-th column---being 
$$
\Phi^{p+1}/\Phi^{p+2}=\pi_*\sh{\pi}(\F^{p+1}[-p-1]).
$$
$$
\xymatrix{
0 & 0 \\
\pi_*\A^{p,0} \ar[u] \ar[r]^{\partial^{p,0,1}} & \pi_*\A^{p+1,0} \ar[u]\\
{\phantom{X}}\ar[u]^{\partial^{p,-1,0}} 
& {\phantom{X}}\ar[u]_{\partial^{p+1,-1,0}}\\
{\phantom{X}}\ar@{.}[u] & {\phantom{X}}\ar@{.}[u] \\
\pi_*\A^{p,-d+1} \ar[u]^{\partial^{p,-d+1,0}} 
\ar[r]^{\partial^{p,-d+1,1}}
& \pi_*\A^{p+1,-d+1} \ar[u]_{\partial^{p+1,-d+1,0}} \\
\pi_*\A^{p,-d} \ar[u]^{\partial^{p,-d,0}} \ar[r]^{\partial^{p,-d,1}}
& \pi_*\A^{p+1,-d} \ar[u]_{\partial^{p+1,-d,0}} \\
0 \ar[u] & 0 \ar[u] \\
\Phi^p/\Phi^{p+1} & \Phi^{p+1}/\Phi^{p+2}
}
$$

Now $\A^{p,q}=0$ for $q>0$ and hence $\Phi^p/\Phi^{p+1}$ has no terms
in degrees $k>p$. Therefore we have a natural surjective map
\begin{equation}\label{eq:A-phi}
{\bar\alpha}^p = {\bar\alpha}^p(\F) \colon \pi_*\A^{p,0}=(\Phi^p/\Phi^{p+1})^p
\twoheadrightarrow \Hr^p(\Phi^p/\Phi^{p+1}) = \E^p.
\end{equation}
For future reference we note that
\begin{equation}\label{eq:phi-p-f}
{\bar{\alpha}}^p(\F)={\bar\alpha}^p(\F^p[-p]).
\end{equation}
The map $\alpha\colon \pi_*\sh{\pi}\F\to\E$ at the graded level
is defined in degree $p$ by the composite
$$
(\pi_*\sh{\pi}\F)^p \xrightarrow{\text{projection}} \pi_*\A^{p,0}
\xrightarrow{\psi^p} E_1^{p,0}= \E^p.
$$
We have to check that $\alpha$ is a map of complexes.  Since
$\partial^{p,q,k}$ is zero for negative $k$, we only need to check
that $\partial^p_\E\circ{\bar{\alpha}}^p = 
\partial^{p,0,1}\circ{\bar{\alpha}}^{p+1}$.
Now, $\partial_\E^p$ is the connecting map associated with
\eqref{cplx:connE1}, and therefore this is easy to verify (the
diagram above might help). 

We next define the isomorphism
$\beta\colon\E \iso \F\otimes_Y \Rr^d\pi_*\omega_\pi$.
Let 
$$
{\mathbb{E}}_\pi(\F) \set \Ed{\De'}(\pi^*\F\otimes(\omega_\pi[d])).
$$
Consider the diagram of complexes
$$
\xymatrix{
{\mathbb{E}}_\pi(\F) \ar[rr]^{\Iso}_{L_\F} & & \sh{\pi}\F \\
& \pi^*\F\otimes(\omega_\pi[d]) \ar[ul]^{\eta_{\pi}(\F)} 
\ar[ru]_{\gamma_{\pi}(\F)}
}
$$
where $\eta_\pi(\F)$ is the map in \cite[\S\S\,10.1,\,(87)]{lns}, $L=L_\F$ 
the isomorphism of Cousin complexes is \cite[\S\S\,8.1]{lns} (see second
paragraph of {\it loc.cit.}) and
$\gamma_{\pi}(\F)$ is the composite $L\circ\eta_\pi(\F)$. By definition
of $\gamma_\pi$ the above diagram commutes. Some notations become
necessary at this point. To that end, for every pair of integers $p<q$ define
\begin{equation*}
\begin{split}
\F_{p,q} & \set \sigma_{\ge p}\F/\sigma_{\ge q}\F, \\
\A_p & \set \pi^*\F_{p,p+1}\otimes(\omega_\pi[d]) 
= (\pi^*\F\otimes\omega_\pi)[d-p], \\
\B_p &\set \pi^*(\F_{p,p+2})\otimes(\omega_\pi[d]). 
\end{split}
\end{equation*}
Note that $\sh{\pi}{\F_{p,q}}=F^p/F^q$ \cite[\S\S\,10.2,\,(89)]{lns}. 
The map $\eta_{\pi}(\F)$ is functorial in $\F\in \Cozs{\De}(\Y)$ as is $L_\F$. 
Therefore $\gamma_\pi(\F)$
is also functorial in $\F$. We therefore have a commutative diagram of
complexes of $\cO_\bbP$--modules with exact rows
$$
\xymatrix{
0 \ar[r] &F^{p+1}/F^{p+2} \ar[r] & F^p/F^{p+2} 
\ar[r] & F^p/F^{p+1} \ar[r] & 0\\
0 \ar[r] &\A_{p+1} \ar[r] \ar[u]^{\gamma_\pi(\F_{p+1,p+2})} & \B_p \ar[r] 
\ar[u]_{\gamma_\pi(\F_{p,p+2})}& \A_p \ar[r] \ar[u]_{\gamma_\pi(\F_{p,p+1})}& 0
}
$$
Denote the exact row on the top $E_t$ and the one at the bottom $E_b$. Now,
by \cite[Prop.\,10.1.6]{lns}, $\eta_\pi$ is a quasi-isomorphism and hence
so is $\gamma_\pi$. On applying $\R\pi_*$ to the above diagram and
identifying $\R\pi_*E_t$ with $\pi_*E_t$ (for $E_t$ consists of flasque
complexes) we get a map (in fact an isomorphism) of triangles in $\Dqc(Y)$:
$$
\xymatrix{
\Phi^{p+1}/\Phi^{p+2} \ar[r] &
\Phi^p/\Phi^{p+2} \ar[r] &
\Phi^p/\Phi^{p+1} \ar[r]^-{\R\pi_*\chi} &
\Phi^{p+1}/\Phi^{p+2}[1] \\
\R\pi_*\A_{p+1} \ar[u]^{\R\pi_*\gamma_\pi}_{\wr} \ar[r] & \R\pi_*\B_p 
\ar[u]^{\R\pi_*\gamma_\pi}_{\wr} \ar[r] & \R\pi_*\A_p 
\ar[u]_{\wr}^{\R\pi_*\gamma_\pi} \ar[r]^-{\R\pi_*\chi'}
& \R\pi_*\A_{p+1}[1] \ar[u]_{\wr}^{\R\pi_*\gamma_\pi}
}
$$
where $\chi\colon F^p/F^{p+1} \to F^{p+1}/F^{p+2}[1]$ is the map
in $\Dqc(\bbP)$ arising from the standard triangle associated to
the exact sequence $E_t$ and $\chi\colon \A_p \to \A_{p+1}[1]$ comes
from a similar process involving $E_b$. We point out that the top
triangle is the standard triangle associated to the exact sequence
\eqref{cplx:connE1}.

As a map of graded $\cO_Y$--modules $\beta$ is defined in degree $p$
as the composite
\begin{equation*}
\begin{split}
\E^p=\Hr^p(\Phi^p/\Phi^{p+1}) & \xrightarrow{\Hr^p(\R\pi_*\gamma_\pi)^{-1}}
{\mathrm{R}}^d\pi_*(\pi^*\F^p\otimes\omega_\pi) \\
& \xrightarrow{\phantom{XXX}\Iso\phantom{XX}\,} 
\,\F^p\otimes{\mathrm{R}}^d\pi_*\omega_\pi.
\end{split}
\end{equation*}
We argue as follows to show that $\beta$ is a map of complexes.
First, it is easy to see that 
$$
\Hr^p(\R\pi_*\chi')\colon \Rr^d\pi_*(\F^p\otimes\omega_\pi) \to
\Rr^d\pi_*(\F^{p+1}\otimes\omega_\pi)
$$
is the negative of the map induced by $\pi^*\partial_\F^p\otimes 1\colon
\pi^*\F^p\otimes\omega_\pi \to \pi^*\F^{p+1}\otimes\omega_\pi$. 
On the other hand $\Hr^p(\R\pi_*\chi)$ is the negative of the connecting map
$\E^p=\Hr^p(\Phi^p/\Phi^{p+1}) \to \Hr^{p+1}(\Phi^{p+1}/\Phi^{p+2})=\E^{p+1}$.
It is now clear from the definition of $\partial^\bullet_\E$ that $\beta$
is a map---in fact an isomorphism---of complexes.

It remains to show that the diagram \eqref{diag:alpha-beta} commutes.
To do so, it is enough to assume that {\em $Y$ is affine} and
and
$$
\F = i_y \Gamma_y\F
$$
for some $y\in Y$. The last assumption is equivalent to the assumption that
$$
\F = i_y\F(y)[-p]
$$
where $p=\De(y)$. Let
$$
M=\F(y), \qquad \G=i_yM
$$
and let $x$ be a point on $\bbP$ which is closed in the fiber of $\pi$
over $y$. Let $R$ be the completion of $\cO_Y$ at $y$, $S$
the completion of $\cO_\bbP$ at $x$, and $\varphi\colon R\to S$ the
map induced by $\pi$. We have a map
$$
{\bar{\alpha}}(x): (\sh{\pi}\F)(x) \longrightarrow \Gamma(\bbP,\,\E^p)
$$
defined as the composite
$$
(\sh{\pi}\F)(x) \xrightarrow{\text{natural}} \Gamma(\bbP,\,\A^{p,0})
\xrightarrow{\eqref{eq:A-phi}} \Gamma(\bbP,\,\E^p).
$$
In view of \Pref{prop:residue-thm}, it is enough to prove that
the dotted arrow in the diagram
\begin{equation}\label{diag:dotted}
\xymatrix{
\varphi_{\sharp}M \ar@{=}[r] \ar@{.>}[d] & \sh{\pi}\F(x) 
\ar[r]^{{\bar\alpha}(x)} & \Gamma(\bbP,\,\E^p) \ar[d]^{\Gamma(\bbP,\,\beta)} \\
\Hr_{\fm_S}^{d}(M\otimes\omega_{S/R}) \ar[r]^{\Iso} & 
\Hr_{x}^{d}(\pi^*\G\otimes\omega_\pi) \ar[r] & 
\Hr^d(\bbP,\,\pi^*\G\otimes\omega_\pi)
}
\end{equation}
can be filled by \eqref{iso:sharp-h} to make it commute.
To that end consider the (possibly different) arrow $u:\varphi_{\sharp}M\iso
\Hr_{\fm_S}^{d}(M\otimes\omega_{S/R})$ described as follows. First
we have a composite of maps of complexes
\begin{equation*}
\begin{split}
\Gamma_x\sh{\pi}\F \iso \R\Gamma_x\sh{\pi}\F 
\xrightarrow[Q_\bbP(\gamma_\pi)^{-1}]{\Iso} &
 \R\Gamma_x(\G[-p]\otimes(\omega_\pi[d])) \\
=\qquad &  \R\Gamma_x((\G\otimes\omega_\pi)[d-p])\\
 \iso\quad & \R\Gamma_{\fm_S}((M\otimes\omega_{S/R})[d-p]).
\end{split}
\end{equation*}
The map $u$ is defined as the $p$-th cohomology of this composite.
In view of the definitions of the maps ${\bar{\alpha}}(x)$ and $\beta$,
it is clear that \eqref{diag:dotted} commutes when the dotted arrow
is filled by $u$. Therefore we need to show that
$$
u=\eqref{iso:sharp-h}.
$$
To that end let 
$$
u'\colon {\mathbb{E}_\pi(\F)}(x) \iso \Hr_{\fm_S}^{d}(M\otimes\omega_{S/R})
$$
be defined in a manner completely analogous to the way $u$ was defined;
i.e., $u'$ is the $p$-th cohomology of the composite
\begin{equation*}
\begin{split}
\Gamma_x{\mathbb{E}_\pi(\F)} \iso \R\Gamma_x{\mathbb{E}_\pi(\F)}
 \xrightarrow[Q_\bbP(\eta_\pi)^{-1}]{\Iso} 
& \R\Gamma_x(\G[-p]\otimes(\omega_\pi[d])) \\
 = \qquad & \R\Gamma_x((\G\otimes\omega_\pi)[d-p])\\
\iso \quad & \R\Gamma_{\fm_S}((M\otimes\omega_{S/R})[d-p]).
\end{split}
\end{equation*}
Clearly, by definition of $\eta_\pi$, $\gamma_\pi$ and $L$ we have
$$
u'= u\circ L(x).
$$
Now according to the discussion in \cite[Remark\,10.1.10]{lns} we have
$$
u'= (-1)^{\varepsilon(x)}\eqref{iso:sharp-h}\circ L(x)
$$
where $\varepsilon(x)=\De'(x)\De(y)+\De'(x)$. Since $x$ is closed in its fiber 
over $y$ (so that $\De'(x)=\De(y)$) it is clear that 
$\varepsilon(x)$ is even. It follows that $u=\eqref{iso:sharp-h}$.
\qed

\subsection{The Trace Theorem} Now that we have the Trace Theorem
for relative projective space, the road to the general Trace
Theorem is clear. We begin with a lemma which helps us
reduce the problem to proper maps between ordinary schemes.



\begin{lem}\label{lem:reduce} Let $f\colon (\X,\,\De')\to (\Y,\De)$
be a pseudo-proper map in $\bbFc$ and let $\F \in \Cozs{\De}(\Y)$.
Let $\J\in\cO_\X$ be a coherent ideal contained in the maximal
defining ideal of $\X$. For each positive integer $n$,
let $v_n\colon \X_n\hookrightarrow \X$ be
the closed immersion defined by $\J^n$ and let $f_n\colon \X_n\to \Y$
be the composite $f_n=fv_n$. Then $\Tr{f}(\F)$ is a map of complexes
if and only if $\Tr{f_n}(\F)$ is a map of complexes for every $n$.
\end{lem}

\proof  Let ${\mathcal{C}}_n = \sHomb_\X({v_n}_*\cO_{\X_n},\,\sh{f}\F)$.
We regard ${\mathcal{C}}_n$ as a subcomplex of $\sh{f}\F$. Note that
$\sh{f}\F$ is the union (direct limit) of $\{{\mathcal{C}}_n\}$.
If $\Tr{f_n}(\F)$ is a map of complexes, by \Lref{lem:trThm-climm} and 
\eqref{diag:tr-tr-f}, 
$$
\Tr{f}(\F)|_{{\mathcal{C}}_n}\colon {\mathcal{C}}_n \to \F
$$
is a map of complexes. Taking direct limits over $n$ we conclude that
$\Tr{f}(\F)$ is a map of complexes. Conversely, if $\Tr{f}(\F)$ is
a map of complexes, then \eqref{diag:tr-tr-f} and \Lref{lem:trThm-climm}
give the result.
\qed

\begin{thm}\label{thm:tr-thm} Suppose $f\colon (\X,\,\De') \to (\Y,\,\De)$
is a pseudo-proper map in $\bbFc$ and $\F$ is an object of $\Cozs{\De}(\Y)$.
\begin{itemize}
\item[(a)] {\em (Trace Theorem)} The map 
$$
\Tr{f}(\F)\colon f_*\sh{f}\F \to \F
$$
is a map of complexes which is functorial in $\F$.
\item[(b)] {\em (Transitivity of Traces)} If 
$g\colon (\W,\,\De'')\to (\X,\,\De')$ is a second pseudo-proper
map then
$$
\xymatrix{
(fg)_*\sh{g}\sh{f}\F \ar@{=}[d] \ar[rr]^{\Iso}_{\sh{C}_{g,f}}
& & (fg)_*\sh{(fg)}\F \ar[dd]^{\Tr{fg}} \\
f_*g_*\sh{g}\sh{f}\F \ar[d]_{f_*\Tr{g}} & & \\
f_*\sh{f}\F \ar[rr]_{\Tr{f}} & & \F
}
$$
commutes.
\end{itemize}
\end{thm}

\proof Part\,(b) follows from \eqref{diag:tr-tr-f}. 
For Part\,(a), if $\J$ is a defining ideal of $\X$, then the maps $f_n$ of
\Lref{lem:reduce} are adic, and hence by \Lref{lem:reduce}
we are reduced to the case where $f$ is proper. 

Next consider a defining ideal $\I$ of $\Y$ and set $\J\set \I\cO_\X$.
By our assumption $\J$ is a defining ideal for $\X$. Let $u_n\colon
\Y_n \to Y$ be  the closed subscheme of $\Y$ defined by $\I^n$ and 
$v_n\colon \X_n\to \X$ the closed subscheme of $\J$ defined by $\J^n$. 
Let $h_n$ be the induced map from $\X_n$ to $\Y_n$. If $\Tr{h_n}(\sh{u_n}\F)$
is a map of complexes for each $n$, then by \Lref{lem:trThm-climm}
applied to $u_n$ and by \eqref{diag:tr-tr-f} we see that the Trace
Theorem is true for $f_n=fv_n$ and hence by \Lref{lem:reduce}
it is true for $f$. Now $h_n$ is a proper map of ordinary schemes.
Thus we are reduced to proving the Trace Theorem when $f$ is
a proper map of ordinary schemes in $\bbF$. Since we have proved
the Trace Theorem for relative projective space, the proof of
\cite[p.369,\,Theorem\,2.1]{RD} applies {\em mutatis mutandis} and
we are done.
\qed


\section{\bf The twisted inverse image pseudofunctor}\label{s:suresh}
In this section we summarize Nayak's work on the
pseudofunctor $\shr{\hd}$ on $\bbG$---laid out in \cite{suresh}---which 
extends the definition of $\shr{\hd}$ given by Alonso, Jerem{\'\i}as and 
Lipman in \cite{dfs} for pseudo-proper maps. Here are some
heuristics to orient the reader. Notionally ``upper shriek"
is to be thought of as a right adjoint to the (perhaps
non existent) derived ``direct image with proper supports on torsion
sheaves". This heuristic would dictate, for example, that we set $f^!$ 
equal to $\R\iGp{\X}f^*=f^*\R\iGp{\Y}$ for an 
{\it open immersion} $f:\X\to \Y$ of noetherian
formal schemes---for the corresponding direct image with
proper supports ``is" the functor ``extension by zero". Since the
objects of $\bbG$ are noetherian formal schemes, Deligne's
definition of direct image with proper support (see his
appendix to \cite{RD}) does not apply principally because
we cannot guarantee that a coherent ideal sheaf $\I\subset \cO_\U$
on an open subscheme $\U$ of a noetherian formal scheme
$\V$ extends to a coherent ideal sheaf $\J\subset \cO_\V$ on $\V$.
This is, incidentally, the reason why Verdier's proof of the
localness of $\shr{\hd}$ \cite[p.\,395,\,Corollary\,1]{f!}
does not apply mutatis mutandis to our situation\footnote{The
proof given in \cite[Proposition\,8.3.1]{dfs} is incorrect.
See however a partial resolution to the localness question
in \cite{dfscorr}.}. Nayak establishes the localness result
and pseudofunctoriality of $\shr{\hd}$. Here then is the 
promised summary.

\subsection{Factorizations} For $f:\X\to\Y$ in $\bbG$, we wish 
to define $f^!:\wDqcp(\Y)\to \Dqct^+(\X)$. If $f$ is pseudo-proper, 
$f^!$ is defined to be the functor $f_{\mathbf{t}}^\times|\wDqcp(\Y)$, 
where $f_{\mathbf{t}}^\times$ is right adjoint to 
$\R{f_*}:\Dqct(\X)\to \D(\Y)$---guaranteed to exist by 
\cite[Theorem\,2(a)]{dfs} and \cite[Theorem\,6.1(a)]{dfs} (see 
also the beginning of 1.3 of {\it Ibid} as well as the notation 
before Theorem\,7.4 in {\it Ibid}). If $f$ is an open immersion, then in 
light of the above discussion, we set 
$$
f^!\set \R\iGp{\X}f^* (=\R f^*\iGp{\Y}).
$$
If $f$ is a general map in $\bbG$, then $f$ can be written as a composite
of open immersions and pseudo-proper maps (in any order), and
the above gives us a clue on how to proceed. The problems that need
to be addressed are compatibilities of different definitions
corresponding to different factorizations of the same map.
Nayak's Theorem\,7.1.3 in \cite{suresh} as well as the remark
in 7.2.4 of {\em Ibid} answers this in a very satisfactory way.
Since the issues involved impact this paper we give a quick
thumbnail sketch of what is involved.

A sequence $\fc = (f_1,\ldots,f_n)$ of maps in $\bbG$
is called a {\it factorization} if $f_1\circ\ldots \circ f_n$
is defined (the maps are composable) and an individual
$f_i$ in the sequence is either pseudo-proper or an
open immersion. If $f=f_1\circ\ldots\circ f_n$, then $\fc$ is
called a {\it factorization of $f$}. We often write
$\lvert\fc\rvert = f_1\circ\ldots\circ f_n$. We define
$$
\shr{\fc}\set f_n^!\ldots f_1^! .
$$
If $\fc=(f_1,\ldots,f_n)$ and $\gc=(g_1,\ldots,g_m)$ are factorizations
such that $f_n\circ g_1$ is defined then we set
$$
\fc*\gc \set (f_1,\ldots,f_n,g_1,\ldots,g_m).
$$
Note that $\fc*\gc$ is a factorization of 
$\lvert\fc\rvert\circ\lvert\gc\rvert$ and that 
$$
\shr{(\fc*\gc)}= \shr{\gc}\shr{\fc}.
$$ 
In \cite{suresh} Nayak proves
that if $\fca$ and $\fcb$ are two factorizations of $f$, then
there is an isomorphism
$$
\bmu(\fca,\fcb) : \shr{\fcb} \iso \shr{\fca}
$$
satisfying the following properties.
\begin{itemize}
\item[(a)] $\bmu(\fca,\fca)={\boldsymbol{1}}_{\shr{\fca}}$.
\item[(b)] If $\fcc$ is a third factorization of $f$, then
$$
\bmu(\fca,\fcb)\circ\bmu(\fcb,\fcc) = \bmu(\fca,\fcc) .
$$
\item[(c)] If $\fc$ is a factorization of $f$, $\gca$, $\gcb$ of $g$ and
$\hc$ of $h$ where $g$ and $h$ are
such that $f\circ g\circ h$ is defined, then setting
$\eca= \fc*\gca*\hc$ and $\ecb=\fc*\gcb*\hc$, the
isomorphism $\bmu(\eca,\ecb)$ is induced by $\bmu(\gca,\gcb)$, i.e.
$$
\bmu(\eca,\ecb)= \shr{\hc}(\bmu(\gca,\gcb)(\shr{\fc})).
$$
\item[(d)] If
$$
\xymatrix{
\V \ar[r]^{v} \ar[d]_{g} & \X \ar[d]^{f} \\
\U \ar[r]_{u} & \Y
}
$$
is cartesian with $f$ pseudo-proper and $u$ an open immersion, then the 
base change map $\beta:v^*\fs \iso \gs\iGp{\U}u^*=\gs u^*$ of 
\cite[Theorem\,7.4]{dfs} (note that $u$ is adic) is given by
$$
\beta = \bmu(\fca,\fcb)
$$
where $\fca=(u,g)$ and $\fcb = (f,v)$.
\end{itemize}

Now for each map $f$ in $\bbG$ pick a distinguished factorization
$\fb$ of $f$---with the understanding that if $f$ is either an
open immersion or pseudo-proper $\fb=(f)$. Nayak then defines
$\fs$ as
$$
\fs \set \shr{\fb}.
$$
In fact Theorem\,7.1.3 and 7.2.4 of \cite{suresh} show that 
$\shr{\hd}$ is a pseudofunctor\footnote{More
precisely a {\it pre-pseudofunctor} in the sense of \ref{ss:conventions}\,(11).
However if for each $\X \in \bbG$ we set 
$\X^!=\Dqct^+(\X)$ rather than $\wDqcp(\X)$, then $\shr{\hd}$ is a 
pseudofunctor.}, the required isomorphism
$$
C^!_{f,g}: \fs\gs \iso (gf)^!,
$$ 
being (with $h=gf$)
$$
C^!_{f,g} \set \bmu(\hb,\gb*\fb).
$$
For a factorization $\fc$ of $f$, we have a canonical isomorphism
$$
\bnu{F} : \shr{\fc} \iso \fs
$$
given by $\bnu{F}=\bmu(\fb,\,\fc)$.

\subsection{Flat base change} Suppose we have a 
cartesian square of noetherian formal schemes
$$
\xymatrix{
\V \ar[r]^{v} \ar[d]_{g} & \X \ar[d]^{f} \\
\U \ar[r]_{u} & \Y
}
$$
with $u$ {\it flat}. The flat base change theorem of Alonso,
Jerem{\'\i}as and Lipman \cite[Theorem\,7.4]{dfs} states that
if $f$ is pseudo-proper there is a base change isomorphism
(between functors on $\wDqcp(\Y)$)
$$
\beta=\beta(f,u) :\R\iGp{\V}v^*\fs \iso \gs\R\iGp{\U}u^*.
$$
Nayak extends $\beta$ to the case when $f$ is a composite 
of compactifiable maps. In greater detail, if $f$ is an open
immersion, let $\beta(f,u)$ be the obvious isomorphism. For a
general morphism $f$ in $\bbG$, suppose $\fc$ is a factorization
of $f$. Then $\fc$ induces, via the base change map $u$,
a factorization $\U\times_\Y\fc$ of $g$. Setting $\gc=\U\times_\Y\fc$
one gets, by applying $\beta$ successively to the maps in the
sequence $\fc$, an isomorphism
$$
\beta(\fc):\R\iGp{\V}v^*\shr{\fc} \iso \shr{\gc}\R\iGp{\U}u^*.
$$
This gives us an isomorphism
\stepcounter{thm}
\begin{equation*}\label{eq:nay-bc}\tag{\thethm}
\beta=\beta(f,u):\R\iGp{\V}v^*\fs \iso \gs\R\iGp{\U}u^*
\end{equation*}
given by
$$
\beta = \bmu(\gb,\,\U\times_\Y\fb)(\R\iGp{\U}u^*)\circ \beta(\fb).
$$
Nayak shows that the map 
$\beta$ is ``independent" of the chosen 
distinguished factorization $\fb$ of $f$ in the following sense.
If $\fc$ is any factorization of $f$ and $\gc=\U\times_\Y\fc$, then
the diagram
\stepcounter{thm}
\begin{equation*}\label{diag:nay-bc}\tag{\thethm}
\xymatrix{
\R\iGp{\V}v^*\shr{\fc} \ar[d]_{\beta(\fc)} \ar[rr]^{\Iso}_{\bnu{F}} & &
\R\iGp{\V}v^*\fs \ar[d]^{\beta} \\
\shr{\gc}\R\iGp{\U}u^* \ar[rr]^{\Iso}_{\bnu{G}} & & \gs\R\iGp{\U}u^*
}
\end{equation*}
commutes. (Cf. \cite[\S\,2.3 and 7.2.3]{suresh}.)

\subsection{Comparing pseudofunctors}\label{ss:gamma-ext} The
construction of $\shr{\boldsymbol{-}}$ on $\bbG$ is built on
its description on open immersions, pseudo-proper maps and
the fact that it satisfies a compatibility relationship for
open base changes of pseudo-proper maps. This allows one
to compare a pseudofunctor on a subcategory of
$\bbG$ with $\shr{\boldsymbol{-}}$
if we have a comparison for pseudo-proper maps and for open
immersions, and certain compatibility relations. We give more
details below.

Let $\lbbG$ be  a full subcategory of $\bbG$ such that if $\Y\in\lbbG$
and $f\colon \X\to \Y$ is a map in $\bbG$, then $f$ is a map
in $\lbbG$ (i.e. $\X\in \lbbG$). Note that all open subschemes
and all closed subschemes of $\Y$ are in $\lbbG$, and if
$g\colon \W\to \Y$ is another map in $\lbbG$, then $\X\times_\Y\W
\in \lbbG$. Note also that $\cbbF$ and $\cbbF\cap\rbbF$ provide
examples of $\lbbG$.

Let ${\mathbf{P}}$ be the class of pseudo-proper maps in $\lbbG$
and ${\mathbf{F}}$ the class of open immersions in $\lbbG$. Suppose
we have another pseudofunctor $\natb{\boldsymbol{-}}$ on $\lbbG$ such
that for each object $\X$ in $\lbbG$ we have a functor
$$
S_\X\colon \nat{\X} \longrightarrow \X^!
$$
and for each map $f:\X\to \Y$ in ${\mathbf{P}}\cup{\mathbf{F}}$ a
functorial map
$$
\gamma_f \colon S_\X\nat{f} \longrightarrow \fs S_\Y
$$
such that 
\begin{itemize}
\item[(a)] If $f \in {\mathbf{F}}$, $\gamma_f$ is an isomorphism.
\item[(b)] If $\X\xrightarrow{f} \Y\xrightarrow{g} \Z$ are a 
pair of maps such that either both $f$ and $g$ are in ${\mathbf{P}}$
or both $f$ and $g$ are in ${\mathbf{F}}$, then the diagram
\stepcounter{thm}
\begin{equation*}\label{diag:gamma-ext-1}\tag{\thethm}
\xymatrix{
S_\X\nat{f}\nat{g} \ar[d]_{\gamma_f} \ar[rr]^{C^{\natural}_{f,g}} & &
S_\X\nat{(gf)} \ar[dd]^{\gamma_{gf}} \\
\fs S_\Y\nat{g} \ar[d]_{\fs\gamma_g} & & \\
\fs\gs S_\Z \ar[rr]_{C^!_{f,g}} & & (gf)^!S_\Z
}
\end{equation*}
commutes.
\item[(c)] If
\stepcounter{thm}
\begin{equation*}\label{diag:gamma-ext-2}\tag{\thethm}
\xymatrix{
\V \ar[r]^v \ar[d]_{g} & \X \ar[d]^{f} \\
\U \ar[r]_u & \Y
}
\end{equation*}
is cartesian with $u\in {\mathbf{F}}$ and $f\in {\mathbf{P}}$, then
the following diagram commutes
\stepcounter{thm}
\begin{equation*}\label{diag:gamma-ext-3}\tag{\thethm}
\xymatrix{
S_\V\nat{g}\nat{u} \ar[d]_{\nat{\varphi}}^{\wr} \ar[rr]^{\gamma_g} & &
g^!S_\U\nat{u} \ar[d]^{\gamma_u}_{\wr} \\
S_\V\nat{v}\nat{f} \ar[d]_{\gamma_v}^{\wr} & & 
g^!u^!S_\Y \ar[d]^{\varphi^!}_{\wr} \\
v^!S_\X\nat{f} \ar[rr]_{v^!\gamma_f} & & v^!f^!S_\Y
}
\end{equation*}
where $\nat{\varphi}=(C^{\natural}_{v,f})^{-1}C^{\natural}_{g,u}$ and
$\varphi^!=(C^!_{v,f})^{-1}C^!_{g,u}$.
\end{itemize}

In \cite[Theorem\,7.2.5]{suresh}, Nayak proves
\begin{thm}\label{thm:gamma-ext}{\em (Nayak)}
Under the above hypotheses on $(\natb{\boldsymbol{-}},\,S,\,\gamma)$, 
the association $f\mapsto \gamma_f$, $f\in {{\mathbf{P}}\cup{\mathbf{F}}}$,
can be extended {\em uniquely} to all morphisms 
$$
f\colon \X \to \Y
$$ 
in $\lbbG$ in such a way that
\begin{itemize}
\item[(a)] If $\Y\xrightarrow{g}\Z$ is a second map 
in $\lbbG$, then the diagram \eqref{diag:gamma-ext-1} commutes.
\item[(b)] The natural transformation $\gamma_f$ is compatible with open
immersions into $\Y$, i.e., whenever we have a cartesian square 
\eqref{diag:gamma-ext-2}
with $u\in{\mathbf{F}}$ and $f$ {\it any} map in $\lbbG$ the diagram 
\eqref{diag:gamma-ext-3} commutes.
\end{itemize}
\end{thm}

Note that the uniqueness assertion is clear. Indeed if $\fc=(f_1,\dots,f_n)$
is any factorization of a map $f:\X\to \Y$ in $\lbbG$, and 
$$
\natb{\fc} \set \nat{f_1}\dots\nat{f_n}
$$
then we have a map
$$
{\boldsymbol{\gamma}}(\fc)\colon S_\X\natb{\fc} \to \shr{\fc}S_\Y
$$
given by the composite (where we are suppressing the functors $S_\X$ etc.)
\begin{align*}
\natb{\fc} = \nat{f_n}\ldots\nat{f_1} & 
\xrightarrow{\,\phantom{XX}\gamma_{f_n}\phantom{XX}}
       f_n^!\nat{f_{n-1}}\ldots\nat{f_1} \\
& \xrightarrow{\phantom{X}f_n^!\gamma_{f_{n-1}}\phantom{X}} 
f_n^!f_{n-1}^!\nat{f_{n-2}}\ldots\nat{f_1}\\
& \dots \\
& \xrightarrow{\phantom{X}f_n^!\ldots f_2^!\gamma_{f_1}\,} f_n^!\ldots f_1^! = 
\shr{\fc},
\end{align*}
and $\gamma_f$ is characterized by the commutativity of
\stepcounter{thm}
\begin{equation*}\label{diag:gamma-ext-4}\tag{\thethm}
\xymatrix{
S_\X\natb{\fc} \ar[d]_{\wr} \ar[rr]^{{\boldsymbol{\gamma}}(\fc)} & &
\shr{\fc}S_\Y \ar[d]^{\wr} \\
S_\X\nat{f} \ar[rr]_{\gamma_f} & & f^!S_\Y
}
\end{equation*}
where the downward arrows are the natural ones obtained from
$\natb{\boldsymbol{-}}$ and $\shr{\boldsymbol{-}}$.


\section{\bf The comparison map}
In this section we define the fundamental comparison map 
which serves as the bridge between the
pseudofunctors $\sha{\hd}$ and $\shr{\hd}$ (see \Dref{def:compare}).
This comparison is easily defined when $f$ is pseudo-proper, but
has a more involved definition for a general $f$ in $\cbbFc$.

\subsection{Pseudo-proper maps} Suppose 
$f:(\X,\De')\to (\Y,\De)$ is a pseudo-proper map in $\bbFc$. 
For $\F\in\Coz(\Y)$, the trace map
$$
\Tr{f}(\F): f_*\sh{f}\F \longrightarrow \F
$$
gives, on applying $Q_\Y$, a map
$$
Q_\Y(\Tr{f}(\F)) : \R{f_*}Q_\X(\sh{f}\F) \to Q_\Y(\F)
$$
whence we have, by the universal property of $(\fs,\ttr{f})$,
a unique map $\ga{f}(\F): Q_\X\sh{f}\F \to \fs Q_\Y\F$ such
that $Q_\Y\Tr{f}(\F)=\ttr{f}\R f_*(\ga{f}(\F))$. This
map is functorial in $\F$ giving a {\it comparison map} of
functors
\stepcounter{thm}
\begin{equation*}\label{eq:compare}\tag{\thethm}
\ga{f}: Q_\X\sh{f} \longrightarrow \fs Q_\Y.
\end{equation*}
Note that the source of the functors (being compared by $\ga{f}$)
is $\Cozs{\De}(\Y)$ and the target is $\Dqct^+(\X)$. We would
like to extend the definition of $\ga{f}$ to the case where
$f\in \cbbFc$ (i.e. $f$ is a composite of pseudo-proper maps and
open immersions in any order). We will, from now on, often
suppress the localization functors $Q_\X$, $Q_\Y$ etc..

We begin with two Lemmas. Let $f$ be as above (i.e. pseudo-proper
and in $\bbFc$). Let $u:\U\to \Y$ be an open immersion,
$\V=\U\times_\Y\X$, $v:\V\to \X$ and $g:\V\to \U$ the two projections.
We have---from the local nature of direct images---$u^*\R{f_*}=\R{g_*}v^*$. 
The sheafified version of duality on
formal schemes \cite[Theorem\,8.2]{dfs} gives us
$v^*\fs\F=\gs u^*\F$ and $u^*\ttr{f}(\F)=\ttr{g}(u^*\F)$
for $\F\in\wDqcp(\Y)$. We also have, by the construction of
$\sh{f}$ and $\sh{g}$, $v^*\sh{f}=\sh{g}u^*$. 

\begin{lem}\label{lem:openimm} In the above situation the diagram
$$
\xymatrix{
v^*\sh{f}\F \ar[r]^-{v^*\ga{f}} \ar@{=}[d] & v^*\fs\F 
\ar@{=}[d]\\
\sh{g}{u^*\F} \ar[r]_{\ga{g}} & \gs{u^*\F}
}
$$
commutes for $\F\in\Cozs{\De}(\Y)$.
\end{lem}
\proof From the punctual nature of the definition of $\Tr{f}$ it
is clear that $\Tr{f}$ behaves well with respect to open immersions
into $\Y$, i.e. $u^*\Tr{f}(\F)=\Tr{g}(u^*\F)$. The Lemma is a simple
consequence of this. (See also \Pref{prop:basech} where a more
general result is proved somewhat elaborately.)
\qed

\begin{lem}\label{lem:ga-trans} Let $(\X,\,\De_1)\xrightarrow{f}
(\Y,\,\De_2)\xrightarrow{g} (\Z,\,\De_3)$ be a pair of pseudo-proper maps
in $\bbFc$. Then the diagram
$$
\xymatrix{
Q_\X\sh{f}\sh{g} \ar[d]_{\ga{f}(\sh{g})} \ar[r]^-{C^{\sharp}_{f,g}}
& Q_\X\sh{(gf)} \ar[dd]^{\ga{gf}} \\
\fs Q_\Y\sh{g} \ar[d]_{f^!(\ga{g})} & \\
\fs\gs Q_{\Z} \ar[r]^{C^!_{f,g}} & (gf)^!Q_{\Z}
}
$$
commutes.
\end{lem}
\proof We will suppress the symbols $Q_\X$ and $Q_\Y$. By the universal
property of $((gf)^!,\,\ttr{gf})$ it is enough to show
\begin{equation}\label{eq:ga-trans}
\ttr{gf}\circ\R(gf)_*(C^!_{g,f}\circ\fs(\ga{g})\circ\ga{f}(\sh{g}))
=\ttr{gf}\circ\R(gf)_*(\ga{gf}\circ C^{\sharp}_{f,g})
\end{equation}
In what follows, we save on notation by writing $f_*$, $g_*$ and
$(gf)_*$ for $\R{f_*}$, $\R{g_*}$ and $\R(gf)_*$. 

Consider the following diagram
$$
\xymatrix{
g_*(f_*\fs)\sh{g} \ar@{=}[ddd] \ar[rr]^{\ttr{f}}
& & g_*\sh{g} \ar@{=}[d] \ar@{}[dr]|{\square_3} \ar@{}[dl]|{\square_2}
\ar[r]^{\ga{g}}
& g_*\gs \ar[r]^{\ttr{g}} & Q_\Y \ar@{=}[dl] \\
& g_*(f_*\sh{f})\sh{g} \ar[ul]_{\ga{f}}\ar[r]_{\Tr{f}} \ar@{}[dl]|{\square_1} 
\ar@{}[drr]|{\square} \ar@{=}[d] & g_*\sh{g} \ar[r]_{\Tr{g}}
& Q_\Y & \\
& (gf)_*\sh{f}\sh{g} \ar[dl]_{\ga{f}} \ar[rr]_{C^{\sharp}_{f,g}}
& & (gf)_*\sh{(gf)} \ar[u]_{\Tr{gf}} \ar[dr]^{\ga{gf}}
\ar@{}[ur]|{\square_4}   &\\
(gf)_*\fs\sh{g} \ar[rr]^{\ga{g}} & & (gf)_*\fs\gs \ar[rr]^{C^!_{f,g}}
& & (gf)_*(gf)^! \ar[uuu]_{\ttr{gf}}
}
$$
If we can show that $\square_1$, $\square_2$, $\square_3$, $\square_4$,
$\square$ and the outer rectangle commute then \eqref{eq:ga-trans}
follows easily. Now $\square_1,\ldots,\square_4$ and $\square$
clearly commute ($\square$ commutes by the transitivity of the
trace). It remains to show that the outer rectangle
commutes. To that end consider the diagram
$$
\xymatrix{
g_*(f_*f^!)\sh{g} \ar@{=}[dd] \ar[dr]^{\ga{g}} \ar[r]^{\ttr{f}} &
g_*\sh{g} \ar[r]^{\ga{g}} & g_*\gs \ar@{=}[d] \ar[r]^{\ttr{g}} &
Q_\Y \\
& g_*f_*\fs\gs \ar@{=}[d] \ar[r]^{\ttr{f}} & g_*\gs & \\
(gf)_*\fs\sh{g} \ar[r]^{\ga{g}} & (gf)_*\fs\gs \ar[rr]^{C^!_{f,g}}
& & (gf)_*(gf)^! \ar[uu]_{\ttr{gf}}
}
$$
The subdiagram on the left obviously commutes. The one on top
commutes by the functoriality of $\ga{g}$ and $\ttr{f}$. The
remaining (L-shaped) subdiagram commutes by the definition of
$C^!_{f,g}$. 
\qed

We are now in a position to prove

\begin{thm}\label{thm:compare} There exists a unique family of natural
transformations
\begin{equation}\label{def:compare}
\ga{f}\colon Q_\X\sh{f} \longrightarrow \fs{Q_\Y},
\end{equation}
one for each morphism $f\colon (\X,\,\De') \to (\Y,\,\De)$ in $\cbbFc$
such that
\begin{itemize}
\item[(a)] If $f$ is pseudo-proper $\ga{f}$ is given by \eqref{eq:compare}.
\item[(b)] If $f$ is an open immersion then $\ga{f}$ is the identity
transformation on $f^*$. 
\item[(c)] The map $\ga{f}$ is compatible with open immersions into $\Y$, i.e.
\Lref{lem:openimm} holds with the weaker hypothesis that $f\in \cbbFc$.
\item[(d)] \Lref{lem:ga-trans} holds under the weaker hypothesis that
$f$ and $g$ are in $\cbbFc$.
\end{itemize}
\end{thm}

\proof Let the triple $(\lbbG,\,\{S_\X\},\,\{\gamma_f\}_{{\mathbf{P}}\cup 
{\mathbf{F}}})$ in \Ssref{ss:gamma-ext} be given by the data
$(\cbbFc,\,\{Q_\X\},\,\{\ga{f}\}_{{\mathbf{P}}\cup{\mathbf{F}}})$.
By \Lref{lem:openimm} and \Lref{lem:ga-trans} we see that \Tref{thm:gamma-ext}
applies to our situation. 
\qed
\begin{rem}
{\em We remind the reader how the map $\ga{f}$ is computed in practice
for $f\colon \X \to \Y$  a general maps in $\cbbFc$.
First, if $f$ is an open immersion, then
recall that $f^!=f^*$ and $\sh{f}=f^*|\Coz(\Y)$, and so we set 
$\ga{f}$ equal to the identity map $f^*|_{\Coz(\Y)}\to f^*|_{\Coz(\Y)}$.
For general $f \in \cbbFc$ we proceed as follows:
if $\fc=(f_1,\ldots,f_n)$ is a factorization of $f$, we write
$\sha{\fc}=\sh{f_n}\ldots\sh{f_1}$. 
We define
$$
\gab(\fc):\sha{\fc} \longrightarrow \shr{\fc}
$$
as the composite 
\begin{align*}
\sha{\fc} = \sh{f_n}\ldots\sh{f_1} & \xrightarrow{\,\phantom{XX}\ga{f_n}\phantom{XX}}
       f_n^!\sh{f_{n-1}}\ldots\sh{f_1} \\
& \xrightarrow{\phantom{X}f_n^!\ga{f_{n-1}}\phantom{X}} 
f_n^!f_{n-1}^!\sh{f_{n-2}}\ldots\sh{f_1}\\
& \vdots \\
& \xrightarrow{\phantom{X}f_n^!\ldots f_2^!\ga{f_1}\,} f_n^!\ldots f_1^! = 
\shr{\fc}.
\end{align*}
Then
$$
\ga{f}\colon \sh{f} \to \fs
$$
is characterized by the commutativity of}
$$
\xymatrix{
\sha{\fc} \ar[d]_{\gab(\fc)} \ar[r]^-{\Iso} & \sh{f} \ar[d]^{\ga{f}} \\
\shr{\fc} \ar[r]^-{\Iso} & \fs
}
$$

\end{rem}

%
%



\subsection{{\'E}tale base change}\label{ss:basech}
Suppose we have a cartesian square \stepcounter{thm}
\begin{equation*}\label{diag:basech}\tag{\thethm}
\xymatrix{
\V \ar[r]^{v} \ar[d]_g & \X \ar[d]^f \\
\U \ar[r]^u & \Y
}
\end{equation*}
with $u$ (and hence $v$) {\it flat and adic}. If $\F\in \Dqct^+(\Y)$ and 
$\G\in \Dqct^+(\X)$, then $\R\iGp{\U}u^*\F=u^*\R\iGp{\Y}\F=u^*\F$
and $\R\iGp{\V}v^*\G=v^*\iGp{\X}\G=v^*\G$ (see 
\cite[Corollary\,5.2.11(c)]{dfs}). By \eqref{eq:nay-bc}
we have a base change isomorphism
$$
\R\iGp{\V}v^*\fs \iso \gs\R\iGp{\U}u^*,
$$ 
and in light of the
relations we have given (in our special case), this translates
to a functorial isomorphism
$$
\beta_\F:v^*\fs\F \iso \gs u^*\F \qquad (\F\in\Dqct^+(\Y)).
$$
If $f$ is pseudo-proper the map $\beta_\F$ is obtained (when, 
as we have assumed, $u$ is adic) in the following manner. By 
\cite[Proposition\,7.2(b)]{dfs} we have a functorial isomorphism
$$
\theta_\G: u^*\R f_*\G \iso \R g_*v^*\G \qquad (\G\in\D(\X))
$$
where the map $\theta_\G$ is adjoint to the composite
$$
\R f_*\G \to \R f_*\R v_*v^*\G \to \R u_*\R g_*v^*\G.
$$
The map $\beta_\G$ is then defined via the universal
property of $(g^!,\ttr{g})$ as the map adjoint to
the composite
$$
\R g_*v^*\fs\F \xrightarrow{\theta^{-1}_{\fs\F}} u^*\R{f_*}\fs\F
\xrightarrow{u^*\ttr{f}} u^*\F.
$$

The following proposition relates the above base change with
a base change for the pseudofunctor $\sha{\hd}$ in a special
situation. We treat the isomorphisms in \cite[Prop.\,10.3.6,\,(90),(91)]{lns}
as equalities.

\begin{prop}\label{prop:basech}
Suppose $f:(\X,\De')\to (\Y,\De)$ is a map in $\cbbFc$, and $u:\U\to \Y$ 
is an {\'e}tale adic map such that for $q\in\U$ the
natural map ${\widehat{\cO}}_{\Y,u(q)}\to {\widehat{\cO}}_{\U,q}$ is an
isomorphism (equivalently, the residue fields at $u(q)$ and
$q$ are isomorphic via the natural map between them). Consider the cartesian 
square \eqref{diag:basech} above. Then for $\F\in\Coz(\Y)$ the following 
diagram commutes
$$
\xymatrix{
v^*\sh{f}\F \ar[r]^{v^*\ga{f}} \ar@{=}[d] & v^*\fs\F 
\ar[d]^{\beta_\F}_{\simeq} \\
\sh{v}\sh{f}\F \ar[d]_{\beta'_\F}^{\simeq} & \gs u^*\F \ar@{=}[d] \\
\sh{g}\sh{u}\F \ar[r]_{\ga{g}} & \gs\sh{u}\F
}
$$
where $\beta'_\F: \sh{v}\sh{f}\F \iso \sh{g}\sh{f}\F$ is the
functorial isomorphism arising from the pseudofunctorial nature
of $\sha{\hd}$, and the equalities are the identifications given in
\cite[Prop.\,10.3.6,\,(90),(91)]{lns}.
\end{prop}

\proof  The proposition is trivially true if $f$ is an open immersion,
for, in this case, $\ga{f}:f^*\to f^*$ and $\ga{g}:g^*\to g^*$ are
the identity maps.

Suppose the proposition is true when $f$ is pseudo-proper. We claim
that then it is true for $f$ a composite of pseudo-proper and
open maps. To see this, fix $f$ in $\cbbFc$. Suppose $\fc=(f_1,\ldots,f_n)$
is a factorization of $f$. Recall that we require each $f_i$ to
be either pseudo-proper or an open immersion. The map $u:\U\to \V$
induces, via base change, a factorization $\gc=(g_1,\ldots,g_n)$
of $g$---the factorization $\U\times_\Y\fc$ in \Sref{s:suresh}---such 
that with $u=v_0$, $v=v_n$, $\U=\V_0$, $\V=\V_n$,
$\Y=\X_0$ and $\X=\X_n$ we have, for $k=1,\ldots,n$ a cartesian
square
$$
\xymatrix{
\V_k \ar[d]_{g_k} \ar[r]^{v_k} & \X_k \ar[d]^{f_k} \\
\V_{k-1} \ar[r]_{v_{k-1}} & \X_{k-1}
}
$$
It is easy to check that each $v_k$, $k=0,\ldots,n$ satisfies the
same hypotheses that $u$ satisfies, viz. $v_k$ is {\'e}tale, adic
and ${\widehat{\cO}}_{\V_k,p} \to {\widehat{\cO}}_{\X_k,v_k(p)}$
is an isomorphism for every $p\in \V_k$. Since the proposition
is true with $f=f_k$, therefore the diagram
\begin{equation}\label{eq:bold}
\xymatrix{
v^*\sha{\fc} \ar[d]_{\beta'(\fc)}^{\simeq} \ar[r]
& v^*\shr{\fc} \ar[d]_{\simeq}^{\beta(\fc)} \\
\sha{\gc}u^* \ar[r] & \shr{\gc}u^*
}
\end{equation}
commutes, where the two horizontal arrows are $v^*\gab(\fc)$ and
$\gab(\gc)u^*$, the map $\beta'(\fc):v^*\sha{\fc} = \sh{v}\sha{\fc} \iso
\sha{\gc}\sh{u}=\sha{\gc}u^*$ is the isomorphism that arises
from the pseudofunctoriality of $\sha{\hd}$ (indeed $\sh{v_k}\sh{f_k}\iso
\sh{g_k}\sh{v_{k-1}}$) and $\beta(\fc):v^*\shr{\fc}\iso \shr{\gc}u^*$
is the natural isomorphism arising from the various base change isomorphisms
$v_k^*f_k^! \iso g_k^!v_{k-1}^*$ (the obvious one if $f_k$ is an open
immersion, and the map $\beta$ in \cite[Theorem\,7.4]{dfs} when
$f_k$ is pseudo-proper). Consider the diagram (with the inner square
being \eqref{eq:bold} and the slanted arrows being the natural
isomorphisms arising from the pseudofunctorial nature of $\sha{\hd}$
and $\shr{\hd}$)
$$
\xymatrix{
v^*\sh{f} \ar[ddd]_{\beta'} \ar[rrr]^{v^*\ga{f}} 
& & & v^*\fs \ar[ddd]^{\beta} \\
& v^*\sha{\fc} \ar@{}[dl]|{\square_1} \ar@{}[ur]|{\square_2}
\ar[ul]^{\simeq} \ar[d]^{\simeq} \ar[r]
& v^*\shr{\fc} \ar@{}[dr]|{\square_3} \ar[ur]_{\simeq} \ar[d]_{\simeq} & 
\\
& \sha{\gc}u^* \ar@{}[dr]|{\square_4} \ar[dl]_{\simeq} \ar[r]
& \shr{\gc}u^* \ar[dr]^{\simeq} & \\
\sh{g}u^* \ar[rrr]_{\ga{g}u^*} &&& \gs u^*
}
$$
If we show that each of the subdiagrams $\square_1,\ldots,\square_4$ commutes, 
then, as the inner square (i.e. \eqref{eq:bold}) commutes, the outer square 
will commute and the proposition will follow.

The subdiagram $\square_1$ commutes by the pseudo functoriality of $\sha{\hd}$
(recall that $u^*|\Coz(\Y)=\sh{u}$ and $v^*|\Coz(\X)=\sh{v}$). 
The subdiagrams $\square_2$ and $\square_4$ commute by the definitions
of $\ga{f}$ and $\ga{g}$. Subdiagram $\square_3$ commutes by
\eqref{diag:nay-bc}.

It remains to prove the proposition when $f$ is pseudo-proper. By the
universal property of $(\gs,\ttr{g})$ it is enough to show
\begin{equation}\label{eq:tau-beta}
\ttr{g}(u^*\F)\circ\R{g_*}(\ga{g}(\sh{u}\F)\circ\beta'_\F)=
\ttr{g}(u^*\F)\circ\R{g_*}(\beta_\F\circ v^*\ga{f}(\F)).
\end{equation}
Part of the proof (commutativity of \eqref{diag:sh-theta} below)
rests on the fact that at the {\it punctual level} the punctual
trace $\Tr{g,x'}$ at a point $x'\in\V$ equals the punctual trace 
$\Tr{f,v(x')}$. We elaborate below.

Let $x\in\X$ and $\G=i_x\sh{f}\F(x)$. We will examine the map
$\theta_\G$ under the identifications $\R{g_*}v^*\G=g_*v^*\G$
and $u^*\R{f_*}\G=u^*f_*\G$ (note that $v^*\G=\sh{v}\G$ is also
flasque). It can be checked that $\theta_\G:u^*f_*\G\to g_*v^*\G$
is adjoint to $f_*\G\to f_*v_*v^*\G=u_*g_*v^*\G$.

Let $y=f(x)$, $A={\widehat{\cO}}_{\Y,y}$, $R={\widehat{\cO}}_{\X,x}$,
$M=\F(y)$ and $\varphi:A\to R$ the map induced by $f$. Let
$-_{\#}$ be the pseudofunctor on zero-dimensional modules over
local rings constructed by Huang in \cite{I-C}. Then $(\sh{f}\F)(x)
=\varphi_{\#}M$. If $v^{-1}(x)$ is empty (so that $u^{-1}(y)$ is
also empty) then $u^*f_*\G$ and $g_*v^*\G$ are both zero and
$\theta_\G$ is the obvious map. The interesting case is when
we have $p\in\V$ such that $x=v(p)$. Let $q=g(p)$. Using
our hypothesis on $u:\U\to\Y$ we make the identifications
${\widehat{\cO}}_{\U,q}=A$ and ${\widehat{\cO}}_{\V,p}=R$ 
(since the isomorphisms underlying these are canonical).
With these identifications, the map ${\widehat{\cO}}_{\U,q}
\to {\widehat{\cO}}_{\V,p}$ induced by $g$ is again
$\varphi:A\to R$.

Suppose $v^{-1}(x)=\{p_1,\ldots,p_n\}$. Set $q_i=g(p_i)$ for
$i=1,\ldots,n$. Then $u^{-1}(y)=\{q_1,\ldots,q_n\}$. We have
$f_*i_x\varphi_{\#}M=i_y\varphi_{\#}M$, $g_*i_{p_k}\varphi_{\#}M
=i_{q_k}\varphi_{\#}M$, $v_*i_{p_k}\varphi_{\#}M=i_x\varphi_{\#}M$,
$u_*i_{q_k}\varphi_{\#}M=i_y\varphi_{\#}M$, $v^*i_x\varphi_{\#}M
=\oplus_ki_{p_k}\varphi_{\#}M$ and finally $u^*i_y\varphi_{\#}M=
\oplus_ki_{q_k}\varphi_{\#}M$. Thus in every case (whether $v^{-1}(x)$
is empty or not) we have
$$
v^*f_*\G=g_*u^*\G.
$$
It is now easy to see that $\theta_\G:v^*f_*\G\to g_*u^*\G$ is the
identity map. Indeed, the only interesting case is when
$v^{-1}(x)$ is non-empty, and in this case, with previous notations, 
the map
$$
f_*\G \to f_*v_*v^*\G = u_*g_*v^*\G
$$
is the diagonal map $i_y\varphi_{\#}M\to \oplus_k i_y\varphi_{\#}M$.
This proves the assertion on $\theta_\G$.

Now suppose $p,x,q,y$ satisfy $x=v(p)$, $y=f(x)$ and $q=g(p)$. The
$R$-modules $(\sh{f}\F)(x)$, $(\sh{v}\sh{f}\F)(p)$ and $(\sh{g}\sh{u})\F(p)$
are all equal to $\varphi_{\#}M$ and the natural isomorphism 
$\beta'_\F(p): (\sh{v}\sh{f}\F)(p)\iso (\sh{g}\sh{u}\F)(p)$ is the
identity map on $\varphi_{\#}M$. If $x$ is closed in the fiber $f^{-1}(y)$
(equivalently $p$ is closed in the fiber $g^{-1}(q)$) then, from the
definitions of $\Tr{f,x}$ and $\Tr{g,p}$ we have
$$
\Tr{f,x}(\F) = \Tr{g,p}(\sh{u}\F) = \Tr{R/A,M} : \varphi_{\#}M \to M.
$$

The above assertions imply that the diagram
\begin{equation}\label{diag:sh-theta}
\xymatrix{
u^*f_*\sh{f}\F \ar[r]^-{\Iso}_{\theta_{\sh{f}\F}} \ar[dd]_{u^*\Tr{f}}
& g_*v^*{\sh{f}\F} \ar@{=}[d] \\
& g_*\sh{v}\sh{f}\F \ar[d]^{g_*\beta'_\F} \\
u^*\F=\sh{u}\F & g_*\sh{g}\sh{u}\F \ar[l]^-{\,\,\Tr{g}}
}
\end{equation}
commutes. 

Consider the following diagram:
$$
\xymatrix{
F' \ar[ddd]_{u^*\ttr{f}} \ar[rrr]^-{\Iso}_\theta & & & G'
\ar[ddd]^{\R{g_*}\beta_\F} \\
& F \ar[ul] \ar[d] \ar[r]^-{\Iso}_\theta
& G \ar[ur] \ar[d] & \\
& E \ar@{=}[dl] & H \ar[l] \ar[dr] & \\
E' &&& H' \ar[lll]^{\ttr{g}}
}
$$
where
\begin{align*}
E  & = u^*\F=\sh{u}\F \\
F  & = u^*f_*\sh{f}\F = u^*\R{f_*}\sh{f}\F \\
G  & = g_*v^*\sh{f}\F = g_*\sh{v}\sh{f}\F = \R{g_*}v^*\sh{f}\F \\
H  & = g_*\sh{g}\sh{u}\F = \R{g_*}\sh{g}\sh{u}\F \\
E' & = u^*\F=\sh{u}\F = E \\
F' & = u^*\R{f_*}\fs\F \\
G' & = \R{g_*}v^*\fs\F \\
H' & = \R{g_*}\gs\sh{u}\F = \R{g_*}\gs u^*\F .
\end{align*}
The inner square is \eqref{diag:sh-theta}.
The arrows pointing northeast and northwest arise from $\ga{f}(\sh{f}\F)$
while the southeast pointing arrow is from $\ga{g}(\sh{u}\F)$.
The inner square, being \eqref{diag:sh-theta}, commutes.
The outer square commutes by the definition of $\beta_\F$
(see \cite[Definition\,7.3]{dfs}). Now consider the four
trapezoids squeezed between the two squares. The top-most
trapezoid commutes by the functoriality of $\theta$. The one's
on the left and the bottom commute by the definition of
$\ga{f}$ and $\ga{g}$. This does not allow us (a priori) to conclude
that the trapezoid on the right commutes---and we are not
interested in showing this---but it does allow us to conclude
that \eqref{eq:tau-beta} holds (since the arrows labeled $\theta$
can be reversed), which is what we wished to show. For completeness
we point out (after the fact) that hence the remaining trapezoid also
commutes.
\qed

\begin{rem}{\em We have not stated \Pref{prop:basech} in its full
generality; only as much as we need for the main conclusions of this
paper. The assumption that the residue field of a point $q\in \V$
coincides with the residue field at $u(q)$ can be relaxed, but
this would involve proving a base change for the punctual trace
$\Tr{f,x}$ for $x\in \X$ closed in $f^{-1}(f(x))$ and we lead us away
from what we have set out to prove. One can also relax 
the assumption that $u$ is adic. In that case the conclusion would be that
$\iGp{\V}(v^*\ga{f}(\F))$ is equivalent to $\ga{g}(\sh{u}\F)$.
In greater detail, if the only assumption on $u$ is that
it is {\'e}tale, and $\beta_\F:\R\iGp{\V}v^*\fs\F\to\gs\R\iGp{\U}u^*\F$
is the base change map in \cite[Theorem\,7.4]{dfs},
then, making the identifications $\sh{v}\G=\R\iGp{\V}v^*\G$ ($\G\in\Coz(\X)$) 
and $\sh{u}\F=\R\iGp{\U}u^*\F$, one can prove that
$$
\beta_\F\circ\R\iGp{\V}v^*\ga{f}(\F)=\ga{g}(\sh{u}\F)\circ\beta'_\F.
$$
We leave the details (modulo the assertion on the behavior of
the punctual trace) to the reader.  We will have no occasion to
use the more general assertion made in this remark. }
\end{rem}


\section{\bf Smooth maps}\label{s:smooth}
In this section we prove that $\ga{f}:Q_\X\sh{f}\to \fs Q_\Y$ is an isomorphism
when $f:\X\to \Y$ is a smooth pseudo-finite type map. The proof we give can be
broken into two parts: we show if there is {\it any} functorial 
isomorphism $\zeta_f:Q_\X\sh{f} \iso \fs Q_\Y$ (one for
each smooth $f$) then $\ga{f}$ is also an isomorphism (this is essentially
done in subsection \ref{ss:spf-maps}); and we also show (essentially
in \Tref{thm:verdier}) that there is indeed such an isomorphism
$\zeta_f$. Indeed, in this paper, that is the only
role played by \Tref{thm:verdier}---it provides a $\zeta_f$
as the inverse of $\xi_f$ (see \eqref{eq:v}). We do not, 
in this paper, examine the relationship between the $\zeta_f$ 
we produce (through Verdier's classic trick) and $\ga{f}$ (are they equal?). 
The principal difficulty in investigating this is our lack of knowledge 
of the composite
$$
\R{f_*}\R\iGp{\X}\omega_f[d] \iso \R{f_*}\fs\F \xrightarrow{\ttr{f}} \cO_\Y,
$$
---the isomorphism coming from \Tref{thm:verdier}---for $f$ pseudo-proper
smooth of relative dimension $d$. Here and below $\omega_f$ is as in 
\cite[Definition\,2.6.4]{lns}.

\subsection{Verdier's isomorphism}\label{ss:verdier}
We will be modifying Verdier's argument
in \cite[Theorem\,3,\,p.397]{f!} to formal schemes to connect
differential forms to duality. We begin with the formal scheme
version of the fundamental local isomorphism 
\cite[p.180,\,Corollary\,III.7.3]{RD}.

\begin{lem} Let $f:\X\to \Y$ be a closed immersion of 
noetherian formal schemes, and $\J \subset \cO_{\Y}$ the coherent
ideal which is the kernel of the natural surjective map
$\cO_\X \twoheadrightarrow f_*\cO_\Y$. Suppose $\J$ is locally generated by
a regular sequence of $d$ elements so that $\J/\J^2$ can be regarded
as a locally free $\cO_\X$--module of rank $d$. Then for $\F \in \wDqcp(\Y)$
we have a functorial isomorphism
$$
f^!\F \simeq \iGp{\X}(\bL{f}^*\F)\otimes_{\X}
\wedge^d_{\cO_\X}(\J/{{\J}^2})^*[-d].
$$
\end{lem}

\proof The natural map $f^!\R\iGp{\Y}\to f^!$ is an isomorphism
by \cite[Corollary\,6.1.5\,(b)]{dfs}. Using this and \cite[Example\,6.1.3\,(4)]{dfs}
we conclude that for $\F\in \wDqcp(\Y)$, we have a functorial isomorphism
$$
f^!\F \iso {\bar{f}}^*\R\sHomb_\Y(f_*\cO_{\X},\,\R\iGp\Y\F).
$$
The rest of the proof is as as in \cite[p.180,\,Corollary\,III.7.3]{RD} 
using ``way-out". Note that since $f$ is adic we have 
$\bL f^*\R\iGp\Y=\R\iGp\X\bL f^*$ (cf. \cite[Corollary\,5.2.11(c)]{dfs}).
\qed

For the next result, recall that $\bbG$ is the category whose
objects are noetherian formal schemes and whose morphisms
are composites of pseudo-proper maps and open immersions.

\begin{thm}\label{thm:verdier} Suppose a morphism $f:\X\to \Y$ in
$\bbG$ is smooth of relative dimension $d$. For $\F \in \wDqcp(\Y)$
we have a functorial isomorphism
\begin{equation}\label{eq:verdier}
\R\iGp{\X}(f^*\F\otimes_{\X}\omega_f[d]) \iso f^!\F.
\end{equation}
\end{thm}

\proof Let $\X^{2}= \X\times_{\Y}\X$, $p_1,p_2:\X^{2}\to \X$ the
two projections and $\delta:\X \hookrightarrow \X^{2}$ the diagonal (cf.
commutative diagram below).
$$
\xymatrix{
 & \X^{2} \ar[d]^{p_1} \ar[r]^{p_2} & \X \ar[d]^f \\
\X \ar@{^{(}->}[ur]^{\delta} \ar[r]_{{\boldsymbol{1}}} & \X \ar[r]_f & \Y
}
$$
Note that $\delta$ is {\it adic} and is a closed immersion. If 
$\J_\delta \subset \cO_{\X^{2}}$ is the ideal of the closed immersion 
$\delta$ then $\J_\delta$ is locally given by a regular sequence,
and one checks, as follows, that $\wedge_{\cO_{\X}}^d(\J_\delta/\J_\delta^2) 
= \omega_f$. First, by \cite[Prop.\,2.6.8]{lns} applied to the composite
\[
\X \xrightarrow{\delta} \X^2 \xrightarrow{p_1} \X
\]
we see that $\delta^*{\widehat{\Omega}}^1_{p_1} = \J_{\delta}/\J_{\delta}^2$.
Next, by applying \cite[Prop.\,2.6.6]{lns} to the product $(\X^2,\,p_1,\,p_2)$,
we see that ${\widehat{\Omega}}^1_{\X/\Y}$ can be identified with 
$\J_\delta/\J_{\delta}^2$, whence 
$\omega_f =\wedge_{\cO_{\X}}^d(\J_\delta/\J_\delta^2)$. The above results
on differentials are probably true in greater generality (as they are
for ordinary schemes), and can presumably be proven via the techniques
in \cite[p.\,126]{ega4}.
 
Consider the commutative diagram above
(with ${\boldsymbol 1} = {\boldsymbol 1}_{\X}$, the identity map).
According to \cite[Theorem\,3]{dfs} (see also \cite[Theorem\,7.4]{dfs})
we have
$$
p_1^!f^*\F \iso \R\iGp{\X^{2}}p_2^*f^!\F.
$$
Since $f$ is not assumed to be pseudo-proper, we have implicitly used 
the localization theorems for $\shr{\boldsymbol{-}}$ proved by
S. Nayak \cite{suresh} (See also \eqref{eq:nay-bc}). Applying 
$\delta^!$ to both sides of the above
isomorphism we get, using $\delta^!p_1^! \simeq {\boldsymbol 1}^! = \R\iGp\X$,
\begin{align*}
\R\iGp{\X}f^*\F & \iso \delta^!\R\iGp{\X^{2}}p_2^*f^!\F \\
 & \iso \left(\R\iGp{\X}\bL\delta^*p_2^*f^!\F\right)\otimes_\X
\omega_f^{*}[-d])
\end{align*}
(the second isomorphism is via the Lemma). This means we have
\begin{align*}
\R\iGp{\X}(f^*\F\otimes \omega_f[d]) 
& \iso \R\iGp{\X}\bL\delta^*p_2^*f^!\F \\
& \iso \R\iGp\X f^!\F  = f^!\F \qquad{\text{(since\,$f^!\F\in \Dqct(\X)$).}}
\end{align*}
\qed

\begin{rem} {\em Let $(\X,\,\De')\xrightarrow{f} (\Y,\,\De)$ be a smooth map 
in $\cbbFc$. Assume that the relative dimension is constant, and is $d$.
By \cite[Lemma\,5.1.3]{lns} and \cite[Main Theorem\,(iii)]{lns}
for $\F\in \Cozs{\De}(\Y)$ we have an isomorphism
\begin{equation}\label{eq:lns}
\iGp{\X}(f^*\F\otimes_\X\omega_f[d)] \iso Q_{\X}\sh{f}\F
\end{equation}
and hence a functorial isomorphism
\begin{equation}\label{eq:v}
\xi_f(\F): \fs\F \iso Q_\X\sh{f}\F
\end{equation}
given by $\xi_f = \eqref{eq:lns}\circ\eqref{eq:verdier}^{-1}$.
This gives a unique map
(functorial in $\F\in\Cozs{\De}(\Y)$) in $\Cozs{\De'}(\X)$
\begin{equation}\label{eq:gammE}
\gas{f}(\F): \sh{f}\F \longrightarrow {\sh{f}\F}
\end{equation}
such that $Q_\X(\gas{f}) =\xi_f\circ\ga{f}$.
It is elementary to see that $\ga{f}(\F)$ is an isomorphism if and
only if $\gas{f}(\F)$ is an isomorphism. Note that $\gas{f}$ is defined
even when the relative dimension is not constant (by defining it
separately on each connected component), and in this case also
$\gas{f}(\F)$ is an isomorphism if and only if $\ga{f}(\F)$ is.
In this paper we do not attempt to nail down $\gas{f}$. This
question is, at the bottom, the same as the question of the
relationship between $\xi_f$ and $\ga{f}$ discussed
briefly at the beginning of this section. 
}
\end{rem}

\subsection{Smooth pseudo-finite maps}\label{ss:spf-maps}
Suppose $f:(\X,\,\De')\to (\Y,\,\De)$ 
is a map in $\cbbFc$ which is smooth and pseudo-finite. 
In this situation we note that $f_*\F$ is in $\Coz(\Y)$ for
every $\F\in \Coz(\X)$. Moreover, in addition to the trace map
$\Tr{f}:f_*\sh{f} \to {\boldsymbol 1}$ between functors on $\Coz(\Y)$
we have another trace map (arising from Grothendieck duality)
$$
\vTr{f}:f_*\sh{f} \to {\boldsymbol 1}
$$
defined as follows. We have the map
$$
Q_\Y(f_*\sh{f}\F) = \R f_*Q_\X(\sh{f}\F) \xrightarrow{\xi_f^{-1}} 
\R f_*\fs\F 
\xrightarrow{\ttr{f}(\F)} \F = Q_\Y(\F).
$$
Since $f_*\sh{f}\F$ and $\F$ lie in $\Cozs{\De}(\Y)$, we define
$\vTr{f}(\F)$ as the unique map of Cousin complexes satisfying
$\ttr{f}(\F)=Q_\Y(\vTr{f}(\F))$.
It is easy to check that $\vTr{f}$ is functorial
in $\F\in \Coz(\Y)$. The two traces are clearly related via the equation
\stepcounter{thm}
\begin{equation*}\label{eq:vTr}\tag{\thethm}
\vTr{f}\circ f_*\gas{f} = \Tr{f}.
\end{equation*}
This new trace $\vTr{f}$ has the expected universal property, viz.:

\begin{lem}\label{lem:psfin} For $f$ as above and $\F$ a Cousin
complex on $(\Y,\,\De)$  the pair $(\sh{f}\F,\,\vTr{f}(\F))$ 
represents the functor $F(\G)=\Hom_\De(f_*\G,\,\F)$ of Cousin 
complexes $\G$ on $(\X,\,\De')$.
\end{lem}
\proof Suppose we have a map of complexes $\varphi: f_*\G \to \F$. 
By the universal property of $(f^!,\,\ttr{f})$ we get a unique
map in $\varphi':\G \to \sh{f}(\F)$ in $\Dqc^+(\X)$
such that $\vTr{f}\R{f}_*\varphi' = \varphi$.
Since $\G$ and $\sh{f}\F$ are in $\Coz(\X)$ therefore
$\varphi'$ has a unique representative (which we also
denote $\varphi'$) in $\Coz(\X)$. Since $f_*\G$, $f_*\sh{f}\F$
and $\F$ are all in $\Cozs{\De}(\Y)$, by Suominen's results
in \cite{suominen}, the relationship between
$\R{f}_*\varphi'$ and $\varphi$ translates to an equality
(in $\Coz(\Y)$) $\vTr{f}f_*\varphi'=\varphi$.
This proves the Lemma.
\qed

\begin{exams}\label{ex:smooth} We give a few examples of $f$ smooth and
pseudo-finite for which $\gas{f}$ 
is an isomorphism. These are  more or less an immediate 
consequence of \Lref{lem:psfin}. The idea is to show that
$(\sh{f}\F,\,\Tr{f}(\F))$ represents the functor $F$ in the
Lemma. Note that if we show this for  $f$ and $\F$
as above we would indeed have proved (via \eqref{eq:vTr})
that $\gas{f}(\F)$ (and hence $\ga{f}(\F)$) is an isomorphism.

(1) Let $(\X,\,\De) \in \bbFc$ with $\J\in \cO_\X$ a coherent ideal
containing an ideal of definition of $\X$. Denote by ${\X^*}$ the
completion of $\X$ with respect to $\J$ and by $\kappa:\X^*\to \X$
the resulting map. The map $\kappa$ is pseudo-finite and smooth of
relative dimension zero (i.e. $\kappa$ is {\'e}tale and pseudo-proper).
For $\F\in \Coz(\X)$ and $\G\in \Coz(\X^*)$ one checks via
\cite[Lemma\,10.3.8(i)]{lns}, [{\it Ibid.}, Prop.\,10.3.9] and
{\it especially} [{\it Ibid.},\,Remark,\,10.3.10]
that
$\R\iG\J\F = \iG\J\F$, $\sh{\kappa}\F = \kappa^*\iG\J\F$,
$\kappa_*\sh{\kappa}\F=\iG\J\F$ and 
$\kappa^*\kappa_*\G= \G$ (see \cite[Proposition\,5.2.8]{dfs} for
the last relation). Needless to say, we are using equality signs for
many well-known canonical isomorphisms.

The Trace map
$\Tr{\kappa}(\F): \kappa_*\sh{\kappa}\F\to \F$ is the natural map
$\iG\J\F\hookrightarrow \F$. Let $F$ be the functor
$$
F = \Hom_\De(\kappa_*-,\,\F)
$$
on $\Cozs{\sh\kappa\De}({\X^*})$. If $\G\in \Coz({{\X^*}})$,
then $\iG\J\kappa_*\G = \kappa_*\G$ and any map 
$\kappa_*\G \xrightarrow{\psi} \F$ in $\Coz(\X)$ must factor (uniquely)
as
\begin{equation}\label{eq:kap-fac}
\kappa_*\G\xrightarrow{\psi'}\iG\J\F \xrightarrow{\Tr{\kappa}}\F.
\end{equation}
Now $\kappa^*(\psi')$ is a map from $\G$ to $\kappa^*\iG{\J}F=\sh{\kappa}\F$
and $\kappa_*\kappa^*(\psi')=\psi'$. Thus $\varphi=\kappa^*(\psi')$
solves $\Tr{\kappa}(\F)\kappa_*\varphi = \psi$, and this is the only solution. 
Indeed for any solution $\varphi$ we have $\kappa_*\varphi=\psi'$ by the 
uniqueness of the factorization \eqref{eq:kap-fac} above, and hence 
$\varphi=\kappa^*\kappa_*(\varphi) =\kappa^*(\psi')$. 
We have therefore proven that 
$(\sh{\kappa}\F,\,\Tr{\kappa}(\F))$
represents $F$ and hence that $\gas{\kappa}$ and $\ga{\kappa}$
are isomorphisms. 

If $\X=\Spf(R)$, $J\subset R$ the open ideal corresponding to $\J$,
$F=\Gamma(\X,\,\F)$, ${R^*}$ the $J$-adic completion of $R$,
then $\Gamma({\X^*},\,\sh\kappa\F)=\Gamma_{J}F$, 
$\sh\kappa\F=\Gamma_{J}F^{\,\sim{R^*}}$ and $\kappa_*\sh\kappa\F
=\Gamma_{J}F^{\,\sim {R}}$. The map $\Tr{\kappa}$ corresponds
to the inclusion $\Gamma_{J}F\subset F$. 

(2) Let $(R,\fm,k)$ be a complete noetherian local ring, $\X=\Spf{(R,\fm)}$
and $x\in \X$ the unique point in $\X$. Suppose that $(\Y,\,\De)$ is
an object in $\bbFc$ and that $f:\X\to \Y$ a smooth pseudo-proper map of 
relative dimension $d$. Note that $f$ is pseudo-finite. 
Let $y=f(x)$ and $\De(y)=p$ so that $\sh{f}\De(x) =p$. Note that $y$ is
a closed point of $\Y$. Set $A=\cO_{\Y,y}$ and 
$\omega_R=\Gamma(\X,\,\omega_f)$. 
For $\F\in\Cozs{\De}(\Y)$ we may---via \eqref{iso:sharp-h}---make the 
identification
$$
\sh{f}\F = i_x\left(\Hr^d_\fm(\F(y)\otimes\omega_R)\right)[-p].
$$
If 
$$
\res{R/A}: \Hr^d_\fm(\F(y)\otimes\omega_R) \to \F(y)
$$
is the residue map \eqref{map:res-ic} then the trace map
$\Tr{f}(\F):f_*\sh{f}\F \to \F$ is the composite
$$
i_y\left(\Hr^d_\fm(\F(y)\otimes\omega_R)\right)[-p]
\xrightarrow{{\text{via}}\,\res{R/A}} i_y\F(y)[-p]
\xrightarrow{\text{natural}} \F.
$$
We want to show that $(\sh{f}\F,\,\Tr{f}(\F))$ represents
the functor $F$ in \Lref{lem:psfin}. Without loss of generality
we assume that $p=0$. By \eqref{map:ic-phi} $(\sh{f}\F,i_y\res{R/A})$
represents the functor $\Hom_\De(f_*-,\,\F(y)^s)$ on
$\Cozs{\sh{f}\De}(\X)$. It follows easily that 
$(\sh{f}\F,\,\Tr{f}(\F))$ represents the functor
$$
F=\Hom_\De(f_*-,\,\F)
$$
on $\Cozs{\sh{f}\De}(\X)$. Thus this is another instance where
$\gas{f}$, and hence $\ga{f}$ is an isomorphism.
\end{exams}

The following lemma, based on the above examples, is very useful
in establishing that $\ga{f}$ is an isomorphism for smooth $f$.

\begin{lem}\label{lem:sp-case-smooth}
Let $(\X,\,\De)$ be an object of $\bbFc$, $x\in \X$ a closed 
point, $R$ the completion
of $\cO_{\X,x}$ at its maximal ideal, and
$$
\kappa\colon \X^*\set \Spf{(R,\,\fm_R)} \longrightarrow \X
$$ 
the resulting pseudo-finite {\'e}tale map. For a map $\alpha\colon \F \to \G$
in $\Cozs{\De}(\X)$ the following are equivalent
\begin{itemize}
\item[(a)] $\kappa^!(Q_\X(\alpha))$ is an isomorphism;
\item[(b)] $\alpha(x)\colon \F(x) \to \G(x)$ is an isomorphism.
\end{itemize}
\end{lem}

\proof Without loss of generality, we may assume that $\De(x)=0$.
Note that $\kappa$ is a special case of (1) and (2) in
\ref{ex:smooth} and hence we have isomorphisms
$$
\kappa^!\circ Q_\X \iso i_x\left[(\cdot)(x)\right] 
$$
(and $\Tr{\kappa}(\F)$ may be identified with
the inclusion $i_x\F(x)\hookrightarrow \F$
for $\F\in \Coz(\X)$). The Lemma follows.
\qed

\subsection{The isomorphism theorem for smooth maps}
Let $(A,\fm)$ be a local ring and assume $A$ is complete with respect
to an ideal $I$ (not necessarily $\fm$-primary) of $A$. Let $\Y=\Spf(A,I)$,
$y\in\Y$ the unique closed point, $\De$ a codimension function on
$\Y$, and $p=\De(y)$. Assume $\Y\in\bbF$. Let $K$ be an $A$-module
which is an injective hull of the residue field of $A$, and consider
the $\cO_\Y$-module $K^{\sim{A}}\in \Aqct(\Y)$. Set $\F=K^{\sim{A}}[-p]$.
Note that $\F\in\Coz(\Y)$.

\begin{lem}\label{lem:inj-hull} In the above situation, if 
$f:(\X,\,\De')\to (\Y,\,\De)$ is a smooth map in $\cbbFc$, then
$\gas{f}(\F)$ is an isomorphism.
\end{lem}

\proof Without loss of generality we may assume that $p=\De(y)=0$.
To reduce notational clutter we write $\cR=\sh{f}\F$ and
$\lambda=\gas{f}(\F)$. Since the
question is local on $\X$ we assume $\X=\Spf(R,J)$. Note that
if $\J'$ is the coherent $\cO_\X$ ideal corresponding to the
open ideal $\fm R+J$ of $(R,J)$ then $\iGp{\J'}\cR=\cR$ and
therefore
\begin{equation}\label{eq:J-gamma}
\iGp{\J'}(\lambda) = \lambda
\end{equation}

We first prove the Lemma under the assumption that $I$ is $\fm$-primary,
i.e. $(A,\fm)$ is a complete local ring and $\Y=\Spf(A,\fm)$. 

Now $\F$ is a residual complex on $\Y$ (cf. \cite[\S\S\,9.1]{lns}). 
By \cite[Proposition\,9.1.4]{lns} we conclude that $\cR$ is a residual
complex on $\X$. 

We have a natural isomorphism of $R$-modules
(from the definition of $t$-dualizing complexes, and from the nature 
of maps between Cousin complexes)
$$
R \iso \Hom_{\De'}(\cR,\,\cR) =: [\cR,\,\cR]
$$
given by $r\mapsto \mu_r$, $\mu_r$= multiplication by $r$. In
particular $\lambda=\gas{f}(\F)$ corresponds to a unique element
$r_f$ in $R$. We have to show that $r_f$ is a unit in $R$ (this is
equivalent to showing that $\lambda$ is an isomorphism). This is
the same as showing that the image of $r_f$ in ${\widehat{R}}_\fn$
is a unit for every maximal ideal $\fn$ of $R$ (${\widehat{R}}_\fn$
= completion of $R$ at $\fn$). Pick such an $\fn$; it corresponds to
a closed point $x$ of $\X$. With $[\cR(x),\,\cR(x)]\set 
\Hom_R(\cR(x),\,\cR(x))$, we have a commutative diagram
$$
\xymatrix{
R \ar[d] \ar[r]^-{\Iso} & [\cR,\,\cR] \ar[d]\\
{\widehat{R}}_\fn \ar[r]^-{\Iso} & [\cR(x),\,\cR(x)]
}
$$
where the left column is the usual completion map and the right
column is the result of applying the functor $(\cdot)(x)$ on
morphisms in $\Coz(\X)$. The horizontal arrow at the bottom
is $\hat{r}\mapsto \mu_{\hat{r}}$, where $\mu_{\hat{r}}$
is multiplication by ${\hat{r}}$. (Note that $\cR(x)$ is
an ${\widehat{R}}_\fn$ module being the injective hull
of the residue field of $R_{\fn}$.) This bottom arrow
is an isomorphism by Matlis duality. Thus the image
of $r_f$ in ${\widehat{R}}_\fn$ is a unit if and only if
$\lambda(x): \cR(x) \to \cR(x)$ is an isomorphism.

Let $\X^*=\X^*_{\fn}=\Spf({\widehat{R}}_\fn)$ (cf. \Ssref{ss:localrings})
and $\kappa:\X^*\to \X$
the natural map. By \Tref{thm:compare}\,(d) we have
a commutative diagram
$$
\xymatrix{
\sh{\kappa}\cR \ar[d]_{\ga{\kappa}} \ar[r]^-{\Iso} & \sh{(f\kappa)}\F
\ar[dd]^{\ga{f\kappa}} \\
\kappa^!\cR \ar[d]_{\kappa^!(\ga{f})} & \\
\kappa^!f^!\F \ar[r]^-{\Iso} & (f\kappa)^!\F
}
$$
According to \ref{ex:smooth}\,(2), $\ga{f\kappa}$ and $\ga{\kappa}$
are isomorphisms. Therefore $\kappa^!(\ga{f})$---or, equivalently
$\kappa^!(Q_\X\lambda)$---is an isomorphism.
By \Lref{lem:sp-case-smooth}, this means $\lambda(x)$ is an isomorphism.
Thus we are done if $\Y=\Spf(A,\fm)$. (See also 
\cite[p.124,\,Lemma\,3]{ast-208}). 

Suppose $I$ is not $\fm$-primary. 
Let ${\widehat{A}}$ be the $\fm$-adic completion of $A$, $S$ the
$\fm R+J$-adic completion of $R$ and set $\U\set 
\Spf({\widehat{A}},\fm{\widehat{A}})$, $\V\set \Spf(S,{\fm}S+JS)$.
Let $u:\U\to\Y$, $v:\V\to \X$, $g:\V\to \U$ be the natural maps.
Note that $\V=\U\times_\Y \X$ and $v,g$ are the two projections.
The two maps $u$ and $v$ both satisfy the hypothesis on the
map $\kappa$ in \ref{ex:smooth}\,(1). Hence $\ga{u}$ and
$\ga{v}$ are isomorphisms. Now $\sh{u}\F = K^{\sim{{\widehat{A}}}}$,
and therefore (by what we proved for the case where $I$ is $\fm$-primary)
$\ga{g}(\sh{u}\F)$ is an isomorphism. Using \Tref{thm:compare}(d)
we see that $\ga{ug}(\F)$ is an isomorphism (since $\ga{u}$ and
$\ga{g}(\sh{u}\F)$ are isomorphisms). In other words $\ga{fv}(\F)$
is an isomorphism. This coupled with the fact that $\ga{v}$ is
an isomorphism gives (via another application of {\it loc.cit.})
that $v^!(\ga{f}(\F))$ is an isomorphism. Since 
$Q_\X(\lambda)=\xi_f\circ\ga{f}$ (see \eqref{eq:gammE}), this
implies that $v^!Q_\X(\lambda)$ is an isomorphism. Now $\ga{v}$
is an isomorphism---by \ref{ex:smooth}\,(1)--- and hence
$Q_\V\sh{v}(\lambda)$ is an isomorphism. In other words 
$\sh{v}(\lambda)$ is an isomorphism in $\Coz(\V)$.
Now, on $\Coz(\X)$ we have the identity of functors
$v_*\sh{v}=\iG{\J'}$. Therefore using \eqref{eq:J-gamma} we get
$\lambda=\iG{\J'}(\lambda) =v_*\sh{v}(\lambda)$ (the last via
\cite[10.3.8(i),\,10.3.9,\,10.3.10]{lns}) proving 
that $\lambda$ is an isomorphism.
\qed

\begin{thm}\label{thm:heart} Let $f:(\X,\,\De')\to (\Y,\,\De)$
be a smooth map in $\cbbFc$. Then the functorial maps
$\ga{f}: Q_\X\sh{f} \to \fs Q_\Y$ and 
$\gas{f}: \sh{f} \to \sh{f}$
are isomorphisms.
\end{thm}

\proof It is enough to show that $\gas{f}$ is an isomorphism.
The question is clearly local on $\X$ and by \Tref{thm:compare}(c)
it is local on $\Y$ too. So without loss of generality we assume
that $\X=\Spf(R,J)$ and $\Y=\Spf(A,I)$ with $\J=J^{\sim{R}}$ and 
$\I=I^{\sim{A}}$ being the corresponding coherent ideal sheaves
on $\X$ and $\Y$ respectively.

We reduce to the case where the Krull dimension of $\Y$ is finite.
To that end, let $y\in\Y$ be a point, $A^*$ the
$\I_y$-adic completion of $\cO_{\Y,y}$, $\U=\Spf(A^*,IA^*)$ and
$u:\U\to \Y$ the resulting adic {\'e}tale map. Consider
the cartesian square \eqref{diag:basech}. To emphasize the
role of $y$, write $v_y=v$ and $g_y=g$. Now $\ga{f}$ is an
isomorphism if and only if $v^*_y(\ga{f})$ is an isomorphism
for every $y\in \Y$. By \Pref{prop:basech} this is true if
$\ga{g_y}$ is an isomorphism for every $y\in\Y$. Since $\U$
has finite Krull dimension, the theorem is true if it is
true for bases $\Y$ of finite Krull dimension and we restrict
ourselves to this case for the rest of the proof.

Since $\Y$ is finite dimensional, all Cousin complexes on
$\Y$ are bounded. For $\F\in \Coz(\Y)$ and $n\in \mathbb{Z}$
we let $\F^n$ denote the degree $n$ component of $\F$.
For $\F\ne 0$ let $a(\F)=\max\{n|\F^n\ne 0\} - \min\{n|\F^n\ne 0\}+1$.
Note that for $\F\ne0$, $a(\F)< \infty$. The theorem is clearly
true for $\F=0$. We prove the result for non-zero $\F$ by induction
on $a(\F)$. So suppose the theorem is true for $\G\in\Coz(\Y)$
such that $1\le a(\G) < l$. Suppose $a(\F)=l$. Let
$p=\min\{n|\F^n\ne 0\}$. We have a short exact sequence of
Cousin complexes on $(\Y,\,\De)$
$$
0 \longrightarrow \G \longrightarrow \F \longrightarrow \F^p[-p]
\longrightarrow 0
$$
where $\F \to \F^p[-p]$ is the natural projection. By \cite[Lemma\,10.2.4]{lns}
$$
0\longrightarrow \sh{f}\G \longrightarrow \sh{f}\F \longrightarrow 
\sh{f}(\F^p[-p]) \longrightarrow 0
$$
is exact. Now $a(\F^p[-p])=1 <l$  and $a(\G)\le l-1$
and and hence by our induction hypothesis $\gas{f}(\G)$ and 
$\gas{f}(\F^p[-p])$ are isomorphisms. Since $\gas{f}$ is functorial, 
it follows that $\gas{f}(\F)$ is an isomorphism.

It remains to prove the theorem when $a(\F)=1$. 
In this case $\F= {\mathcal H}[-p]$ where ${\mathcal H}
=\oplus_{\De(y)=p}i_yM_y$ where $M_y$ is a zero dimensional
$\cO_{\Y,y}$--module. So, without loss of generality, we
assume that $\F=i_yM_y[q]$ ($M_y$ a zero dimensional
$\cO_{\Y,y}$--module) for some point $y\in\Y$ with $\De(y)=-q$.
We again set $\U=\Spf(A^*,IA^*)$ where $A^*$ is the
completion of $\cO_{\Y,y}$ with respect to the ideal $\I_y$.
Let $\V$, $u$, $v$ be as in \eqref{diag:basech}. Since
$\F$ is concentrated at $y$, therefore, by the construction
of $\sh{f}\F$, we have $\sh{f}\F(x)=0$ for every $x$ not in 
$f^{-1}(y)$. This means that $\sh{f}\F(x)$ is an
${\widehat{\cO}}_{\Y,y}$--module for every $x\in\X$.
It follows that $\sh{f}\F=v_*v^*\sh{f}\F$. Since $v$ is adic
we have $v^*\sh{f}\F=\sh{v}\sh{f}\F$, and hence $\sh{f}\F
=v_*\sh{v}\sh{f}\F$ and $\gas{f}(\F)=v_*\sh{v}(\gas{f}(\F))$.
According to \Pref{prop:basech}, $v^*\ga{f}(\F)$ is an
isomorphism if and only if $\ga{g}(\sh{u}\F)$ is an
isomorphism, i.e. $\sh{v}(\gas{f}(\F))$ is an isomorphism
if and only if $\gas{g}(\sh{u}\F)$ is an isomorphism. Since
$\gas{f}(\F)=v_*\sh{v}(\gas{f}(\F))$, it is enough to
show that $\gas{g}(\sh{u}\F)$ is an isomorphism. 

Thus we are reduced to the case $A$ is a local ring (so that
$\Y$ has only one closed point $y$) and $\F$ is concentrated
at $y$, i.e. $\F=i_yM[q]$, $M$ a zero dimensional $A$--module.
Let $\fm$ be the maximal ideal of $A$.
We can find an exact sequence
$$
0 \longrightarrow M \longrightarrow E^0 \longrightarrow E^1
$$
where $E^0$ and $E^1$ are injective $A$--modules with
$\Gamma_{\fm}E^i=E^i$ for $i=0,1$. Let $\G^i=(E^i)^{\sim{A}}[q]$,
$i=0,1$. Note that each $\G^i$ is an object of $\Cozs{\De}(\Y)$.
Since $\sh{f}$ is an exact functor (see \cite[Lemma\,10.2.4]{lns})
and $\gas{f}$ is functorial, we have a commutative diagram
with exact rows
$$
\xymatrix{
0 \ar[r] & \sh{f}\F \ar[d]_{\gas{f}} \ar[r] & \sh{f}\G^0
\ar[d]_{\gas{f}} \ar[r] & \sh{f}\G^1 \ar[d]_{\gas{f}} \\
0 \ar[r] & \sh{f}\F \ar[r] & \sh{f}\G^0 \ar[r] & \sh{f}\G^1
}
$$
Now each $E^i$, being injective and zero-dimensional, is of
the form $\oplus_\alpha K$, where $K$ is an injective hull
of the residue field. By \Lref{lem:inj-hull}, $\gas{f}(\G^0)$
and $\gas{f}(\G^1)$ are isomorphisms. Hence $\gas{f}(\F)$
is also an isomorphism.
\qed


\section{\bf The Cousin of the comparison map}

In the last section we showed that if $f:\X\to \X$ is smooth
and is a composite of compactifiable maps then $\sh{f}$ is a
concrete model for $\fs|{\cm}(\X)$. For non-smooth
$f$ simple counter-examples exist showing that $\sh{f}$
cannot model $\fs|{\cm}(\X)$ in general. Indeed,
let $\Y=\Spec{R}$ where $R$ is a discrete valuation ring
with residue field $k$, and let $f:\X\set\Spec{k}\to \Y$
be the natural closed immersion. Then $\fs{k^{\sim_R}} \simeq
k\oplus k[-1]$. Now $k^{\sim_R}\in {\cm}(\Y,\,\De)$
where $\De$ is the codimension function on $\Y$ which is $0$ on the
closed point. But $k\oplus k[-1]$ is not Cohen-Macaulay with
respect to $\sh{f}{\De}=0$. Note however that 
$\Ed{\sh{f}(\De)}(k\oplus k[-1])
\simeq k \simeq \sh{f}(k^{\sim_R})$. It is worth asking---for a general
map $f:(\X,\,\De')\to (\Y,\De)$ in $\cbbFc$---if $\sh{f}$ is isomorphic
to $\Ed{\De'}(\fs Q_\X)$. One of the principal results of this paper
is that this is so (see \Tref{thm:heart'}).

\subsection{Definitions and notations} Let $(\X,\De') \to (\Y,\De)$
be a map in $\cbbFc$. Define a functor $\fE:\Cozs{\De}(\Y) \to 
\Cozs{\De'}(\X)$ by setting
\stepcounter{thm}
\begin{equation*}\label{eq:fE}\tag{\thethm}
\fE \set \Ed{\De'}(\fs Q_\Y)
\end{equation*}

The functors $\sh{f}$ and $f^E$ can be compared via the map
$\Ed{\De'}(\ga{f})$. More precisely we have a functorial map
\stepcounter{thm}
\begin{equation*}\label{eq:gafE}\tag{\thethm}
\gae{f}: \sh{f} \to \fE
\end{equation*}
defined by the composite
$$
\sh{f} \iso \Ed{\De'}Q_\X\sh{f} \xrightarrow{\Ed{\De'}(\ga{f})}
\Ed{\De'}(\fs Q_\Y) = \fE.
$$
We will show that $\gae{f}$ is a functorial isomorphism. We have seen
this is true when $f$ is smooth (cf. \Tref{thm:heart}). Locally
$f$ can be factored as a closed immersion followed by a smooth map.
We therefore turn our attention to closed immersions.

\subsection{Closed immersions} Suppose $f:(\X,\,\De')\to (\Y,\De)$
is a closed immersion in $\bbFc$ and suppose $\F\in\Cozs{\De}(\Y)$.
Recall from \cite[Example\,6.1.3(4)]{dfs} that we have an isomorphism
$$
{\bar{f}}^*\R\sHomb_\Y(f_*\cO_{\X},\,\F)\iso \fs{\F}
$$
induced by the universal property of $(\fs,\,\ttr{f})$ and the
map
$$
f_*{\bar{f}}^*\R\sHomb_\Y(f_*\cO_{\X},\,\F) = \R\sHomb_\Y(f_*\cO_\X,\,\F)
\xrightarrow{{\text{evaluation\,at\,1}}} \F.
$$
Recall also from \cite{lns} that we have an isomorphism of Cousin
complexes
$$
{\bar{f}}^*\sHomb_\Y(f_*\cO_\X,\,\F) \iso \sh{f}\F
$$
such that 
$$
\xymatrix{
f_*{\bar{f}}^*\sHomb_\Y(f_*\cO_{\X},\,\F) \ar@{=}[d] \ar[rr]^-{\Iso} 
 & &
f_*\sh{f}\F \ar[d]^{\Tr{f}} \\
\sHomb_\Y(f_*\cO_{\X},\,\F) \ar[rr]_{\boldsymbol{e}} &  & \F
}
$$
commutes where ${\boldsymbol{e}}$ is ``evaluation at 1". Let
\stepcounter{thm}
\begin{equation*}\label{eq:QR1}\tag{\thethm}
Q_\Y\sHomb_\Y(f_*\cO_\X,\,\F) \lra \R\sHomb_\Y(f_*\cO_\X,\,\F)
\end{equation*}
be the obvious map obtained by applying 
$\sHomb_\Y(f_*\cO_\X,\,{\boldsymbol{-}})$ to an $\Aqct(\Y)$--injective
resolution $\F\to \J$ of $\F$. Clearly the diagram
\stepcounter{thm}
\begin{equation*}\label{diag:natural}\tag{\thethm}
\xymatrix{
Q_\X{\bar{f}}^*\sHomb_\Y(f_*\cO_{\X},\,\F)
\ar[d]_{\eqref{eq:QR1}} \ar[rr]^-{\Iso} && Q_\X\sh{f}\F 
\ar[d]^{\ga{f}}\\
{\bar{f}}^*\R\sHomb_\Y(f_*\cO_{\X},\,\F) \ar[rr]^-{\Iso}
& & \fs{\F}
}
\end{equation*}
commutes, i.e. \eqref{eq:QR1} is an aspect of $\ga{f}$. We wish
to examine a similar phenomenon at a punctual level, culminating
in a tractable description of functor $\R\Gamma_x(\ga{f}) 
(\simeq \R\Gamma_{f(x)}(f_*\ga{f}))$
for a point $x\in\X$. To that end suppose $x\in \X$, $y=f(x)$, $M=\F(y)$,
$p=-\De(y)$, $R=\cO_{\Y,y}$ and $S=\cO_{\X,x}$. There is an isomorphism 
(in $\D^+(R)$)
\stepcounter{thm}
\begin{equation*}\label{iso:ext}\tag{\thethm}
\R\Homb_R(S,\,M[-p])\iso \R\Gamma_y\R\sHomb_\Y(f_*\cO_\X,\,\F)
\end{equation*}
defined as follows. Let ${\mathcal P}$ be the open prime ideal sheaf of 
$\cO_\Y$ corresponding to the point  $y\in\Y$. Then 
$\Gamma_y=(\iG{\mathcal P})_y=\Gamma_{\fm_y}(\boldsymbol{-})_y$. 
The natural map $\iG{\mathcal P}\F\to \F$
induces a map
$$
\R\sHomb_\Y(f_*\cO_\X,\,\R\iG{\mathcal P}\F) \lra \R\sHomb_\Y(f_*\cO_\X,\,\F)
$$
and since the source of this map is $\R\iG{\mathcal P}$ stable, it factors
through the map 
$$
\R\iG{\mathcal P}\R\sHomb_\Y(f_*\cO_\X,\,\F)
\to \R\sHomb_\Y(f_*\cO_\X,\,\F).
$$ 
It is not hard to see that the induced map
$$
\R\sHomb_\Y(f_*\cO_\X,\,\R\iG{\mathcal P}\F) \to 
\R\iG{\mathcal P}\R\sHomb_\Y(f_*\cO_\X,\,\F)
$$
is an isomorphism (for $f_*\cO_\X$ is coherent and $f$ is adic). 
Now $f_*\cO_\X$ is coherent on $\Y$, and so 
on taking stalks at $y$ we get an isomorphism
$$
\R\Homb_R(S,\,\R\Gamma_y\F) \iso \R\Gamma_y\R\sHomb_\Y(f_*\cO_\X,\,\F).
$$
Since $\R\Gamma_y\F =M[p]$, we obtain the map \eqref{iso:ext}. 

Next, let
\stepcounter{thm}
\begin{equation*}\label{eq:QR2}\tag{\thethm}
Q_R\Homb_R(S,\,M) \lra \R\Homb_R(S,\,M)
\end{equation*}
be the obvious map obtained by applying $\Homb_R(S,\,{\boldsymbol{-}})$
to an $R$-injective resolution $M\to I^\bullet$ of $M$. Note that
the $0$th cohomology of \eqref{eq:QR2} is an isomorphism.

\begin{prop}\label{prop:QR} With notations as above, the following diagram
commutes:
$$
\xymatrix{
Q_R\Homb_R(S,\,M[p]) 
\ar@{=}[r] \ar[dd]_{{\eqref{eq:QR2}}[p]} &
Q_R\Gamma_y\sHomb_\Y(f_*\cO_\X,\,\F) 
\ar[r]^-{\Iso} &
\R\Gamma_yQ_\Y\sHomb_\Y(f_*\cO_\X,\,\F) 
\ar[dd]^{\R\Gamma_y{\eqref{eq:QR1}}} \\
& & \\
\R\Homb_R(S,\,M[p]) 
\ar[rr]^{\Iso}_{\eqref{iso:ext}} 
&& \R\Gamma_y\R\sHomb_\Y(f_*\cO_\X,\,\F)
}
$$
\end{prop}

\proof The proof is a straightforward unraveling of definitions.
We point out that if $\mu:\F\to \J$ is a resolution of $\F$ by
injectives in $\Aqct(\Y)$ then $\Gamma_y(\mu):M[p]\to \Gamma_y(\J)$
is an injective resolution of $M[p]$. We leave the details to the
reader.
\qed

\begin{scor}\label{cor:closed-imm} Let $f:(\X,\,\De')\to (\Y,\,\De)$ 
be a closed immersion in $\bbFc$. Then $\gae{f}:\sh{f}\to f^E$ is an 
isomorphism.
\end{scor}

\proof It is enough to show that for $x\in\X$, $\Hr_x^{\De'(x)}(\ga{f})$
is an isomorphism, or equivalently $\Hr^{-p}_y(f_*\ga{f})$ is an
isomorphism where $y=f(x)$ and $p=-\De(y)$. By \eqref{diag:natural}
and the Proposition, this is so if and only if $\Hr^{-p}(\eqref{eq:QR2}[p])$
is an isomorphism, i.e. if and only if $\Hr^0\eqref{eq:QR2}$ is an
isomorphism. But $\Hr^0\eqref{eq:QR2}$ is obviously an isomorphism.
\qed

\begin{prop}\label{prop:ext} Let $f:(\X,\De')\to (\Y,\De)$ be a closed
immersion in $\bbFc$ and $\F$ an object in $\Cozs{\De}(\Y)$. Then $\fs{\F}$
is Cohen-Macaulay with respect to $\De'$ if and only if for every
$x\in\X$ and every $i\ne 0$
$$
{\mathrm {Ext}}^i_R(S,\,M)=0
$$
where $M=\F(f(x))$, $S=\cO_{\X,x}$ and $R=\cO_{\Y,f(x)}$.
\end{prop}

\proof According to \eqref{iso:ext},
$$
\Hr^{i+\De'(x)}_x(\fs{\F}) \simeq {\mathrm {Ext}}^i_R(S,\,M),
$$
giving the proposition.
\qed

\begin{prop}\label{prop:injective} Suppose $f:(\X,\,\De')\to (\Y,\,\De)$
is a map in $\cbbFc$ and $\F\in \Cozs{\De}(\Y)$ is a complex of injective
objects of $\Aqct(\Y)$. Then $\fs\F$ is Cohen-Macaulay with respect to
$\De'$.
\end{prop}

\proof Suppose first that $f$ is a closed immersion. For $y\in \Y$,
$M=\F(y)$ is an injective $R=\cO_{\Y,y}$--module. The result follows
from \Pref{prop:ext}.

If $f$ is not a closed immersion then, locally, $f=gh$ with 
$h$ a closed immersion
and $g$ a smooth map in $\cbbFc$. By \Tref{thm:heart} we see that
$\gs{\F}\simeq \sh{g}\F$. If the source of $g$ is $({\mathscr P},\,\sh{g}\De)$,
then the Cousin complex $\sh{g}\F$ is a complex of 
$\Aqct({\mathscr P})$--injectives. By what just proved in the previous
paragraph, $h^!\sh{g}\F$ is Cohen-Macaulay. But $h^!\sh{g}\F\simeq f^!\F$
in $\Dqct^+(\X)$ and we are done.
\qed

\subsection{General maps} We are now in a position to prove

\begin{thm}\label{thm:heart'} Suppose $f:(\X,\,\De')\to (\Y,\,\De)$ is
a map in $\cbbFc$. Then the functorial map
$$
\gae{f}: \sh{f} \lra f^E
$$
is an isomorphism.
\end{thm}

\proof The question is local and therefore we may assume that $f=gh$
where $h$ is a closed immersion and $g$ is a smooth map in $\cbbFc$. By
\Tref{thm:compare}\,(d) we have a commutative diagram
$$
\xymatrix{
\sh{h}\sh{g} \ar[d]_{\ga{h}} \ar[rr]^{\Iso}& & \sh{f} \ar[dd]^{\ga{f}} \\
h^!\sh{g} \ar[d]_{h^!\ga{g}} & \\
h^!g^! \ar[rr]^{\Iso} & & f^! 
}
$$
Applying $\Ed{\De'}$ to the above diagram we get a commutative diagram
$$
\xymatrix{
\sh{h}\sh{g} \ar[rr]^{\Iso} \ar[d]_{\gae{h}} & & \sh{f} \ar[dd]^{\gae{f}} \\
\Ed{\De'}(h^!\sh{g}) \ar[d]_{{\text{via}}\,\ga{g}} & & \\
\Ed{\De'}(h^!g^!) \ar[rr]^{\Iso} & & f^E 
}
$$
Since $h$ is a closed immersion, $\gae{h}$ is an isomorphism by
\Cref{cor:closed-imm}, and by \Tref{thm:heart} $\ga{g}$ is an isomorphism.
It follows that $\gae{f}$ is an isomorphism.
\qed


\begin{thm}\label{thm:CM} Let $\F\in\Cozs{\De}(\Y)$. The following are
equivalent:
\begin{enumerate}
\item[(i)] The map $\ga{f}(\F)$ is an isomorphism.
\item[(ii)] The twisted inverse image $\fs\F$ is Cohen-Macaulay with respect
to $\De'$.
\end{enumerate}
\end{thm}

\proof Clearly (i)\,$\Rightarrow$\,(ii). To go the other way, suppose
the complex $\fs\F$ is Cohen-Macaulay. Then $\ga{f}(\F)$ is a morphism in the
category $\cm(X,\,\De')$. Since the functor $E:\cm(\X,\De') \to \Cozs{\De}(\X)$
is an equivalence of categories, it is enough to prove that $E(\ga{f}(\F))$
is an isomorphism. This follows from the Theorem.
\qed

\begin{thm}\label{thm:injective} Let $\F$ be an object in
$\Cozs{\De}(\Y)$. Then the following are equivalent
\begin{itemize}
\item[(a)] $\F$ is a complex of $\Aqct(\Y)$--injectives.
\item[(b)] $\fs\F$ is Cohen-Macaulay with respect to $\sh{f}(\De)$
for every map $f$ in $\cbbF$ with target $\Y$.
\item[(c)] $\ga{f}(\F)$ is an isomorphism for
every map $f$ in $\cbbFc$ with target $(\Y,\,\De)$.
\item[(d)] $\ga{f}(\F)$ is an isomorphism for
every closed immersion $f$ in $\bbFc$ with target $(\Y,\,\De)$.
\item[(e)] The map $\ga{f}(\F)$ is an isomorphism for every
closed immersion of the form $f\colon (X,\,\De')\to (\Y,\,\De)$ with
$X$ an ordinary integral scheme.
\end{itemize} 
\end{thm}

\proof From \Pref{prop:injective} we get (a)\,$\Rightarrow$\,(b). 
By \Tref{thm:CM} we have (b)\,$\Leftrightarrow$\,(c). Clearly
we have a chain of implications (c)\,$\Rightarrow$\,(d)\,$\Rightarrow$\,(e).
It remains to show that (e)\,$\Rightarrow$\,(a). So suppose (e) is true.
Let $y$ be a point in $\Y$ and $R=\cO_{\Y,y}$. 
We have to show that $M\set \F(y)$ is an injective $R$--module. Since
$\Gamma_{\fm_R}M = M$ it enough to show that $\mu_i(\fm_R,\,M)=0$ for
$i>0$ where, for $\fp\in\Spec{R}$, $\mu_i(\fp,\,M)$ is the $i$th Bass
number of $M$ at $\fp$. Let $X$ be the closed integral subscheme of
$\Y$ defined by the closure of $y$ in $\Y$, i.e. 
$X={\mathbf{Spec}}(\cO_\Y/\I)$ where $\I$ is the open prime ideal sheaf
corresponding to the point $y\in\Y$. Let $f\colon X \hookrightarrow \Y$
be the resulting closed immersion. Let $x\in X$ be the unique point
such that $f(x)=y$. Note that $x$ is the generic point of $X$ and
the local ring of $X$ at $x$ is $k_R$.
According to our hypothesis $\fs\F$ is Cohen-Macaulay
with respect to $\sh{f}\De$. By \Pref{prop:ext} 
$$
{\mathrm{Ext}}^i_R(k_R,\,M) = 0
$$
for $i > 0$. It follows that $\mu_i(\fm_R,\,M)=0$ for $i>0$.
\qed

\begin{rem}\label{rem:gorenstein}
{\it Note that by \cite[Theorem\,4.3.1\,IV.]{lns}, if $\F$ is
a Cousin complex on $(\Y,\De)$ consisting of $\Aqct(\Y)$
injectives (so that for $y\in\Y$, $\F(y)$ is a direct sum of injective
hulls of the residue field at $y$) then $\sh{f}\F$ is complex
of $\Aqct(\X)$--injectives for every map $f:\X\to \Y$ in
$\cbbF$.}
\end{rem}

One can eliminate references to $\ga{f}$ in the above Theorem
and state it completely in terms of Grothendieck Duality and
Gorenstein complexes---a notion which we now define. Let
$(\Y,\De)$ be as \Tref{thm:injective}. A complex
$\F$ in $\Dqct^+(\Y)$ is said to be {\em Gorenstein} with
respect to $\De$ if it
is Cohen-Macaulay with respect to $\De$ and if its Cousin
complex with respect to $\De$ consists of injective objects
in $\Aqct(\Y)$. \Tref{thm:injective} and \Rref{rem:gorenstein}
give us the following (where we decided to keep matters simple
and not list all possible equivalences obtainable from \Tref{thm:injective}
and \Rref{rem:gorenstein}):

\begin{thm}\label{thm:gorenstein} Let $(\Y,\,\De)$ be as in
\Tref{thm:injective}, and let $\F$ be an object of $\Dqct(\Y)$.
Then the following are equivalent:
\begin{itemize}
\item[(a)] $\fs\F$ is Cohen-Macaulay with respect to $\sh{f}\De$ for
every map $f$ in $\cbbF$ with target $\Y$.
\item[(b)] $\fs\F$ is Gorenstein with respect to $\sh{f}\De$
for every map $f$ in $\cbbF$ with target $\Y$.
\item[(c)] $\F$ is Gorenstein with respect to $\De$.
\end{itemize}
\end{thm}


\section{\bf The Comparison map for flat morphisms}

\Tref{thm:CM} shows that $\ga{f}$ is an isomorphism of functors if
and only if $\fs$ takes Cohen-Macaulay objects to Cohen-Macaulay
objects. This characterization, unfortunately, does not give us
much information about the map $f$. In this section we prove that
$\ga{f}$ is an isomorphism of functors if and only if $f$ is {\em flat}
(cf. \Tref{thm:flat}).

\subsection{Tor and Ext} Consider a commutative diagram in $\bbF$
$$
\xymatrix{
\X \ar[dr]_f \ar[r]^h & {\mathscr P} \ar[d]^{g}\\
& \Y
}
$$
such that $g$ is smooth and $h$ is a closed immersion (so that $h$
is adic).

\begin{prop}\label{prop:tor-ext} Let $x\in \X$ and let $S$, $R$ and
$A$ be the local rings at $x$, $h(x)$ and $f(x)$ respectively.
Let $k$ be the residue field of $A$, $\overline{R}=R\otimes_Ak$,
$K$ an $R$-injective hull of the residue field of $R$, and
$\varphi:\widehat{A}\to \widehat{R}$ the map of complete local
rings induced by $g$.
\begin{enumerate}
\item[(a)] For every integer $i$ there is an isomorphism
of $R$--modules
$$
{\mathrm {Ext}}^i_R(S,\,\varphi_{\sharp}(k)) \iso 
\Hom_R({\mathrm {Tor}}_i^R(S,\,\overline{R}),\,K),
$$
i.e. ${\mathrm{Ext}}^i_R(S,\,\varphi_{\sharp}(k))$ is the Matlis dual
of the finitely generated $R$--module ${\mathrm{Tor}}_i^R(S,\,\overline{R})$.
\item[(b)] Let $i$ be an integer. Then the following are equivalent:
\begin{enumerate}
\item[(i)] ${\mathrm{Ext}}^i_R(S,\,\varphi_{\sharp}(M))
=0$ for every finitely generated $0$-dimensional $A$--module $M$;
\item[(ii)] ${\mathrm{Ext}}^i_R(S,\,\varphi_{\sharp}(M))=0$
for every $0$-dimensional $A$--module $M$;
\item[(iii)] ${\mathrm{Tor}}_i^A(S,k)=0$.
\end{enumerate}
\end{enumerate}
\end{prop}

\proof The statements are a trifle disingenuous since both statements
are trivially true if $i$ is negative. However stating matters the
way we have avoids annoying trivialities later.

Let $d=\dim{\overline{R}}$, $\fm$ the maximal ideal of $R$
and ${\overline{\fm}}$ the maximal ideal of ${\overline{R}}$. Since
$\omega_{g,x}\simeq R$ we have $R$--isomorphisms
$$
\varphi_{\sharp}(k) \simeq \Hr^d_\fm(k\otimes_A\omega_{g,x})\simeq
\Hr^d_\fm({\overline{R}})\simeq \Hr^d_{{\overline{\fm}}}({\overline{R}}).
$$
Since ${\overline{R}}$ is a regular local ring, 
$\Hr^d_{{\overline{\fm}}}({\overline{R}})$ is an ${\overline{R}}$--injective
hull of the residue field of ${\overline{R}}$. This means we have an
$R$--isomorphism
$$
\varphi_{\sharp}(k) \simeq \Hom_R({\overline{R}},\,K)
$$
since the right side is also an ${\overline{R}}$--injective hull of the
residue field of ${\overline{R}}$. Let $F^\bullet \to S$ be an
$R$--free (and hence $A$--flat) resolution of $S$. We have
\begin{align*}
{\mathrm {Ext}}^i_R(S,\,\varphi_{\sharp}(k)) & \simeq 
\Hr^i(\Hom_R(F^\bullet,\,\Hom_R({\overline{R}},\,K))) \\
& \simeq
\Hr^i(\Hom_R(F^\bullet\otimes_R{\overline{R}},\,K)) \\
& \simeq
\Hom_R(\Hr^{-i}(F^\bullet\otimes_R{\overline{R}}),\,K) \qquad 
{\text{(since\,$K$\,is\,$R$--injective)}} \\
& \simeq
\Hom_R({\mathrm{Tor}}^R_i(S,\,{\overline{R}}),\,K),
\end{align*}
thus proving part (a).

In order to prove part (b), first note that 
${\mathrm{Tor}}_i^R(S,\,{\overline{R}}) \simeq
{\mathrm{Tor}}_i^A(S,\,k)$, since $F^\bullet\otimes_R{\overline{R}}
=F^\bullet\otimes_Ak$. 

Next note that since ${\mathrm{Ext}}^i_R(S,-)$ commutes with direct
limits, and since $\varphi_{\sharp}\simeq \Hr_{\fm}^d(-\otimes_AR)$
commutes with direct limits, (i) and (ii) are equivalent (for every
module is the direct limit of finitely generated modules).
Now suppose that for a given $i$ condition (ii) holds: 
${\mathrm{Ext}}_R^i(S,\,\varphi_{\sharp}(M))=0$
for every $M\in A_{\sharp}$. 
Since $k \in A_{\sharp}$ this implies by part (a) that 
the Matlis dual of the finitely generated $R$--module 
${\mathrm{Tor}}_i^R(S,\,{\overline{R}})$ is zero.
This means ${\mathrm{Tor}}_i^R(S,\,{\overline{R}})=0$, whence
${\mathrm{Tor}}_i^A(S,\,k)=0$.

We have to show that if $i$ satisfies condition (iii), then it
satisfies condition (i). So suppose that $i$ is such that 
${\mathrm{Tor}}^A_i(S,\,k)=0$, i.e. ${\mathrm{Tor}}^R_i(S,\,{\overline{R}})=0$. Let $F$ be the functor on finitely generated $A$--modules in $A_{\sharp}$ 
given by
$$
F \set {\mathrm{Ext}}_R^i(S,\,\varphi_{\sharp}{\boldsymbol{-}}).
$$
We have to show that $F(M)=0$ for every $M\in A_{\sharp}$ which is
finitely generated.  We proceed by induction on the
length $\ell(M)$ of $M$. If $\ell(M)=1$, then $M\simeq k$, and by
part (a) $F(k)=0$ since we have ${\mathrm{Tor}}_i^R(S,\,{\overline{R}})=0$.
If $\ell(M)>1$ we have a short exact sequence of $A$--modules
$$
0 \lra M' \lra M \lra M'' \lra 0
$$
with $\ell(M'), \ell(M'') < \ell(M)$. By induction hypothesis
$F(M')=F(M'')=0$. Now $\varphi_{\sharp}$ is exact (see 
\cite[Lemma\,10.2.4]{lns}) and hence
$$
F(M') \lra F(M) \lra F(M'')
$$
is exact. It follows that $F(M)=0$.
\qed

\subsection{Local cohomology and the twisted inverse image} In this subsection
we examine the relationship between certain local cohomology modules
associated with the twisted inverse image functor and ${\mathrm{Tor}}$ modules.
%
%

\begin{prop}\label{prop:shriek-tor} Let $f: (\X,\,\Delta')\to (\Y,\,\Delta)$
be a map in $\cbbFc$. Let $x\in\X$, $y=f(x)$, $S=\cO_{\X,x}$, $A=\cO_{\Y,y}$
and $k$ the residue field of $A$. Let $p=\De(y)$ and $q=\De'(x)$. Then
for a fixed integer $i$ the following are equivalent:
\begin{enumerate}
\item[(i)] $\Hr_x^{i+q}(\fs(i_yM[-p]))=0$ for every $M\in A_{\sharp}$.
\item[(ii)] $\Hr_x^{i+q}(\fs(i_yM[-p]))=0$ for every finitely generated
$A$--module in $A_{\sharp}$.
\item[(iii)] $\Hr_x^{i+q}(\fs\F)=0$ for every Cousin complex 
$\F\in\Cozs{\De}(\Y)$.
\item[(iv)] ${\mathrm{Tor}}_i^A(S,\,k)=0$
\end{enumerate}
\end{prop}

\proof Since the statements are local in a neighborhood of $x$, we assume
that $f=gh$, where $h:\X\to {\mathscr P}$ is a closed immersion and
$g$ is smooth (and a composite of compactifiable maps). Let $R$ be the
local ring at $h(x)$ and suppose $\varphi:{\widehat A}\to {\widehat R}$
is as in the statement of \Pref{prop:tor-ext}.

(i)\,$\Leftrightarrow$\,(ii). By \eqref{iso:ext} $\Hr^{i+q}_x(\fs i_yM[-p])$
is isomorphic to ${\mathrm{Ext}}^i_R(S,\,\varphi_{\sharp}(M))$. As in the
proof of \Pref{prop:tor-ext}(b), by taking direct limits we see that
(i) and (ii) are indeed equivalent.

(i)\,$\Rightarrow$\,(iii). Suppose $\F\in\Cozs{\De}(\Y)$. Let $M=\F(y)$,
$\G_1=\sh{g}\F$, $\G_2=\sh{g}(i_yM[-p])$. Let $z=h(x)$. Then
$\G_1(z)=\G_2(z)=\varphi_{\sharp}(M)$. By \eqref{iso:ext} applied
to the map $h$ we see that
\begin{equation*}
\begin{split}
\R\Gamma_z\R\sHomb_{\mathscr{P}}(h_*\cO_\X,\,\G_1)
& \xleftarrow[\eqref{iso:ext}]{\Iso}
\R\Homb_R(S,\,\varphi_{\sharp}M[-q]) \\
& \xrightarrow[\eqref{iso:ext}]{\Iso}
\R\Gamma_z\R\sHomb_{\mathscr{P}}(h_*\cO_\X,\,\G_2).
\end{split}
\end{equation*}

Thus (since $Q_\X\sh{g}\simeq \gs{Q_\Y}$ on $\Coz(\Y)$)
$$
\R\Gamma_xh^!\gs\F \simeq \R\Gamma_xh^!\gs(i_yM[-p])
$$
i.e.
$$
\R\Gamma_x\fs\F \simeq \R\Gamma_x\fs(i_yM[-p]).
$$
By (i) it follows that $\Hr_x^{i+q}(\fs\F)=0$.

(iii)\,$\Rightarrow$\,(iv). Set $\F=(i_yk)[-p]$ and let $K$ be
an $R$--injective hull of the residue field of $R$. Then
$$
\Hom_R({\mathrm{Tor}}_i^R(S,\,R\otimes_Ak),\,K) 
\xrightarrow[\ref{prop:tor-ext}]{\Iso}
{\mathrm{Ext}}^i_R(S,\,\varphi_{\sharp}(k))
\xrightarrow[\eqref{iso:ext}]{\Iso}
\Hr_x^{i+q}(\fs\F). 
$$
But by (iii) the last $S$--module is zero. Hence 
${\mathrm{Tor}}_i^R(S,\,R\otimes_Ak)=0$ (since its
$R$--Matlis dual is zero and it is finitely generated as
an $R$--module). In other words ${\mathrm{Tor}}_i^A(S,\,k)
\simeq {\mathrm{Tor}}_i^R(S,\,R\otimes_Ak) = 0$.

(iv)\,$\Rightarrow$\,(i). By \Pref{prop:tor-ext}\,(b), if
${\mathrm{Tor}}_i^A(S,\,k)=0$ then 
${\mathrm{Ext}}_R^i(S,\,\varphi_{\sharp}(M)=0$
for every $M\in A_{\sharp}$. By \eqref{iso:ext} this means that
$\Hr^{i+q}_x(h^!\sh{g}(i_yM[-p]))=0$ for every $M\in A_{\sharp}$. Since
$\sh{g}\(i_yM[-p])\simeq \gs\(i_yM[-p])$ (cf. \Tref{thm:heart}) we are
done.
\qed

\begin{thm}\label{thm:flat} Let $f: (\X,\,\De')\to (\Y,\,\De)$
be a map in $\cbbFc$. The following are equivalent:
\begin{enumerate}
\item[(i)] The map $\ga{f}: Q_\X\sh{f} \to \fs Q_\Y$ is an isomorphism
of functors.
\item[(ii)] If $\F\in\Cozs{\De}(\Y)$ then $\fs{Q_\Y}\F\in \cm(\X,\,\De')$.
\item[(iii)] The functor $\fs|_{\cm(\Y,\,\De)}$ takes values in 
$\cm(\X,\,\De')$.
\item[(iv)] The map $f$ is flat.
\end{enumerate}
\end{thm}

\proof By \Tref{thm:CM} (i) and (ii) are equivalent. Moreover (ii) and
(iii) are clearly equivalent.
We will show that (ii)\,$\Leftrightarrow$\,(iv).
By \Pref{prop:shriek-tor}, $\fs\F$ is Cohen-Macaulay with respect to
$\De'$ for every $\F\in\Cozs{\De}(\Y)$ if and only if 
for every $x\in\X$ and every 
$i\ne 0$ ${\mathrm{Tor}}_i^{\cO_{\Y,y}}(\cO_{\X,x},\,k)=0$ where $y=f(x)$
and $k$ is the residue field at $y$. By 
\cite[p.\,174,\,Theorem 22.3\,(i) and (iii)]{matsumura} 
this is equivalent to $\cO_{\X,x}$ being flat over $\cO_{\Y,y}$. 
\qed


\section{\bf The universal property of the trace}\label{s:u-prop}

If $f:(\X,\,\De')\to (\Y,\,\De)$ is a {\em pseudo-proper} map
in $\bbFc$ and $\rho:f_*\C\to \F$ a map of complexes 
where $\C\in \Cozs{\De'}(\X)$ and $\F\in \Cozs{\De}(\Y)$, then
the resulting map $C\to \fs\F$ in $\Dqct^+(\X)$ induces---on applying
the Cousin functor $\Ed{\De}$ and the inverse of $\gae{f}(\F)$---a
map of Cousin complexes $\alpha(\rho):\C\to \sh{f}\F$. This suggests
that $(\sh{f}\F,\,\Tr{f}(\F))$ represents the functor $\Hom_\Y(f_*\C,\,\F)$
of Cousin complexes $\C\in\Cozs{\De'}(\X)$. What is required is to
show that $\rho=\Tr{f}(\F)\circ\alpha(\rho)$ and that $\delta=\alpha(\rho)$
is the only solution of the equation $\rho=\Tr{f}(\F)\circ\delta$.
This section proves these assertions and hence proves that 
$(\sh{f}\F,\,\Tr{f}(\F))$ has a universal property giving us a duality 
theory for Cousin complexes.

Throughout this section we fix a morphism
$f:(\X,\,\De')\to (\Y,\,\De)$ in $\bbFc$.

\subsection{Duality for Cousin complexes}

For the rest of this section the map $f$ is assumed to be {\em pseudo-proper}.
For $\C\in\Cozs{\De'}(\X)$, $\F\in\Cozs{\De}(\Y)$ and a map of complexes
$\rho:f_*\C\to \F$ we define a map
\stepcounter{thm}
\begin{equation*}\label{eq:alpha-shr}\tag{\thethm}
{\tilde{\alpha}}(\rho):Q_\X\C\to \fs\F
\end{equation*}
as the unique map such that $\ttr{f}(\F)\circ\R{f_*}({\tilde{\alpha}}(\rho)) 
=Q_\Y(\rho)$ (we are implicitly using $\R{f_*}\C\simeq Q_\Y f_*\C$.)
The natural isomorphism $\C\iso \Ed{\De'}Q_\X\C$ followed by
$\Ed{\De'}({\tilde{\alpha}}(\rho))$ gives us a map
\stepcounter{thm}
\begin{equation*}\label{eq:alphaE}\tag{\thethm}
\alpha'(\rho):\C\to f^E\F
\end{equation*}
in $\Coz(\X)$. Since $\gae{f}$ is an isomorphism, we can define a map
\stepcounter{thm}
\begin{equation*}\label{eq:alpha-sha}\tag{\thethm}
\alpha(\rho): \C\to \sh{f}\F
\end{equation*}
in $\Coz(\X)$ as the map satisfying 
$\gae{f}(\F)\circ\alpha(\rho)=\alpha'(\rho)$.

For the rest of this section we fix the data $(\C,\,\F,\,\rho)$ where
$\C\in\Cozs{\De'}(\X)$, $\F\in\Cozs{\De}(\Y)$ and $\rho$ is a map
of complexes $\rho:f_*\C\to \F$. Note that if $x\in\X$ and $y\in\Y$
then $\rho$ induces a map $\C(x)\to \F(y)$. By \cite[Lemma\,10.2.1]{lns} this
map is zero unless $y=f(x)$ and $x$ is closed in $f^{-1}(y)$. For
$x\in \X$ closed in $f^{-1}(f(x))$ let
\stepcounter{thm}
\begin{equation*}\label{eq:rho-x}\tag{\thethm}
\rho(x): \C(x) \to \F(f(x))
\end{equation*}
be the map induced by $\rho:f_*\C\to \F$. Then 
$$
\rho = \sum_x{i_{f(x)}}\rho(x)
$$
where the sum is taken over points $x$ which are closed in their fibers over
$\Y$.

\begin{lem}\label{lem:alpha-inj} Suppose $\delta\in \Hom_{\De'}(\C,\,\sh{f}\F)$
is such that $\Tr{f}(\F)\circ f_*\delta=\rho$. Then $\delta=\alpha(\rho)$.
\end{lem}

\proof By the universal property of $(\fs,\,\ttr{f})$ we see that
$\ga{f}(\F)\circ Q_\X(\delta) = \tilde{\alpha}(\rho)$. This implies,
by the definition of $\alpha'(\rho)$ (cf. \eqref{eq:alphaE}) and
of $\gae{f}$ that $\gae{f}(\F)\circ\delta=\alpha'(\rho)$. It follows
that $\delta=\alpha(\rho)$.
\qed

Let $F^p(\C)$, $G_p(\C)$ be the complexes in \cite[\S\S\,10.2]{lns}
(cf. especially the discussion immediately following the proof of
Lemma\,10.2.4 in {\em Ibid.}.

\begin{lem}\label{lem:alpha-phi} Let $\F'$ and $\C'$ be Cousin complexes
on $(\Y,\,\De)$ and $(\X,\,\De')$ respectively and suppose
$\varphi\in\Hom_{\De'}(\C',\,\C)$, $\psi\in\Hom_\De(\F',\,\F)$ and
$\rho'\in \Hom_\Y(f_*\C',\,\F')$ are such that the diagram of
complexes
$$
\xymatrix{
f_*\C' \ar[d]_{\rho'} \ar[r]^{f_*\varphi} & f_*\C \ar[d]^{\rho} \\
\F' \ar[r]_{\psi} & \F
}
$$
commutes. Then
$$
\alpha(\rho)\circ\varphi = \sh{f}(\psi)\circ \alpha(\rho')
$$
in $\Cozs{\De'}(\X)$.
\end{lem}

\proof By the universal property of $(\fs,\,\ttr{f})$ we have
${\tilde{\alpha}}(\rho)\circ Q_\X\varphi = Q_\X\sh{f}(\psi)\circ
\tilde{\alpha}(\rho')$. The Lemma follows.
\qed

\medskip
For a Cousin complex $\C$ on $(\X,\De')$, let $\{F^p\C\}$ and $\{G_p\C\}$
be the filtrations in \cite[\S\S\,10.2]{lns} (see material immediately
following the proof of [{\it Ibid.},\,10.2.4]). Recall from 
\cite[\S\S\,10.2,(89)]{lns} that $F^p\F=\sh{f}(\sigma_{\ge p}\F)$
and $G_p(\F)=\sh{f}(\sigma_{\le p}\F)$.
Note that by \cite[Lemma\,10.2.1]{lns} the maps $\rho(x)$ of \eqref{eq:rho-x}
induce  maps
\begin{equation}
\begin{split}
\rho^+=\rho^+_p: f_*F^p(\C) & \to \sigma_{\ge p}\F \\
\rho^-=\rho^-_p: f_*G_p(\C) & \to \sigma_{\le p}\F
\end{split}
\end{equation}
where, for example, $\rho^+_p=\sum_xi_{f(x)}\rho(x)$---the sum being taken
over points $x$ such that $\De'(x)=\De(f(x))\ge p$.

\begin{lem}\label{lem:F-G-rho} The diagram
\begin{align*}
{\xymatrix{
F^{p}(\C) \ar[d]_{\alpha(\rho^+)} \ar[r] & \C \ar[d]^{\alpha(\rho)} 
\ar[r] & G_p(\C) \ar[d]^{\alpha(\rho^-)}\\
\sh{f}(\sigma_{\ge p}\F)\ar[r] & \sh{f}\F \ar[r] & \sh{f}(\sigma_{\le p}\F) }}
\end{align*}
commutes.
\end{lem}

\proof One checks, by the definition of $\rho^+$, that
$$
\xymatrix{
f_*F^{p}(\C) \ar[r] \ar[d]_{\rho^+} & f_*\C \ar[d]^{\rho} \\
\sigma_{\ge p}\F \ar[r] & \F
}
$$
commutes, where the horizontal arrows are the obvious inclusions.
The rectangle on the left in our assertion therefore commutes by
\Lref{lem:alpha-phi}. A similar argument gives the commutativity of
the rectangle on the right.
\qed

\begin{rem}{\em The above lemma asserts that 
$\alpha(\rho^+_p)=F^p(\alpha(\rho))$ and $\alpha(\rho^-_p)=G_p(\alpha(\rho))$.}
\end{rem}

\begin{prop}\label{prop:rho} $\Tr{f}(\F)\circ f_*(\alpha(\rho))=\rho$.
\end{prop}

\proof Let $x\in\X$ be closed in its fiber. By \cite[Lemma\,10.2.1]{lns} it is
enough to show that
\begin{equation}\label{eq:tr-alpha-rho}
\Tr{f,x}(\F)\circ\alpha(\rho)(x) = \rho(x).
\end{equation}
Let $p=\De'(x)$. We have the identities $\rho(x)=\rho^-_p(x)$, 
$\Tr{f,x}(\F)=\Tr{f,x}(\sigma_{\le p}\F)$ and (by \Lref{lem:F-G-rho})
$\alpha(\rho)(x)=\alpha(\rho^-_p)(x)$. Thus, without loss of generality,
we assume that $\F=\sigma_{\le p}\F$ and $\C=G_p(\C)$, by replacing
$\C$ by $G_p(\C)$, $\F$ by $\sigma_{\le p}\F$ and $\rho$ by $\rho^-_p$. 
Now let 
$$
\C'\set i_x\C(x)[-p].
$$
Since $\C=G_p(\C)$, the natural map $\varphi:\C'\to \C$--induced by the
identity map $\C'(x)\to \C(x)$--is a map of complexes. Moreover,
since $\F=\sigma_{\le p}\F$, the map $\rho(x)$ induces a map of complexes
$\rho':\C'\to \F$ such that $\rho'(x)=\rho(x)$. \Lref{lem:alpha-phi}
applied to $\varphi$ above and $\psi=1_\F$ gives us $\alpha(\rho)(x)
=\alpha(\rho')(x)$. Thus in order to establish \eqref{eq:tr-alpha-rho},
we may (and will) assume that $\C=\C'$ and $\rho=\rho'$, i.e. $\C$ is
concentrated at $x$. Let $A$ and $S$ be the completions of the
local rings at $y=f(x)$ and $x$ respectively, and $h:A\to S$ the map
induced by $f$. The pair $((\sh{f}\F)(x),\,\Tr{f,x}(\F))$ represents
the functor $\Hom_A(M,\,\F(y))$ of $0$--dimensional $S$--modules $M$.
Therefore we have a map $d:\C(x)\to (\sh{f}\F)(x)$ such that
$\Tr{f,x}(\F)\circ d=\rho(x)$. Since $\C=i_x\C(x)[-p]$, the map $d$
gives rise to a map of complexes $\delta:\C\to \sh{f}\F$ given by
the composite
$$
\C \xrightarrow{i_x(d)[-p]} i_x((\sh{f}\F)(x))[-p] \xrightarrow{\text{natural}}
\sh{f}\F.
$$
The second arrow is a map of complexes since the $n$th graded piece of
$\F$ for $n\ge p$ is zero and $x$ is a point closed in its fiber with
$\De'(x)=p$. Clearly $\Tr{f}(\F)\circ\delta=\rho$. By \Lref{lem:alpha-inj}
we are done.
\qed

For future reference we gather the results in \Lref{lem:alpha-inj}
and \Pref{prop:rho} into the following theorem.

\begin{thm}\label{thm:u-prop} Let $f:(\X,\,\De')\to (\Y,\,\De)$ be
a pseudo-proper map in $\bbFc$. For $\C\in\Cozs{\De'}(\X)$ and
$\F\in\Cozs{\De}(\Y)$ the natural map
\begin{equation*}
\begin{split}
\Hom_{\De'}(\C,\,\sh{f}\F) & \to \Hom_\Y(f_*\C,\,\F) \\
\delta & \mapsto \Tr{f}(\F)\circ f_*\delta
\end{split}
\end{equation*}
is a bifunctorial isomorphism. In particular, the pair $(\sh{f}\F,\,\Tr{f}(\F))$
represents the functor $\Hom_\Y(f_*\C,\,\F)$ of Cousin complexes $\C$
on $(\X,\,\De')$.
\end{thm}

\proof \Lref{lem:alpha-inj} proves that the map $\delta\mapsto \Tr{f}(\F)\circ
f_*\delta$ is injective. \Pref{prop:rho} shows that it is surjective.
\qed

\medskip
\Tref{thm:u-prop} has an interesting corollary when $f$ is a {\it pseudo-finite}
map, i.e. $f$ is pseudo-proper and its fibers are finite, or equivalently,
$f$ is pseudo-proper and affine. In this case, $f$ corresponds locally
(on $\Y$) to a homomorphism $\varphi:(R,I) \to (S,J)$ of adic rings
(i.e. $\varphi(I) \subset J$) and $S/J$ is a finite $R$-module. Let
$\U=\Spf{(R,I)}$ and $\V=f^{-1}(\U)=\Spf{(S,J)}$. As in \cite[\S\S\,2.1]{dfs}
set
\[
\Hom_{R,J}(S,\,F) \set \Gamma_J\Hom_R(S,\,F)
\]
where $F=\Gamma(\U,\F)$ ($\F\in \Cozs{\De}(\Y)$ as in \Tref{thm:u-prop}).
The complex $\Hom_{R,J}(S,\,F)$ can also be interpreted as the complex
of $S$--modules of {\it continuous} $R$--maps from $S$ to $F$ when $S$
is $J$--adically topologized and $F$ is discrete. Moreover
\begin{equation}\label{eq:contSF}
\Hom_{R,J}(S,F) = \dirlm{n}\Hom_R(S/J^n,\,F)
\end{equation}
by standard 3-lemma arguments. Since $\cO_\X$ is coherent and $\F$ is
a Cousin complex the complexes $\Hom_{R,J}(S,\,F)^{\widetilde{}S}$ can be 
patched as $\U=\Spf{R}$ varies over an affine open cover of $\Y$ to give a
complex $f^\flat\F$ (cf. \cite[Lemma\,2.3.5(iii)]{lns} applied to
the torsion modules $\Hom_{R,J}(S,\F(y))$ for $y\in\Y$).

If $\J$ is a defining ideal for $\X$ and $X_n=(\X,\,\cO_\X/\J^n)$ then
by \eqref{eq:contSF}
\begin{equation}\label{eq:contfF}
f^\flat\F = \dirlm{n} {i_n}_*f_n^\flat\F
\end{equation}
where $i_n\colon X_n\hookrightarrow \X$ is the natural closed immersion
and $f_n\colon X_n \to \Y$ the resulting finite (= adic and pseudo-finite)
map. It is not hard to see that $f_n^\flat\F$ is Cousin on 
$(X_n,\,\sh{f_n}\De)$, whence $f^\flat\F$ is Cousin on $(\X,\,\De')$.
Moreover, for $x\in\X$ (with $y=f(x)$) one has
\begin{equation}\label{eq:bSF}
(f^\flat\F)(x) = \Gamma_{\fm_{\widehat{S}}}\Hom_{\widehat{R}}(\widehat{S},\F(y))
=\Hom_{\widehat{R}}^c({\widehat{S}},\,\F(y))
\end{equation}
where ${\widehat{R}}$ (resp.~${\widehat{S}}$) is the
completion of the local ring at $y$ (resp.~$x$) and the superscript ``c"
over the $\Hom$ on the right refers to continuous ${\widehat{R}}$--maps with
${\widehat{S}}$ having the $\fm_{\widehat{S}}$--adic topology and $\F(y)$
having the discrete topology. The relation \eqref{eq:bSF} is obtained
by using \eqref{eq:contfF} to reduce to the case where
$f$ is finite where the relation is not hard to establish.
Using \cite[\S\,7]{I-C} (see also \eqref{map:ic-phi}) we have an isomorphism
of ${\widehat{S}}$--modules
$$
\Theta(x)=\Theta_f(x)\colon (\sh{f}\F)(x) \iso 
\Hom_{\widehat{R}}^c({\widehat{S}},\,\F(y)) = (f^{\flat}\F)(x) 
$$
characterized by the relation
$$
{\boldsymbol{e}}(x)\circ\Theta(x) = \Tr{f}(\F)(x)
$$
where ${\boldsymbol{e}}(x)\colon \Hom_{\widehat{R}}^c({\widehat{S}},\,\F(y))
\to \F(y)$ is ``evaluation at $1$". That this is a characterization is
readily seen by observing that $(\sh{f}\F(x),\,\Tr{f}(\F)(x))$ and
$(f^\flat\F(x),\,{\boldsymbol{e}}(x))$ represent the same functor,
viz.~the functor $\Hom_{\widehat{R}}(N,\,\F(y))$
of zero-dimensional ${\widehat{S}}$--modules $N$.

From the $\Theta(x)$ we get in an obvious way a map of graded $\cO_\X$--modules
$$
\Theta_f(\F)\colon \sh{f}\F \iso f^{\flat}\F
$$
which is functorial in $\F\in\Cozs{\De}(\Y)$.
A consequence of \Tref{thm:u-prop} is the following:

\begin{cor}\label{cor:finite} Let $f\colon (\X,\,\De')\to (\Y,\,\De)$
be a pseudo-finite map in $\bbFc$ and $\F$ a Cousin complex on
$(\Y,\,\De)$. Then the graded map $\Theta_f(\F)$ above is a map of
complexes.
\end{cor}

\proof By \eqref{eq:contfF} it is enough to assume $f$ is finite. 
We have a map ${\boldsymbol{e}}\colon f_*f^\flat\F \to \F$ given
by ``evaluation at $1$" (note that $f_*f^\flat\F = \sHomb_\Y(f_*\cO_\X,\,\F)$).
Clearly $(f^\flat\F,\,{\boldsymbol{e}})$ represents the same functor
that $(\sh{f}\F,\,\Tr{f}(\F))$ does. This results in an isomorphism
$$
\Theta' \colon \sh{f}\F \iso f^\flat\F
$$
such that ${\boldsymbol{e}}\circ f_*(\Theta') = \Tr{f}(\F)$. It follows
that for a point $x\in \X$, $\Theta'(x) = \Theta_f(\F)(x)$, whence
$\Theta'=\Theta_f(\F)$.
\qed


\section{\bf Variants}\label{s:(!)}

In this section we construct a variant of $\shr{\boldsymbol{-}}$ on
the full subcategory of $\bbF$ consisting of schemes which admit
a bounded residual complex (cf. \cite[\S\S\,9.1]{lns}, 
\cite[5.9]{Ye} and, for related matters, \cite[\S\S\,2.5]{dfs}).
For every such scheme, the associated category $\shrbr{\X}$ is a
full subcategory $\wDcp(\X)$ of $\wDqcp(\X)$. The method we use
is the Grothendieck's original method via residual complexes developed in
\cite{RD}, but with our canonical Cousin complex valued
pseudofunctor $\sha{\boldsymbol{-}}$ Cousin complexes (more
precisely, its ``restriction" to residual complexes) in place of
${\boldsymbol{-}}^{\boldsymbol{\De}}$ of \cite{RD}. Moreover, we
are dealing with formal schemes rather than ordinary schemes and
hence there is a need to retell the story, albeit in an
abbreviated form. The reader is advised to look at the very careful
account given by Conrad in \cite[Chapter\,3]{conrad} (especially
\S\,3.3 and \S\,3.4) to flesh out missing details in what follows.

\subsection{Preliminaries} Let $\Y$ be a noetherian formal scheme. 
Suppose $\Y$ admits a {\em bounded} residual complex 
$\cR$ \cite[\S\S\,9.1]{lns}\footnote{This forces $\Y$ to have
finite Krull dimension.}. By \cite[9.2.2\,(ii) and (iii)]{lns}
$\cR$ is a $t$--dualizing complex in the sense of \cite[Definition\,2.5.1]{dfs}.
By \cite[Theorem\,5.6]{Ye}, if $\cR'$ is another residual complex
on $\Y$ then $\cR'\simeq\cR\otimes {\mathcal L}[n]$ for an invertible
$\cO_\Y$--module ${\mathcal L}$ and an integer valued locally constant
function $n$. The residual complex
$\cR$ induces a codimension function $\De_\cR$ on $\Y$; for a point
$y$ in $\Y$, $\De_\cR(y)$ is the unique integer $p$ such that
$\Hp{y}{p}(\cR)$ is non-zero.

For any $\E\in \D(\Y)$ set
$$
\cD_\cR(\E)\set \R\sHomb(\E,\,\cR).
$$
Now, $\cR$ is a complex of $\Aqct(\Y)$--injectives (cf. \cite[9.1.3]{lns}
and \cite[2.3.6\,(ii)]{lns}), and therefore, if $\E \in \Dqct^+(\Y)$
we may, and will, make the identification
$$
\cD_\cR(\E) = \sHomb(\E,\,\cR)
$$
(cf. \cite[Proposition\,5.3.1]{dfs}). We refer the reader to
\cite[\S\S\,2.5]{dfs}---especially (a) and (c) of Proposition\,2.5.8---for
further details on $\cD_\cR$.

\subsection{Twisted inverse image via residual complexes}\label{ss:shr-residual}
Let $\rbbF$ be the full subcategory of $\bbF$ consisting of schemes
which admit {\it bounded} residual complexes. For $\Y\in \rbbF$
let $\Dc^*(\Y)$ and $\wDcp(\Y)$ be as in \Ssref{ss:conventions} and
set
$$
\shrbr{\Y}\set \wDcp(\Y).
$$

For the rest of this subsection we fix maps
$$
\V \xrightarrow{h} \W \xrightarrow{g} \X \xrightarrow{f} \Y
$$
in $\rbbF$.

Suppose $\cR$ is a residual complex on $\Y$. Define
$$
\shrbr{f}_\cR\colon \shrbr{\Y} \longrightarrow \shrbr{\X}
$$
by setting
\stepcounter{thm}
\begin{equation*}\label{eq:f(!)-R}\tag{\thethm}
\shrbr{f}_\cR \set \cD_{\sh{f}\cR}\circ\bL f^*\circ\cD_\cR\circ\R\iGp{\Y}.
\end{equation*}
This functor takes values in $\shrbr{\X}$ for the following reasons;
(a) $\E\in \shrbr{\Y} \Rightarrow \R\iGp{\Y}\E \in \Dc^*(\Y)\cap \D^+(\Y)$ by 
{\em definition} of $\shrbr{\Y}$; (b) $\E\in \Dc^*(\Y)\cap \D^+(\Y)
\Rightarrow \cD_\cR(\E)\in \Dc^-(\Y)$ by \cite[Proposition\,2.5.8\,(a)]{dfs},
(c) $\G\in \Dc^-(\Y) \Rightarrow \bL f^*\G\in\Dc^-(\X)$
and (d) $\F\in \Dc^-(\X) \Rightarrow \cD_{\sh{f}\cR}(\F)\in \Dc^*(\X)\cap
\D^+(\X)$ by \cite[Proposition\,2.5.8\,(b)]{dfs}.

Defining a pseudofunctor $\shrbrb{\boldsymbol{-}}$ on $\rbbF$ is now
a formal process given in detail in \cite[\S\,3.3]{conrad}. We
point out the main steps.

1) If $\cR'$ is another residual complex on $\Y$, then as in
\cite[p.\,135,\,(3.3.12)]{conrad} we have a comparison
map (an isomorphism)
$$
\phi_{\cR,\cR'}=\phi_{f,\cR,\cR'}\colon \shrbr{f}_{\cR'} \iso
\shrbr{f}_{\cR}
$$
stemming from $\cR'\simeq \cR\otimes {\mathcal{L}}[n]$ for some
invertible sheaf ${\mathcal{L}}$ which is compatible with a similar
relation between $\sh{f}\cR'$ and $\sh{f}\cR$ via $f^*{\mathcal{L}}[n]$.
This map does not depend on the isomorphism 
$\cR'\simeq \cR\otimes{\mathcal{L}}[n]$. We sketch the idea behind
the independence. The only non-trivial case is when ${\mathcal{L}}=\cO_\Y$
and $n=0$. The issue comes down to this: suppose $\varphi\colon \cR\iso \cR$
is an automorphism, then we have to show that the induced automorphism
$\tilde\varphi\colon f^{(!)}_\cR \iso f^{(!)}_\cR$ is the identity map.
For this assume without loss of generality that $\Y=\Spf(A,I)$. Then
$\varphi$ is given by multiplication by a unit $a\in A$. The map
$\tilde\varphi$ can be obtained by taking any path from the top left
vertex to the bottom right vertex in the commutative diagram (note that
all the arrows are invertible):
\[
\xymatrix{
\cD_{\sh{f}\cR}(\bL f^*\cD_\cR(\F))
\ar[rrr]^{\sHomb(\bL f^*\cD_\cR(\F),\,\sh{f}\varphi)}
& & &
\cD_{\sh{f}\cR}(\bL f^*\cD_\cR(\F)) \\
\cD_{\sh{f}\cR}(\bL f^*\cD_\cR(\F))
\ar[u]_{\wr}^{\cD_{\sh{f}\cR}(\bL f^*\sHomb(\F,\,\varphi))}
\ar[rrr]_{\sHomb(\bL f^*\cD_\cR(\F),\,\sh{f}\varphi)}
& & & \cD_{\sh{f}\cR}(\bL f^*\cD_\cR(\F))
\ar[u]^{\wr}_{\cD_{\sh{f}\cR}(\bL f^*\sHomb(\F,\,\varphi))}
}
\]
But all arrows are multiplication by the unit $a$ of $A$ giving that
$\tilde\varphi$ is the identity arrow.

The map $\phi_{\cR,\cR'}$ satisfies the cocycle condition for 
three residual complexes
$\cR, \cR', \cR''$ given in \cite[p.\,135,\,(3.3.13)]{conrad}. Standard
techniques give a well defined functor
\stepcounter{thm}
\begin{equation*}\label{def:f(!)}\tag{\thethm}
\shrbr{f}\colon \shrbr{\Y} \longrightarrow \shrbr{\X}
\end{equation*}
together with isomorphisms
$$
\theta_\cR = \theta_{f,\cR} \colon \shrbr{f}_\cR \iso \shrbr{f}
$$
satisfying
$$
\phi_{\cR,\cR'}=\theta_\cR^{-1} \circ \theta_{\cR'}.
$$

2) As in \cite[p.\,136,\,(3.3.15)]{conrad} the isomorphism 
$\sh{C}_{g,f}\colon \sh{g}\sh{f}\cR \iso \sh{(fg)}\cR$
gives an isomorphism 
$$
\shrbr{C}_{g,f,\cR} \colon \shrbr{g}_{\sh{f}\cR}\shrbr{f}_\cR
\iso \shrbr{(fg)}_\cR
$$
in such a way that ``associativity" holds---the last because
$\sha{\boldsymbol{-}}$ is a pseudofunctor. One checks that
the isomorphism
\stepcounter{thm}
\begin{equation*}\label{def:C(!)}\tag{\thethm}
\shrbr{C}_{g,f} \colon \shrbr{g}\shrbr{f} \iso \shrbr{(fg)}
\end{equation*}
defined by $\shrbr{C}_{g,f,\cR}$, $\theta_{f,\cR}$, $\theta_{g,\sh{f}\cR}$
and $\theta_{fg,\cR}$ is independent of $\cR$. 

This gives the required (pre)--pseudofunctor on $\rbbF$.

\subsection{Comparison of the two twisted inverse images} We will use
the discussion in \Ssref{ss:gamma-ext} to show that $\shr{\boldsymbol{-}}$
and $\shrbrb{\boldsymbol{-}}$ agree whenever both are defined (cf.
\Tref{thm:(!)!}). With that in mind we examine the two theories for
open immersions and for pseudo-proper maps.

1) Suppose $f\colon \X\to \Y$ is an open immersion in $\rbbF$. Let $\cR$
be a residual complex on $\Y$. Then $\sh{f}\cR=f^*\cR$. Moreover
$f^!=f^*\iGp{\Y}= \iGp{\X}f^*$. If $\F\in \shrbr{\Y}$ we have the
following sequence of isomorphisms
\stepcounter{thm}
\begin{equation*}\label{comp:open(!)!}\tag{\thethm}
\begin{split}
\shrbr{f}_\cR\F  &\, \,=\, \quad 
\sHomb_\X(f^*\sHomb_\Y(\iGp{\Y}\F,\,\cR),\,f^*\cR) \\
 &\, \,=\, \quad \sHomb_\X(\sHomb_\X(f^*\iGp{\Y}\F,\,f^*\cR),\,f^*\cR) \\
 & \,\,=\, \quad \cD_{f^*\cR}\cD_{f^*\cR}f^*\iGp{\Y}\F\\
&  \iso \, f^*\iGp{\Y}\F \\
& \,\,=\,\quad \fs\F
\end{split}
\end{equation*}
Let $\Phi_{f,\cR}\colon \shrbr{f}_\cR\F \iso \fs\F$ be the above composite.
It is easy to verify that 
\stepcounter{thm}
\begin{equation*}\label{map:open(!)!}\tag{\thethm}
\Phi_f \set \Phi_{f,\cR}\circ \theta_\cR^{-1} \colon \shrbr{f} \iso 
\fs|_{\shrbr{\Y}}
\end{equation*}
is independent of $\cR$.

2) Suppose $f\colon \X\to \Y$ is a map in $\rbbF$ which is pseudo-proper.
Let $\cR$ be a residual complex on $\Y$ and let $\F$ be an object in
$\shrbr{\Y}$. One then has the following sequence of maps---the first
arrow arising from the adjoint pair $(\bL f^*,\,\R f_*)$:
\stepcounter{thm}
\begin{equation*}\label{comp:tau-R(!)}\tag{\thethm}
\begin{split}
\R{f_*}\shrbr{f}\F
& \xrightarrow{\,\,\,\,\Iso\,\,\,\,} 
\R{f_*}\R\sHomb_\X(\bL{f^*}\cD_\cR\iGp{\Y}\F,\,\sh{f}\cR) \\
& \xrightarrow{\,\,\,\,\Iso\,\,\,\,} 
\R\sHomb_\Y(\cD_\cR\iGp{\Y}\F,\,\R{f_*}\sh{f}\cR) \\
& \xrightarrow{{\text{via}}\,\Tr{f}} \R\sHomb_\Y(\cD_\cR\iGp{\Y},\,\cR) = \cD_\cR\cD_\cR\iGp{\Y}\F \\
& \xrightarrow{\,\,\,\Iso\,\,\,} \iGp{\Y}\F \\
& \xrightarrow{\,\,\,\text{nat'l}\,\,} \F
\end{split}
\end{equation*}
The above composite gives a trace map
$$
\ttr{\cR}(\F) = \ttr{f,\cR}(\F)\colon \R{f_*}\shrbr{f}\F \longrightarrow \F.
$$
If $g\colon \W \to \X$ is a second pseudo-proper map, one checks using
\eqref{diag:tr-tr-f} (i.e. the transitivity property of traces on Cousin
complexes) that the following relation holds
\stepcounter{thm}
\begin{equation*}\label{eq:tr-tr-R}\tag{\thethm}
\ttr{\cR,fg} = 
\ttr{f,\cR}\circ\ttr{\sh{f}\cR}\circ\R{f_*}\R{g_*}(\shrbr{C}_{g,f,\cR})^{-1}.
\end{equation*}
Define
\stepcounter{thm}
\begin{equation*}\label{def:tau(!)}\tag{\thethm}
\ttr{f}^r \set \ttr{f,\cR}\circ\R{f_*}(\theta_\cR^{-1}).
\end{equation*}
Since $\Tr{f}(\cR)$ is compatible with Zariski localizations of $\Y$, 
functorial with respect to maps of residual complexes {\em with the
same codimension function} (whence compatible with tensoring $\cR$
by an invertible sheaf) and compatible with translations of residual
complexes, therefore one checks that 
$$
\ttr{f}^r:\R{f_*}\shrbr{f} \to {\boldsymbol{1}}_{\shrbr{\Y}}
$$
is independent of the residual complex $\cR$. Since a residual complex
on $\Y$ is a complex of $\Aqct(\Y)$--injectives therefore \Tref{thm:injective}
applies and we make the identifications
\stepcounter{thm}
\begin{equation*}\label{eq:set-up-dfs}\tag{\thethm}
\fs\cR = \sh{f}\cR \qquad \qquad  \ttr{f}(\cR) = \Tr{f}(\cR)
\end{equation*}
By \cite[2.5.12 and 6.1.5(b)]{dfs} and \eqref{eq:set-up-dfs} 
we get an isomorphism
\stepcounter{thm}
\begin{equation*}\label{iso:(!)!}\tag{\thethm}
\Phi_f\colon \shrbr{f} \iso \fs|_{\shrbr{\Y}}
\end{equation*}
---the map $\Phi_f$ being the unique map arising from the universal
property of the pair $(\fs,\,\ttr{f})$ for which the relation
\stepcounter{thm}
\begin{equation*}\label{eq:(!)!-tau}\tag{\thethm}
\ttr{f}^r = \ttr{f}\circ\R{f_*}\Phi_f
\end{equation*}
holds. 

If $g\colon \W\to \X$ is a second pseudo-proper map then the transitivity
relation \eqref{eq:tr-tr-R} above gives us a commutative diagram
\stepcounter{thm}
\begin{equation*}\label{diag:tr-tr-(!)}\tag{\thethm}
\xymatrix{
\R{f_*}\R{g_*}\shrbr{g}\shrbr{f} \ar[d]_{\ttr{g}^r} \ar@{=}[r] &
\R{(fg)_*}\shrbr{g}\shrbr{f} \ar[r]^{\Iso} &
\R(fg)_*\shrbr{(fg)} \ar[d]^{\ttr{fg}^r} \\
\R{f_*}\shrbr{f} \ar[rr]_{\ttr{f}^r} & & {\boldsymbol{1}}_{\shrbr{\Y}}
}
\end{equation*}
We are now in a position to state the following Theorem, which can be
reformulated as stating that appropriate restrictions
of the  pseudofunctors ${\boldsymbol{\shr{-}}}$ and 
${\boldsymbol{\shrbr{-}}}$  are isomorphic. 

\begin{thm}\label{thm:(!)!} There is a unique family of isomorphisms
$$
\Phi_f\colon \shrbr{f} \iso \fs|_{\shrbr{\Y}},
$$
one for each map $f\colon\X\to\Y$ in $\rbbF\cap \cbbF$ such that
\begin{itemize}
\item[(a)] If $f$ is pseudo-proper and $\F\in\shrbr{\Y}$ then $\Phi_f$
is the map \eqref{iso:(!)!}, i.e., it is the unique map satisfying
$$
\ttr{f}^r(\F) = \ttr{f}(\F)\circ\Phi_f(\F).
$$
\item[(b)] If $f$ is an open immersion then $\Phi_f$ is the isomorphism
\eqref{map:open(!)!}.
\item[(c)] If $g\colon \W \to \X$ is a second map in $\rbbF\cap\cbbF$ then
the diagram
\begin{equation}\label{diag:phi-ext}
\xymatrix{
\shrbr{g}\shrbr{f} \ar[d]_{\Phi_g,\Phi_f}^{\wr} \ar[rr]^{\Iso}_{C^{(!)}_{g,f}}
& & \shrbr{(fg)} \ar[d]^{\Phi_{fg}}_{\wr} \\
(\gs\fs)|_{\shrbr{\Y}} \ar[rr]^{\Iso}_{C^!_{g,f}} & & (fg)^!|_{\shrbr{\Y}}
}
\end{equation}
commutes.
\item[(d)] The map $\Phi_f$ is compatible with Zariski localizations
of $\Y$.
\end{itemize}
\end{thm}

\proof We wish to use \Tref{thm:gamma-ext}.
To that end,  for $\X\in \rbbF\cap\cbbF$ let 
$$
S_\X \colon \shrbr{\X}=\wDcp(\X) \to \X^!=\wDqcp(\X)
$$ 
be the natural inclusion. Here is the dictionary to help us pass from
\Ssref{ss:gamma-ext} to the situation we are now in.
The subcategory $\lbbG$ in \Ssref{ss:gamma-ext}
is, for us, $\rbbF\cap\cbbF$, the pseudofunctor $\natb{\boldsymbol{-}}$
is $\shrbrb{\boldsymbol{-}}$ and the maps $\gamma_f$ for $f\in {\mathbf{P}}
\cup{\mathbf{F}}$ are the maps $\Phi_f$ of \eqref{map:open(!)!} and
\eqref{iso:(!)!}. Checking diagram \eqref{diag:gamma-ext-1} commutes
amounts to checking that \eqref{diag:phi-ext} commutes when $f$ and $g$
are either both open or both pseudo-proper. 

In view of \eqref{diag:tr-tr-(!)}, the diagram \eqref{diag:phi-ext}
commutes whenever $f$ and $g$ are both pseudo-proper. If $f$ and $g$
are both open immersions then from the definition in \eqref{comp:open(!)!}
it is clear that \eqref{diag:phi-ext} commutes.

Next suppose $f\colon \X\to \Y$ is a pseudo-proper map in $\rbbF\cap\cbbF$
and suppose $u\colon \U\to \Y$ is an open immersion. Consider the
resulting fiber square diagram \eqref{diag:gamma-ext-2}. Pick a
residual complex $\cR$ on $\Y$. Then $\Tr{f}(\cR)$ is compatible
with open immersions into $\Y$ and hence so is $\ttr{f}^r$. It follows
that \eqref{diag:gamma-ext-3} commutes. 
\qed

\begin{rem}{\em As we mentioned earlier, the theorem is a way of
saying that the pseudofunctors $\shrbrb{\boldsymbol{-}}$ and
$\shr{\boldsymbol{-}}$ are isomorphic when each is appropriately
restricted to the ``domain" where both are meaningful. The reformulation
is a little awkward in view of the fact that there are
numerous ways of ``restricting" pseudofunctors since we are
dealing with families of categories indexed by yet another
category. We state what is needed briefly. Suppose 
${\boldsymbol{\vec{-}^{(!)}}}$ 
is the pseudofunctor obtained by restricting $\shrbrb{\boldsymbol{-}}$
to $\rbbF\cap\cbbF$ and ${\boldsymbol{\vec{-}^!}}$
the pseudofunctor obtained on $\rbbF\cap\cbbF$ via restricting
$\shr{\boldsymbol{-}}$ and by setting $\vec{\X}^!=\shrbr{\X}$
(there are two types of restrictions inherent here). Then \Tref{thm:(!)!}
asserts that ${\boldsymbol{\vec{-}^{(!)}}}$ is isomorphic to
${\boldsymbol{\vec{-}^!}}$.
}
\end{rem}

\begin{thm}\label{thm:(!)flat} Let $f\colon (\X,\,\De') \to (\Y,\,\De)$ 
be a map in $\rbbFc$, and $\De$ a codimensions function on $\Y$. 
Then $\shrbr{f}|_{\cm^*(\Y,\,\De)}$ takes values in $\cm^*(\X,\,\De')$
if and only if $f$ is flat.
\end{thm}

\proof The question is local and locally $f$ can be factored as
a composite of pseudo-proper maps and open immersions 
(cf. \cite[Lemma\,2.4.3]{lns}). Therefore, without loss of
generality, we assume that $f$ is also in $\cbbFc$.
If $f$ is flat then \Tref{thm:(!)!} and \Tref{thm:flat}
imply that $f^{(!)}\F$ is in $\cm^*(\X,\,\De')$ if 
$\F \in \cm^*(\Y,\,\De)$. For the converse we need to argue with
a little care, for \Tref{thm:flat} requires one to test
$f^!\F$ for every $\F\in \cm(\Y,\,\De)$, whereas we are restricting
our $\F$'s to lie in $\cm^*(\Y,\,\De)$. So suppose $f^{(!)}\F \in
\cm^*(\X,\,\De')$ whenever $\F\in\cm^*(\Y,\,\De)$. Let
$x\in X$, $S=\cO_{\X,x}$, $y=f(x)$ and $A=\cO_{\Y,y}$. By restricting
to an open neighborhood of $x$ if necessary, we may assume that
$f=gh$ where $\X\xrightarrow{h}{\mathcal{P}}$ is a closed immersion
and ${\mathcal{P}}\xrightarrow{g}\Y$ is a smooth map in $\cbbF$.
Let $R=\cO_{{\mathcal{P}},h(x)}$ and $z=h(x)$. As in \Pref{prop:tor-ext}, let
$k=k_A$, $K$ an $R$--injective hull of $k_R$ and $\varphi\colon
{\widehat A}\to {\widehat R}$ the natural map induced by $g$. 
In what follows we will be applying \Tref{thm:(!)!} to the
maps $f$, $g$ and $h$ without comment since all three of them
are morphisms in $\rbbF\cap \cbbF$.
Let $\I\subset\cO_\Y$
be the open prime ideal sheaf corresponding to the point $y$.
Let $\cR$ be a residual complex in $\Cozs{\De}(\Y)$. Set
$$
\F\set \sHomb_\Y(\cO_\Y/\I,\,\cR).
$$
Then $Q_\Y\F\in \cm^*(\Y,\,\De)$. Hence by our hypothesis
$\fs\F\in \cm^*(\X,\,\De')$. Therefore
$$
\Hr_x^{i+\De'(x)}(\fs\F) = 0
$$
for $i\ne 0$. By \eqref{iso:ext} and the fact that
$\fs\F \simeq h^!\sh{g}\F$ this amounts to saying
$$
{\mathrm{Ext}}^i_R(S,\,(\sh{g}\F)(z)) = 0
$$ 
for $i>0$. By \Pref{prop:tor-ext}(b) this means that
${\mathrm{Tor}}_i^A(S,\,k)=0$ for $i>0$. Hence we are done by
\cite[p.\,174,\,Thm.\,22.3,\,(i) and (iii)]{matsumura}.
\qed

The following result follows from \Tref{thm:injective}
and \Tref{thm:(!)!}, and the fact that we may, in each
case, Zariski localize the source and the target.

\begin{thm}\label{thm:(!)injective} Let $(\Y,\,\De)\in \rbbFc$ and let
$\F\in \Cozs{\De}(\Y)$. Then the following are equivalent
\begin{itemize}
\item[(a)] $\F$ is a complex of $\Aqct$--injectives;
\item[(b)] $f^{(!)}\F\in \cm^*(\X,\,\De')$ for every morphism $(\X,\,\De')
\xrightarrow{f} (\Y,\,\De)$ in $\rbbFc$;
\item[(c)] $f^{(!)}\F\in \cm^*(X,\,\De')$ for every closed immersion
$(X,\,\De')\xrightarrow{f} (\Y,\,\De)$ such that $X$ is an
integral ordinary scheme.
\end{itemize}
\end{thm}


\begin{ack} This work has been a long time in the making, and
from its very beginning in 1995, Joe Lipman has been a constant
source of stimulation and encouragement. I also wish to thank
Suresh Nayak for stimulating conversations and for being
sensitive to the needs of this paper while Lipman, he and I
were working on \cite{lns}.
\end{ack}

\bibliographystyle{plain}
\end{document}